
\input graphicx.tex
\input amssym.def
\input amssym
\magnification=1000
\baselineskip = 0.20truein
\lineskiplimit = 0.01truein
\lineskip = 0.01truein
\vsize = 8.9truein
\voffset = 0.1truein
\parskip = 0.09truein
\parindent = 0.3truein
\settabs 12 \columns
\hsize = 6.0truein
\hoffset = 0.1truein

\setbox\strutbox=\hbox{%
\vrule height .708\baselineskip
depth .292\baselineskip
width 0pt}
\font\caps=cmcsc10
\font\bigtenrm=cmr10 at 14pt

\def\sqr#1#2{{\vcenter{\vbox{\hrule height.#2pt
\hbox{\vrule width.#2pt height#1pt \kern#1pt
\vrule width.#2pt}
\hrule height.#2pt}}}}
\def\square{\mathchoice\sqr46\sqr46\sqr{3.1}6\sqr{2.3}4}

\centerline{\bigtenrm A POLYNOMIAL UPPER BOUND}
\centerline{\bigtenrm ON REIDEMEISTER MOVES}
\tenrm
\vskip 14pt
\centerline{MARC LACKENBY}
\vskip 18pt
\centerline{\caps Abstract}
{\narrower \noindent
We prove that any diagram of the unknot with $c$ crossings may be reduced to the trivial
diagram using at most $(236 \,c)^{11}$ Reidemeister moves.\par}
\vskip 18pt
\centerline{\caps 1. Introduction}
\vskip 6pt

Turing stated in one of his famous articles [24] that `{\sl No systematic method is
yet known by which one can tell whether two knots are the same.}'
Even the basic case of recognising the unknot is not obviously soluble.
A few years later, in his groundbreaking work on normal surfaces, Haken solved the problem
of recognising the unknot [11] and then made a crucial contribution to the more general problem of whether
two knots are equivalent [12]. This was finally solved by the efforts of several mathematicians,
including Hemion [17] and Matveev [22]. But it remains a major unresolved question
to determine exactly how complex these problems are. The current state of
our knowledge is that unknot recognition is in NP and co-NP.
The fact that it is in NP is due to Hass, Lagarias and Pippenger [14] and that it is in co-NP
was proved by Agol [1], but not written down in detail, and an alternative
solution was given by Kuperberg [21], assuming the Generalised Riemann
Hypothesis.

There are many examples of challenging diagrams of unknots.
In 1934, Goeritz gave an example of a diagram with 11 crossings,
with the property that any sequence of Reidemeister moves taking it
to the trivial diagram must go via a diagram with more than 11 crossings.
Other tricky examples have been given by Thistlethwaite, Haken,
Henrich and Kauffman [18]. We include some of these below in Figures 1-3.
They all point to the probable conclusion that there is no simple
way of recognising the unknot.

The most elementary and natural way of approaching the unknot recognition problem is to try to
find an explicit upper bound on the number of Reidemeister moves
required to turn a given diagram of the unknot with $c$ crossings
into the trivial diagram. It is easy to see that the existence of a computable upper bound
is equivalent to the solvability of the unknot recognition problem.
But of course one wants a bound that is as small a function of $c$
as possible. 

In [13], Hass and Lagarias showed that a diagram of the unknot with
$c$ crossings can be converted into the trivial diagram using
at most $2^{kc}$ Reidemeister moves, where $k = 10^{11}$. In [15], Hass and Nowik proved
that, in general, at least $c^2 / 25$ moves are required. There
is a large gap between these upper and lower bounds, and so it has
remained a basic question: is there a polynomial upper bound on the
number of Reidemeister moves required to turn an unknot diagram
into the trivial diagram? This is what we solve in this paper.

\noindent {\bf Theorem 1.1.} {\sl Let $D$ be a diagram of the unknot
with $c$ crossings. Then there is a sequence of at most $(236 \, c)^{11}$
Reidemeister moves that transforms $D$ into the trivial diagram.
Moreover, every diagram in this sequence has at most $(7 \, c)^2$ crossings.
}

It is worth pointing out that this does not actually improve our knowledge of the complexity class
of the unknot recognition problem. But it does give
an alternative way of establishing that the unknot recognition problem is in NP,
because the sequence
of Reidemeister moves provided by the above theorem gives
a polynomial time certificate of unknottedness. Therefore,
it remains an unsolved problem whether unknot recognition is in P.
Of course, this may be very difficult, because a negative answer
would imply that ${\rm P} \not= {\rm NP}$. Moreover, it is unlikely that a polynomial time algorithm
could be ruled out, even conditional upon the hypothesis
that ${\rm P} \not= {\rm NP}$, because it is widely
conjectured that problems in NP $\cap$ co-NP are not
${\rm NP}$-complete (see p.95 of [10]).

We also have a result for split links.

\noindent {\bf Theorem 1.2.} {\sl Let $D$ be a diagram of a split
link with $c$ crossings. Then there is a sequence of at most $(49 \, c)^{11}$
Reidemeister moves that transforms $D$ into a disconnected diagram.
Moreover, every diagram in this sequence has at most $9 \, c^2$ crossings.
}

Our theorems rely in a crucial way on groundbreaking work of
Dynnikov [8]. He considered a special way of arranging a knot or link called
an arc presentation. One way of visualising is these is via rectangular diagrams (also called grid diagrams),
which are diagrams in the plane consisting of horizontal and vertical
arcs, subject to the condition that the vertical arc always passes
over the horizontal one at a crossing and the condition that no two arcs are
collinear. The number of vertical arcs
equals the number of horizontal arcs, and this is known as the {\sl arc
index} of this presentation. Dynnikov proved the surprising result
that any arc presentation of the unknot can be reduced to the trivial
presentation using a sequence of moves, known as exchange moves, cyclic
permutations and destabilisations (see Figures 5-7). Crucially, the arc index never needs to increase.
This has the striking consequence that if a diagram of the unknot
has $c$ crossings, then there is a sequence of Reidemeister moves taking
it to the trivial diagram, such that all diagrams in this sequence
have at most $2(c+1)^2$ crossings (Theorem 2 in [8]). But this does not give a polynomial
upper bound on the number of such moves.

It is also possible to show that the approach of Hass and Lagarias in [13]
does not provide a polynomial upper bound. They start with a diagram
of the unknot with $c$ crossings, and they use this to build a
triangulation of a convex polyhedron with $t \leq 840 c$ tetrahedra,
each of which is straight in ${\Bbb R}^3$ and which contains the given
unknot in its 1-skeleton. From this, they construct a triangulation
of the knot exterior. By work of Haken [11], the disc that the unknot spans can be
realised as a normal surface with respect to this triangulation, and Hass and Lagarias show that
at most $2^{kt}$ normal triangles and squares are required, where $k = 10^7$. 
They then isotope the unknot across this disc. The projection to the plane
of the diagram then gives a sequence of Reidemeister moves.
The bound on the number of normal squares and triangles gives the
exponential bound on the number of Reidemeister moves. It does not seem
feasible to use this approach
of sliding the knot across a normal spanning disc to obtain a better bound
on Reidemeister moves. This is because Hass,
Snoeyink and Thurston [16] gave examples of unknots consisting of
$10n + 9$ straight arcs, for which any piecewise linear spanning disc must
have at least $2^{n-1}$ triangular faces.

Instead, our approach here is to combine Dynnikov's methods with the use of
normal surfaces. Given an arc presentation for an unknot, Dynnikov
explains how a spanning disc may be placed in what he calls
admissible form. He defines a measure of complexity on such
surfaces. The key part of his argument is to show that an
admissible spanning disc must have at some point a certain
local configuration. This then specifies a way of modifying
the surface and the arc presentation. This has the effect of
performing `generalised exchange moves' on the arc presentation and possibly
destabilisations. He shows that, during this process, either
a destabilisation is performed or the complexity of the spanning
disc has gone down. 

Dynnikov defines a triangulation of the 3-sphere associated to
an arc presentation of a link. If the arc index is $n$, this has
$n^2$ tetrahedra. It turns out that placing the spanning disc or splitting sphere in
admissible form is almost equivalent to placing this surface into
normal form with respect to this triangulation. Moreover, his
measure of complexity is (under reasonable assumptions) just the number of intersections between
the disc or sphere and certain edges of the triangulation. Thus, using
the bound on the complexity of normal surfaces that was proved
by Hass and Lagarias in [13], the complexity of a splitting sphere is at most $n 2^{7n^2}$.
(A similar, but slightly larger bound is required for the spanning disc
of the unknot.) Hence, using Dynnikov's argument, one can show that 
the number of generalised exchange moves that
one needs to perform before one can apply a destabilisation is
at most an exponential function of $n^2$. 

However, this is much larger than a polynomial upper bound.
To obtain this, one needs to go deeper into
normal surface theory. In a triangulated 3-manifold with $n^2$ tetrahedra,
any normal surface consists of at most $5n^2$ types of normal triangles and squares. One can
show that if there is a local configuration of the spanning disc or splitting sphere
which specifies a way of reducing complexity, then one can also reduce
complexity in regions of the surface that are normally parallel.
Thus, one might hope that, using a single generalised exchange
move, one can reduce complexity by a factor of roughly $(1 - n^{-2})$. This
is probably too optimistic, for it may be the case that most of the weight 
of the normal surface is concentrated in regions
where this good configuration does not occur. The key technical
part of this paper is to show that, under this situation,
the surface does not have minimal complexity. In particular,
there is another spanning disc or splitting sphere, with smaller complexity, for the
same arc presentation. This is shown by establishing that some multiple of the
given surface is actually a normal sum of a normal torus
and a multiple of some simpler spanning disc or splitting sphere. The proof of this is somewhat delicate,
and relies on the use of branched surfaces and `first-return maps'.

Thus, the results that we actually prove are as follows. (For the definitions of
trivial and disconnected arc presentations, see Section 2.1.)

\noindent {\bf Theorem 1.3.} {\sl Let $D$ be an arc presentation of
the unknot with arc index $n$. Suppose that the associated rectangular diagram
has writhe $k$. Then there is a sequence of at most
$4 \times 10^{18} \, n^{10}$ exchange moves, at most $6 \times 10^{18} \, n^{9}$ cyclic
permutations, at most $10^{19} \, n^8$ generalised exchange moves, at most
$3 \times 10^{13} \, n^6$ stabilisations and at most $3 \times 10^{13} n^6$
destabilisations taking $D$ to the trivial arc presentation. Moreover, the arc
index is at most $2n+|k|+1$ throughout this sequence of moves.}

\noindent {\bf Theorem 1.4.} {\sl Let $D$ be an arc presentation of
a split link with arc index $n$. Then there is a sequence of at most $3 \times 10^{11} \, n^8$
generalised exchange moves, at most $2 \times 10^{11} \, n^9$ cyclic permutations
and at most $8 \times 10^{10} \, n^{10}$ exchange moves 
that takes $D$ to a disconnected arc presentation.}

Now each generalised exchange move on an arc presentation with arc index $n$
can be expressed as a composition of
at most $(3/2)n^3$ Reidemeister moves (Lemma 2.4). Any exchange move is a product
of at most $n$ Reidemeister moves (Lemma 2.2). A cyclic permutation requires at most
$(n-1)^2$ Reidemeister moves (Lemma 2.3). Also, given any diagram of a knot or link with $c$
crossings, this is isotopic to a rectangular diagram with arc index at most $(81/20)c$ (Lemma 2.1).
Any rectangular diagram with arc index $n$ has at most $(n-1)^2/2$ crossings (see the proof of
Theorem 2 in [8]). These observations, combined with Theorems 1.3 and 1.4, imply Theorems 1.1 and 1.2.

The plan of the paper is as follows. In Section 2, we give some elementary
properties of arc presentations. Section 3 contains an overview of Dynnikov's
proof that arc presentations of the unknot and split links can be
simplified using a sequence of exchange moves, cyclic permutations and destabilisations.
In Section 4, we present an alternative argument, which provides an explicit
upper bound on the number of exchange moves, cyclic permutations, stabilisations and destabilisations
required to trivialise a rectangular diagram of the unknot, given an upper
bound for the complexity of the spanning disc. This is an unsurprising
result, and is required only in the case of the unknot.
In Section 5, we recall some key facts from normal surface theory,
including some results about vertex normal surfaces. We introduce a new
notion of a boundary-vertex normal surface, which is useful in the
parts of the proof dealing with the unknot. In Section 6, we introduce
normal surface theory to arc presentations. We give Dynnikov's triangulation
of the 3-sphere, and explain how surfaces that are normal with respect
to this triangulation have a form that is very close to admissible. 
Section 7 contains the proof of Theorems 1.3 and 1.4, 
assuming the result that the normal spanning disc or splitting sphere
cannot contain large `Euclidean' regions. This is proved in Sections 8
and 9, using branched surfaces. In the final section, we discuss possible improvements
to the degree of the polynomial bound, and we also give some potential directions
for further research.

The presence of surfaces with boundary causes several complications in
these arguments, and so the case of the unknot is more complex than
the case of split links. We therefore suggest that the reader initially concentrates
on the split link case.

I would like to thank the referee for their very careful reading of an earlier version of this
paper.

\vskip 18pt
\centerline{
\includegraphics[width=1.55in]{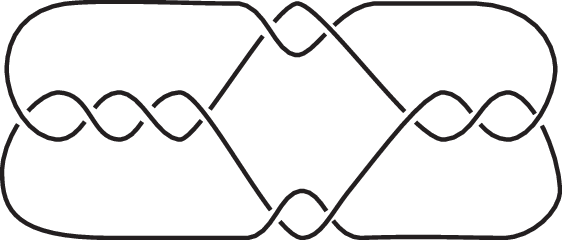}
}
\vskip 6pt
\centerline{Figure 1: Goeritz's unknot}

\vskip 12pt
\centerline{
\includegraphics[width=1.6in]{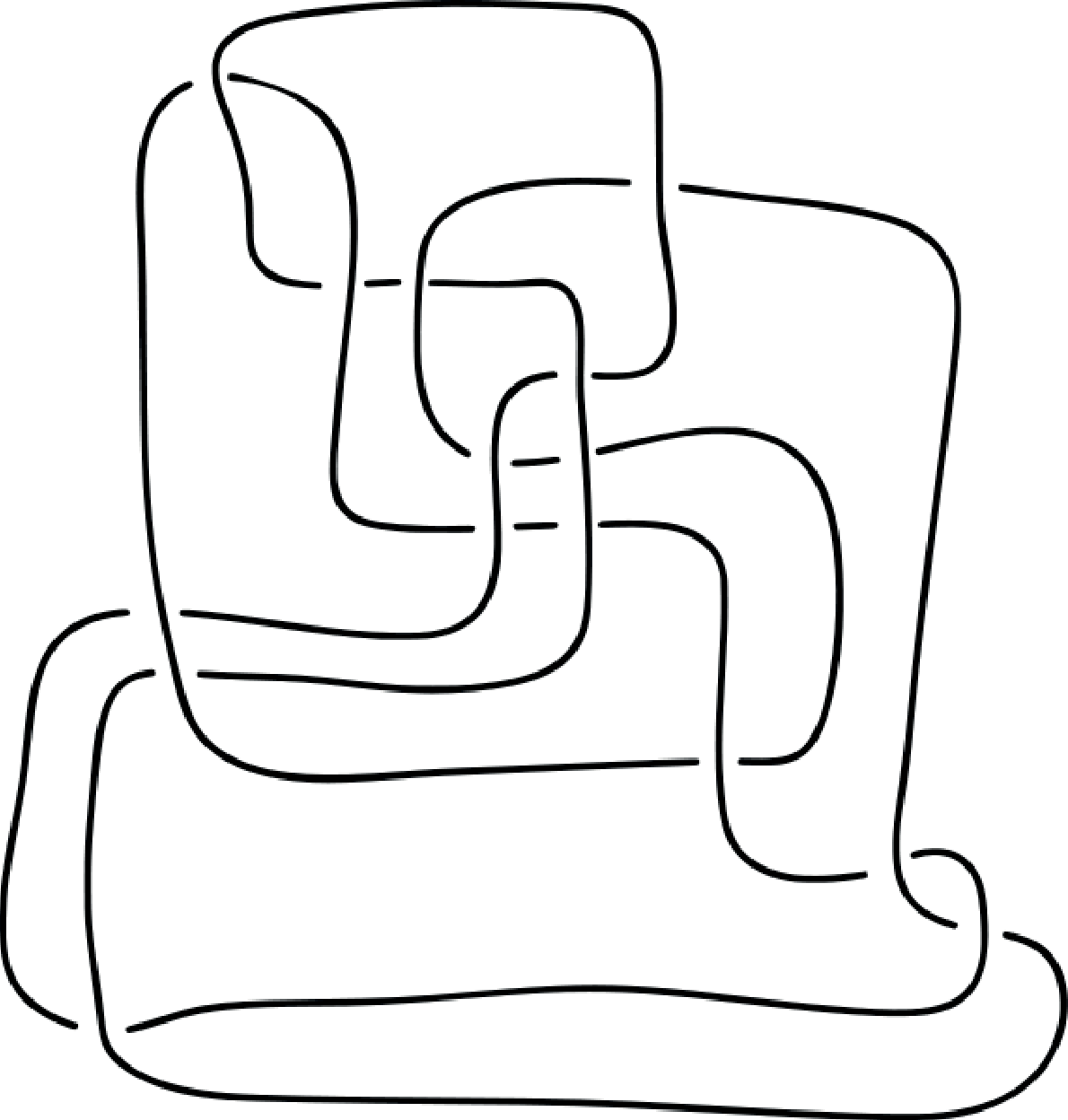}
}
\vskip 6pt
\centerline{Figure 2: Thistlethwaite's unknot}

\vskip 12pt
\centerline{
\includegraphics[width=2.1in]{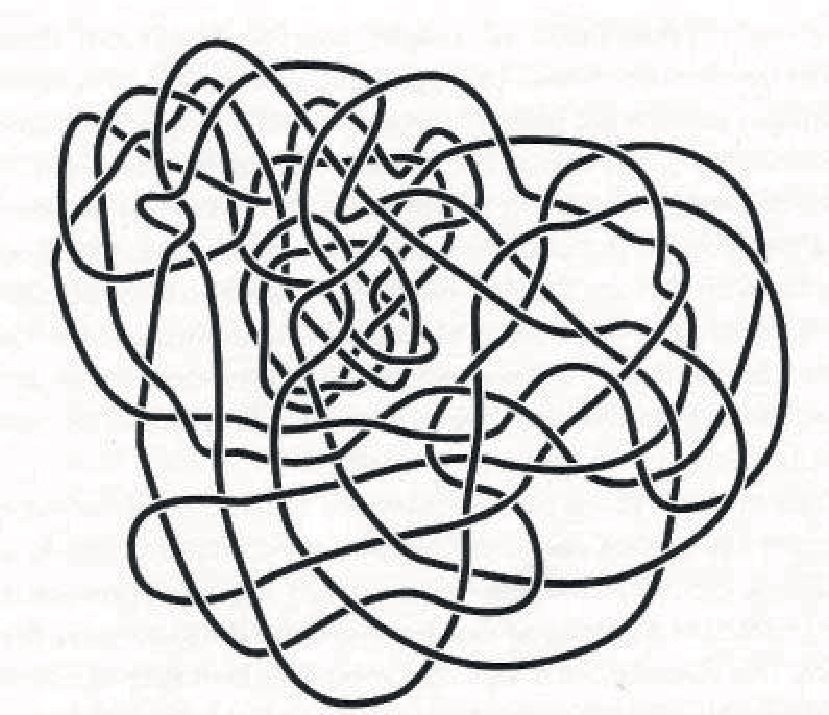}
}
\vskip 5pt
\centerline{Figure 3: One of Haken's unknots (image courtesy of Cameron Gordon)}

\vskip 18pt
\centerline{\caps 2. Basic properties of arc presentations}
\vskip 6pt

In this section, we present some elementary material on arc presentations
and rectangular diagrams. Much of this was first discovered by Cromwell
[6]. We have largely followed Dynnikov's presentation in [8].

\vskip 6pt
\noindent {\caps 2.1. Definition of arc presentations}
\vskip 6pt

We fix a description of the 3-sphere as the join $S^1 \ast S^1$
of two circles. The co-ordinate system $(\phi, \tau, \theta)$ is used,
where $\phi, \theta \in {\Bbb R} / 2 \pi {\Bbb Z}$ are co-ordinates on the circles,
and $\tau \in [0,1]$. Thus, $(\phi, 0, \theta_1)$ and $(\phi, 0, \theta_2)$
are identified for all $\theta_1$ and $\theta_2$. Similarly, $(\phi_1, 1, \theta)$
and $(\phi_2, 1, \theta)$ are identified for all $\phi_1$ and $\phi_2$.
The circles $\tau = 0$ and $1$ are denoted by $S^1_\phi$ and $S^1_\theta$
respectively. The circle $S^1_\phi$ is called the {\sl binding circle}.
The open disc defined by $\theta = t$ and $\tau > 0$ is called a {\sl page}
and denoted ${\cal D}_t$.

Suppose that a link $L$ satisfies the following two conditions: $L \cap S^1_\phi$
is a finite set, called the {\sl vertices} of $L$, and for any $t \in {\Bbb R} / 2 \pi {\Bbb Z}$, the intersection
${\cal D}_t \cap L$ is either empty or an open arc approaching two distinct
vertices. This is called an {\sl arc presentation} of $L$. The number of vertices
equals the number of pages that contain open arcs of $L$. This number is
called the {\sl arc index} of the arc presentation.

We say that an arc presentation is {\sl trivial} if it has arc index 2.
We say that it is {\sl disconnected} if there is a 2-sphere that intersects each page
in a single embedded arc, and which has components of $L$ on both sides of it.

\vskip 6pt
\noindent {\caps 2.2. Rectangular diagrams}
\vskip 6pt

There is an equivalence between arc presentations and rectangular
diagrams, which we now describe. 

A {\sl rectangular diagram} of a link $L$ is a link diagram defined as follows.
The plane of the diagram has a product structure ${\Bbb R} \times {\Bbb R}$. 
We require that the projection of
$L$ is a union of arcs, each of which is of the form $\{ s \} \times [t_1, t_2]$
or $[s_1, s_2] \times \{ t \}$. These are known as {\sl vertical} and
{\sl horizontal arcs}. Whenever the interiors of two arcs of 
the projection intersect, the over-arc at the resulting crossing is 
required to be the vertical arc. Also, no two arcs may be collinear.

\vskip 18pt
\centerline{
\includegraphics[width=1.8in]{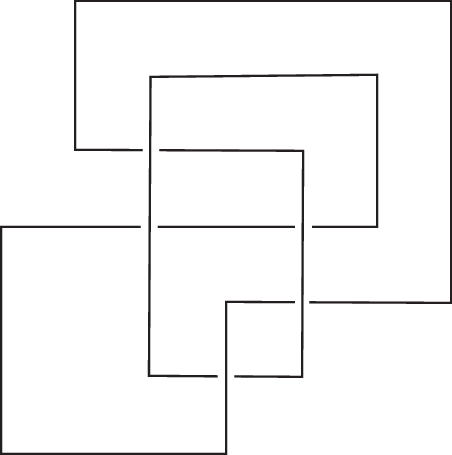}
}
\vskip 6pt
\centerline{Figure 4: A rectangular diagram}

An arc presentation of $L$ determines a rectangular diagram as follows.
The arc presentation can be specified by the following data: the $\phi$-values
of the vertices, the values of $\theta$ where the page contains an arc, and
the vertices at the endpoint of this arc. Let $s_1, \dots, s_n \in [0, 2 \pi)$
denote the $\phi$-values of the vertices, and let $t_1, \dots, t_n \in [0, 2 \pi)$
denote the $\theta$-values of the arcs. For each arc of $L$, lying in ${\cal D}_t$,
joining vertices $s_i$ and $s_j$ where $s_i < s_j$, we insert a horizontal edge of the
rectangular diagram at $[s_i, s_j] \times \{ t \}$. For each vertex $s$ of $L$,
its two adjacent arcs lie in ${\cal D}_{t_i}$ and ${\cal D}_{t_j}$, where $t_i < t_j$. For each
such vertex, we insert a vertical edge of the rectangular diagram at
$\{ s \} \times [t_i, t_j]$. 

We now explain briefly why this is indeed a diagram of $L$. In fact, we will give
a reasonably explicit map from the complement of the link defined by the arc presentation
to the complement of the link defined by the rectangular diagram. (A more complete
explanation is given in [6].)

Consider an arc presentation for $L$. We replace each arc of $L$ in a page ${\cal D}_t$,
joining vertices $s_1$ and $s_2$, where $s_1 < s_2$, by the concatenation of three arcs:
$$\eqalign{
& \{ \phi = s_1, \theta = t, \epsilon \leq \tau \leq 1- \epsilon \} \cr
\cup \ &\{ s_1 \leq \phi \leq s_2, \theta = t, \tau = 1- \epsilon \} \cr
\cup \ &\{ \phi = s_2, \theta = t, \epsilon \leq \tau \leq 1- \epsilon \}.}$$
Here, $\epsilon$ is some fixed real number in the interval $(0, 1/2)$. 
As $L$ approaches a vertex $s$ in pages ${\cal D}_{t_1}$ and ${\cal D}_{t_2}$,
where $t_1 < t_2$, we replace it by an arc
$$\{ \phi = s, t_1 \leq \theta \leq t_2, \tau = \epsilon \}.$$
After this, $L$ lies in the region $\{ \epsilon \leq \tau \leq 1- \epsilon \}$,
which is a thickened torus. If we project onto $\{ \tau = 1/2 \}$,
we obtain a diagram in a torus, and this torus is standardly embedded
in $S^3$. Because we ensured that the arcs did not go beyond $\phi = 0$ and $\theta = 0$,
the diagram lies in the square
$$\{ 0 \leq \phi < 2 \pi, 0 \leq \theta < 2 \pi, \tau = 1/2 \}.$$
If we realise this square as a subset of the plane, we obtain the required rectangular
diagram for $L$.

\vskip 6pt
\noindent {\caps 2.3. From ordinary diagrams to rectangular diagrams}
\vskip 6pt

Cromwell [6] proved that any link $L$ has an arc presentation, by starting with
an arbitrary diagram of $L$ and making it rectangular. In this subsection, we 
will carry out this procedure, but also keep track of an upper bound on the arc index of the
resulting rectangular diagram.

\noindent {\bf Lemma 2.1.} {\sl Let $D$ be a diagram of a link with $c$ crossings.
Then $D$ is isotopic to a rectangular diagram with arc index at most $(81/20) c$.}

\noindent {\sl Proof.} We may clearly assume that $D$ is connected. We may also
assume that $D$ contains no edge loops (which are arcs of the diagram with both
endpoints at the same crossing). For we may remove all such edge loops,
then isotope the resulting diagram so that it is rectangular, and then add back in
the loops in a rectangular fashion.

Let $X$ be the underlying 4-valent planar graph
specified by $D$. This has $2c$ edges. We will modify $X$ by subdividing
its edges. If any pairs of edges are parallel, subdivide
one of the edges from each pair. We may assume that at least 6 edges of the diagram are not parallel 
to any other edge, since otherwise $D$ is a standard diagram of a $(2,n)$-torus link or a
simple type of 2-bridge link, in which
case the lemma is easy to establish.
We deduce that $X$ now has at most $2c + (2c - 6)/2 = 3c -3$ edges. 

In [23], Storer examined the problem of how to arrange a planar graph (with no edge loops or
parallel edges) so that
its edges are horizontal and vertical arcs, possibly after subdividing 
its edges. By Corollary 4 in [23], $X$ may be subdivided so that it has
a total of at most $(17/10)m + 4$ vertices, where $m$ is the original number of edges of $X$,
and then isotoped so that each edge is horizontal or vertical in the plane.
So, the number of 2-valent vertices of $X$ is now
at most $(17/10)(3c - 3) + 4 - c \leq (41/10)c$. 

This diagram might not be a rectangular diagram for two reasons. Firstly,
some edges may be collinear. But if so, then a small modification,
keeping the arcs horizontal and vertical, can made to avoid this.
Secondly, at some crossings, the over-arc may be horizontal,
rather than vertical. But if so, there is an obvious modification which
introduces 8 new 2-valent vertices at such a crossing (see Figure 7 of [6]).
Note that we may assume that at least half the crossings have the correct
behaviour, as otherwise, we can instead just rotate the entire diagram by a quarter turn.
So, the number of 2-valent vertices is at most $(41/10)c + 4c = (81/10)c$.
The arc index of this rectangular diagram is at most half the number of
2-valent vertices, which is less than $(81/20)c$, as required. $\square$

This bound of $(81/20)c$ is obviously not optimal. In fact, Cromwell and Nutt in [7]
show that in many cases, the link specified by $D$ has an arc presentation with arc index
at most $c+ 2$. In the proof of Theorem 2 in [8], Dynnikov states that one can always find an arc
presentation for the link with arc index at most $2c + 2$. However, the resulting
rectangular diagram is not necessarily isotopic to $D$. So, to be able to
use this fact in the proof of Theorems 1.1 and 1.2, one would need to be able to find an upper bound on the
number of Reidemeister moves required to transform $D$ into the new
rectangular diagram. This is surely possible, but it is not completely
straightforward. So, we have chosen to follow the simpler course
of isotoping $D$ so that it is rectangular, even though this might not
lead to the optimal upper bound on arc index.

\vskip 6pt
\noindent {\caps 2.4. Exchange moves, stabilisations and destabilisations}
\vskip 6pt

Cromwell [6] introduced a set of moves, which modify an arc presentation
without changing the link. These are most simply visualised using
rectangular diagrams:
\item{(1)} cyclic permutation of the horizontal (or vertical) arcs;
\item{(2)} stabilisation and destabilisation;
\item{(3)} interchanging parallel edges of the rectangular diagram, as long as 
they have no edges between them, and their pairs of endpoints do not
interleave; this is termed an {\sl exchange move}.

\noindent These are shown in Figures 5-7.

When we use the term exchange move, we assume that the parallel edges
that are moved past each other do not lie either side of $\theta = 0$
or $\phi = 0$. In this case, a cyclic permutation needs to be done first,
before the exchange move can be performed. The reason that we make
this distinction is that an exchange move requires fewer Reidemeister moves
in general than a cyclic permutation. 

\vskip 18pt
\centerline{
\includegraphics[width=2.7in]{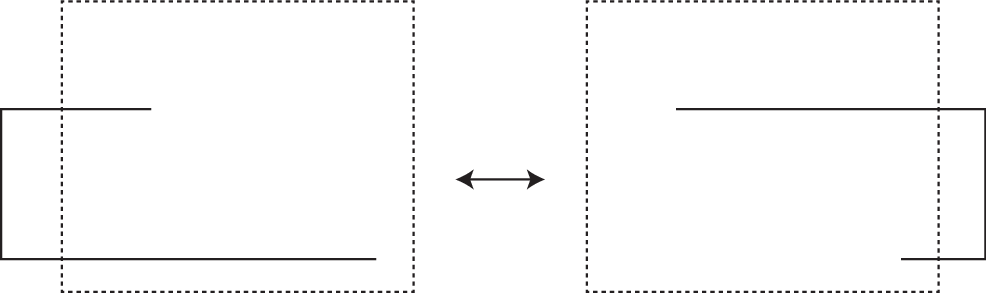}
}
\vskip 6pt
\centerline{Figure 5: Cyclic permutation of the vertical edges}

\vskip 18pt
\centerline{
\includegraphics[width=3.2in]{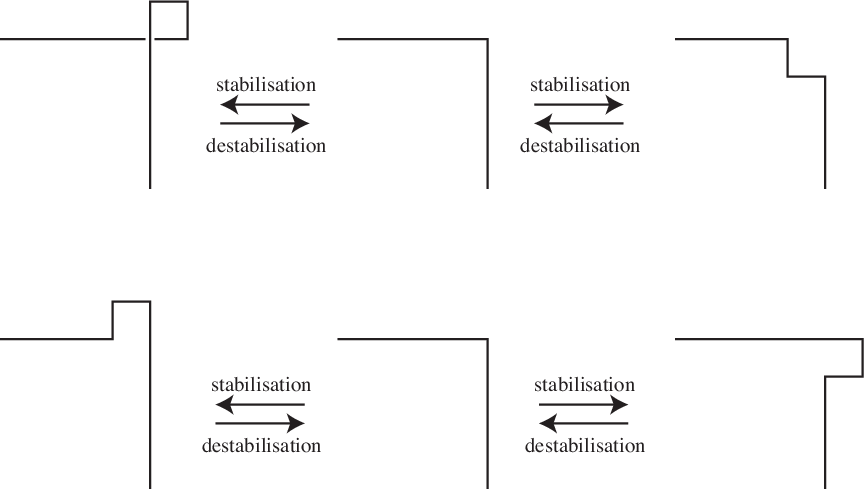}
}
\vskip 6pt
\centerline{Figure 6: Stabilisations and destabilisations}

\vskip 18pt
\centerline{
\includegraphics[width=3.2in]{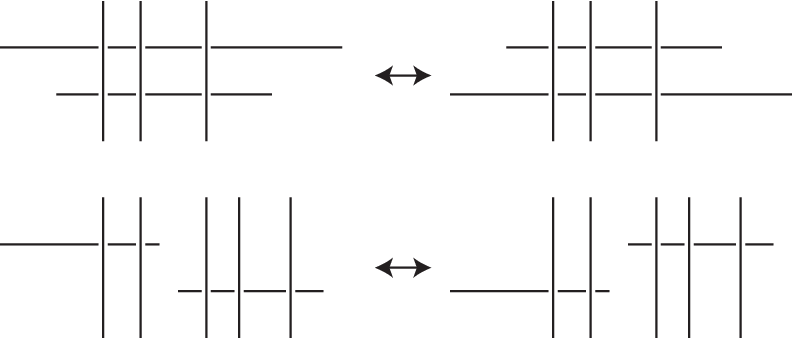}
}
\vskip 6pt
\centerline{Figure 7: Exchange moves}

We now provide upper bounds on the number
of these moves.

\noindent {\bf Lemma 2.2.} {\sl Let $n$ be the arc index
of an arc presentation of $L$, and let $D$ be the resulting
rectangular diagram. Suppose that an exchange move
is performed on this arc presentation, and let $D'$
be the resulting rectangular diagram. Then $D'$ and $D$
differ by a sequence of at most $n$ Reidemeister moves.}

\noindent {\sl Proof.} This is fairly evident from Figure 7.
In the bottom case of Figure 7, no Reidemeister moves
are required. In the top case, one might first need to make a type 2 Reidemeister move
to make the two horizontal edges overlap, then a sequence of
at most $n-2$ type 3 Reidemeister moves, then possibly a type 2 move.
$\square$

\noindent {\bf Lemma 2.3.} {\sl Let $n$ be the arc index
of an arc presentation of $L$. Suppose that a cyclic permutation
is performed on the vertical (or horizontal) arcs. Then the
resulting rectangular diagrams differ by a sequence of at most $(n-1)^2$
Reidemeister moves.}

\noindent {\sl Proof.} In Figure 5, a vertical arc is slid across the diagram
from left to right. As it meets another vertical arc, a type 2 Reidemeister move
might need to be performed, followed by a sequence of at most $(n-2)$
type 3 moves, then possibly a type 2 move if one was not performed at the
beginning. This is at most $n-1$ Reidemeister moves. There are at
most $n-1$ vertical arcs that it is slid across. So, at most
$(n-1)^2$ Reidemeister moves are needed in total. $\square$

\vskip 6pt
\noindent {\caps 2.5. Generalised exchange moves}
\vskip 6pt

A more substantial modification to an arc presentation was introduced in [6],
known as a generalised exchange move. This is defined as follows.

Let $0 < s_1 < s_2 < s_3 < 2 \pi$ be values of $\phi$ which are disjoint from the vertices
of $L$. Let $0 \leq t_1 < t_2 < 2 \pi$ be values of $\theta$ which are disjoint from the
arcs of $L$. Suppose that each horizontal arc $[s,s'] \times \{ t \}$
of the rectangular diagram satisfies the following conditions:
\item{(1)} if $t \in (t_1, t_2)$, then $\{ s , s' \}$ is not interleaved
with $\{ s_2, s_3\}$;
\item{(2)} if $t \in S^1_\theta - (t_1, t_2)$, then $\{ s , s' \}$ is not interleaved
with $\{ s_1, s_2\}$.

\noindent Then one can modify the rectangular diagram by changing
the $\phi$ value of all the vertices between $s_1$ and $s_2$
so that they lie between $s_2$ and $s_3$ in the same order,
and by changing the $\phi$ value of all the vertices between $s_2$ and $s_3$
so that they lie between $s_1$ and $s_2$ in the same order.
This is a {\sl generalised exchange move}.

The effect of a generalised exchange move on the rectangular diagram is shown in
Figure 8, in the case where $t_1 = 0$, where it is evident that it does not change the link type.

\noindent {\bf Lemma 2.4.} {\sl Let $n$ be the arc index of an arc presentation
of $L$. A generalised exchange move on this arc presentation is a composition
of at most $(3/2)n^3$ Reidemeister moves. It is also a composition
of at most $n$ cyclic permutations and at most $(3/4) n^2$ exchange moves.}

\noindent {\sl Proof.} In Figure 8, a generalised exchange move is shown
where $t_1 = 0$. In general, as many as $n/2$ cyclic permutations
may need to be made before $t_1 = 0$ and by Lemma 2.3, these
may require at most $n^3/2$ Reidemeister moves.

Figure 8 shows how the generalised
exchange moves can be divided into three steps. We estimate
the number of Reidemeister moves or exchange moves
required in the first step. Place each horizontal arc $[s, s'] \times \{ t \}$
in one of the following sets:
\item{(1)} In $A_1$ if $s, s' \in (s_1, s_2)$ and $t \in S^1_\theta - (t_1, t_2)$;
\item{(2)} In $A_2$ if $s, s' \not \in (s_1, s_2)$ and $t \in S^1_\theta - (t_1, t_2)$;
\item{(3)} In $A_3$ if $s, s' \not \in (s_2, s_3)$ and $t \in (t_1, t_2)$;
\item{(4)} In $A_4$ if $s, s' \in (s_2, s_3)$ and $t \in (t_1, t_2)$.

\noindent So, the first step of the generalised exchange move slides the $A_1$ arcs
past those in $A_4$ and some of those in $A_2$. It also slides the
$A_4$ arcs past some of those in $A_3$. The number of exchange moves is therefore
at most $|A_1||A_4| + |A_1||A_2| + |A_3||A_4| \leq (|A_1| + |A_3|) (|A_2| + |A_4|)
\leq n^2 /4$. The other two
steps are similar, and so we obtain the required bound of
$(3/4) n^2$ exchange moves. 
By Lemma 2.2, the first and third steps
each require at most $n^3/4$ Reidemeister moves. The second step
evidently needs no Reidemeister moves. 

Finally, we reverse the cyclic permutations that were made initially. This is necessary because
the generalised exchange move does not change the $\theta$-value of any arc. Again, by Lemma 2.3,
these require at most $n^3/2$ Reidemeister moves. So, in total, at most
$(3/2)n^3$ Reidemeister moves are needed.  $\square$

\vskip 18pt
\centerline{
\includegraphics[width=3in]{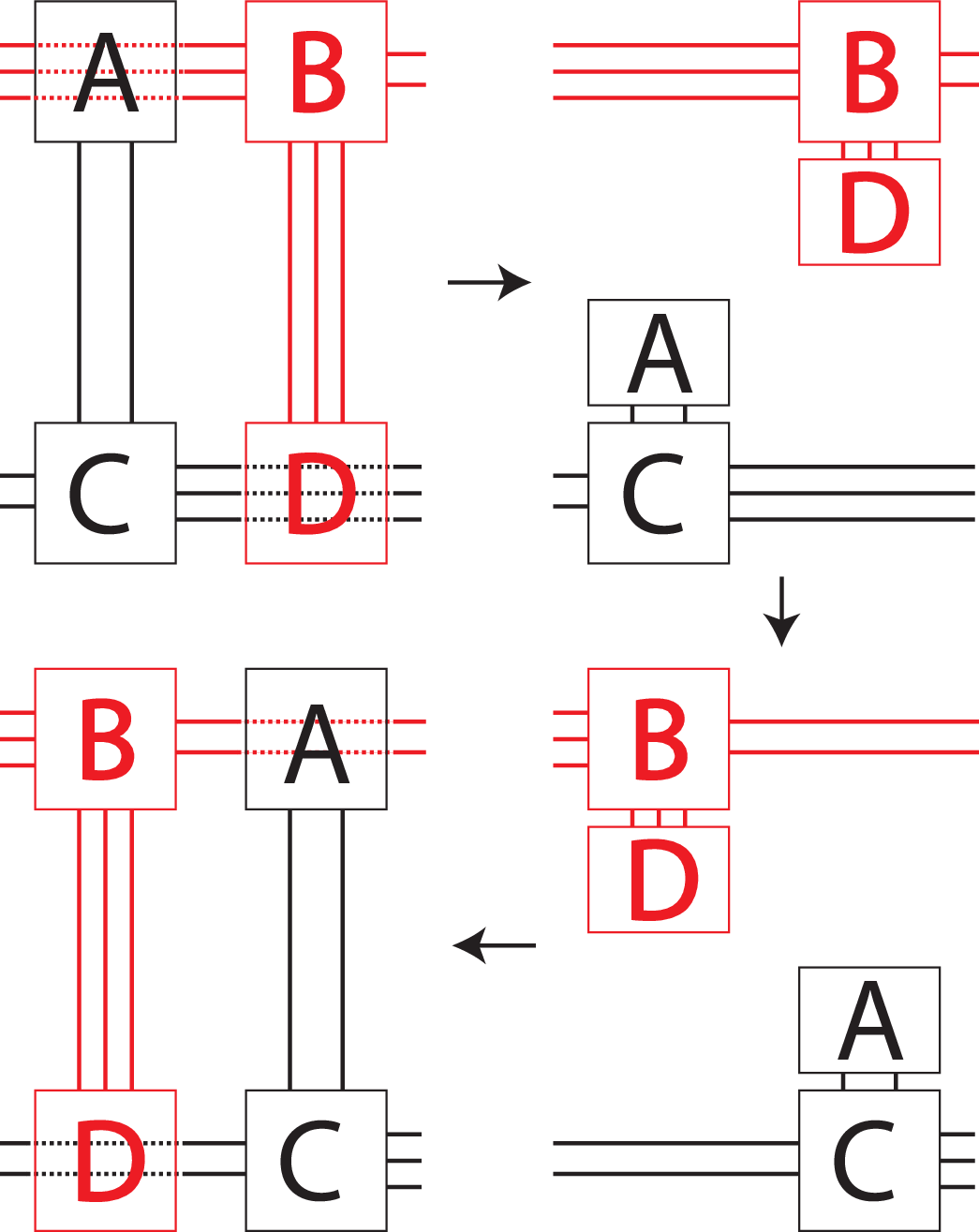}
}
\vskip 6pt
\centerline{Figure 8: A generalised exchange move}

\vskip 6pt
\noindent {\caps 2.6. Generalised destabilisations}
\vskip 6pt

Dynnikov also introduces another move called a {\sl generalised destabilisation}.
Here, one assumes that there are two arcs of $L$, one running
from a vertex $s_1$ to a vertex $s$, and the second running from
$s$ to a vertex $s_2$. Let the $\theta$-values
of these two arcs be $t_1$ and $t_2$. One assumes that
there are no arcs of $L$ with $\theta$ values in $(t_1, t_2)$.
Then, the generalised destabilisation replaces these
two arcs of $L$ by a single arc, running from $s_1$ to $s_2$,
at height $t_2$, say.

\vskip 18pt
\centerline{
\includegraphics[width=5.9in]{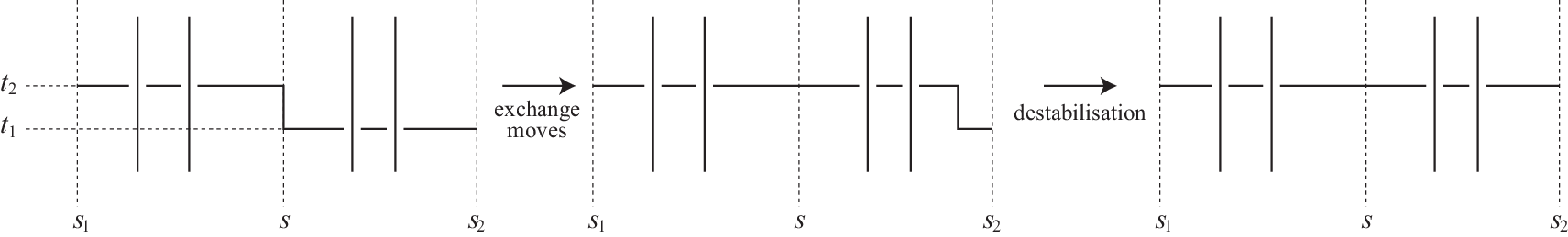}
}
\vskip 6pt
\centerline{Figure 9: A generalised destabilisation}

In Figure 9, a generalised destabilisation is expressed as composition
of exchange moves and a destabilisation. The following is clear.

\noindent {\bf Lemma 2.5.} {\sl Let $n$ be the arc index of an arc presentation
of $L$. Then a generalised destabilisation is the composition of 
at most $n$ exchange moves, followed by a destabilisation.}

\vskip 18pt
\centerline{\caps 3. A summary of Dynnikov's methods}
\vskip 6pt

In this section, we give an overview of Dynnikov's work on monotonic
simplification of arc presentations in [8]. This was highly influenced by
Cromwell's initial investigations into arc presentations in [6]. In turn,
this was influenced by the development of braid theory by Birman and
Menasco (see [4] for example, or the survey in [3]) and Bennequin [2].
Our presentation in this section is substantially based on [8].

\vskip 6pt
\noindent {\caps 3.1. Admissible form for characteristic surfaces}
\vskip 6pt

When $L$ is the unknot or a split link, there is an associated
surface, that Dynnikov refers to as a {\sl characteristic surface}.
In the case of the unknot, this is a spanning disc. For a split link, it is
a 2-sphere disjoint from the link, and with link components on both sides
of it.

This surface $S$ inherits a singular foliation ${\cal F}$ on 
$S - S^1_\phi$ defined by $d \theta = 0$. The intersection points $S \cap S^1_\phi$
are called the {\sl vertices} of $S$.

Dynnikov places the characteristic surface $S$ into {\sl admissible form}, which is
defined as follows:
\item{(1)} The surface $S$ is smooth everywhere, except at $\partial S \cap S^1_\phi$.
\item{(2)} $S - \partial S$ intersects the binding circle $S^1_\phi$ transversely
at finitely many points.
\item{(3)} The foliation ${\cal F}$
has only finitely many singularities, which are points of tangency of $S$ with the
pages ${\cal D}_t$.
\item{(4)} All singularities of ${\cal F}$ are of Morse type, ie local maxima, local minima
or saddle critical points.
\item{(5)} Near any point of $(\partial S) \cap S^1_\phi$, the foliation ${\cal F}$ is
radial.
\item{(6)} There is at most one point $p \in (\partial S) \cap S^1_\phi$ at which
$|\int_\gamma d \theta| > 2 \pi$, where $\gamma \subset S$ is a properly
embedded arc in a small neighbourhood of $p$ such that the
endpoints of $\gamma$ in $\partial S$ lie on different sides of $p$.
Such a point $p$ is called a {\sl winding vertex}. The quantity $|\int_\gamma d \theta|$
is the {\sl winding angle} at this vertex.
\item{(7)} There is at most one point $p \in (\partial S) - S^1_\phi$
at which the surface $S$ is not transverse to the corresponding page
${\cal D}_{\theta(p)}$. At the exceptional point, the foliation ${\cal F}$
must have a saddle critical point. If such a saddle and a winding vertex
are both present, then the winding vertex is an endpoint of the edge
containing the saddle.
\item{(8)} Each page ${\cal D}_t$ contains at most one arc of $L$
and at most one singularity of ${\cal F}|_{S - \partial S}$, but not both.

Consider an arc of $L$, which is the intersection with some page, and suppose
that it does not contain a saddle of $S$. Suppose that $\partial S = L$ (and so
we are in the case where $L$ is the unknot).
Then, near this arc, except at the endpoints, all points of $S$ satisfy
one of the following:
\item{(1)} they have $\theta$-values slightly greater than that of the arc, or
\item{(2)} they have $\theta$-values slightly smaller than that of the arc.

\noindent We term this an {\sl up} or {\sl down} arc, respectively.

Now consider two incident arcs of $L$, neither of which contains a
saddle of $S$. Then, by examining their common vertex, we see
that one must be an up arc and one must be a down arc. So,
as one travels along $L$, one meets up and down arcs alternately,
with the possible exception of an arc containing a saddle.
As a consequence, when the arc index of an unknot $L$ is odd, then
the characteristic surface must have a saddle somewhere on its boundary.

In the case where $L$ is a split link, placing the characteristic 2-sphere
into admissible form is a simple application of general position. However,
when $L$ is the unknot, a little more work is required. 
One first declares that the arcs of $L$ are alternately
up and down arcs, plus possibly one arc that contains a saddle of $S$. This controls
the location of $S$ near these arcs. Near each vertex of $L$, the
foliation is required to be radial. When the vertex is not a winding vertex,
this determines the behaviour of $S$ near that vertex. At the winding
vertex, the amount that the surface winds is chosen so that the curve
$\partial N(L) \cap S$ has zero linking number with $L$.
Thus, one first specifies the location of $S$ near $L$, using
this recipe. Then a small isotopy supported away from a
small neighbourhood of $L$ moves $S$ into admissible form. More details
can be found in the proof of Lemma 1 of [8].

\vskip 6pt
\noindent {\caps 3.2. The structure of admissible surfaces}
\vskip 6pt

Near a singular point of ${\cal F}$ or a vertex of $S$,
there are the following possible local pictures:

\vskip 18pt
\centerline{
\includegraphics[width=4.5in]{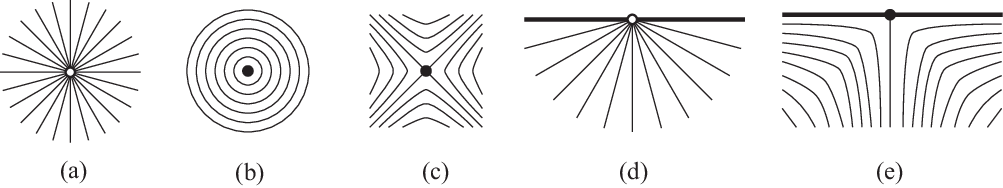}
}
\vskip 6pt
\centerline{Figure 10: Singularities of the foliation}

In (a), the behaviour near a point of $(S - \partial S) \cap S^1_\phi$ is shown.
This is termed an {\sl interior vertex} of ${\cal F}$. A {\sl boundary vertex}
is shown in (d), which is a point of intersection $\partial S \cap S^1_\phi$.
The singularities shown in (b), (c) and (e) are called a {\sl pole}, an {\sl interior
saddle} and a {\sl boundary saddle}. 
We follow Dynnikov by denoting a vertex of $S$ by a hollow dot, and
a Morse singularity by a solid dot.

\vfill\eject

When the singularities are removed from ${\cal F}$, the result is
a genuine foliation on $S - S^1_\phi$. Each leaf is known as a {\sl fibre}.
(Dynnikov also calls the singularities of ${\cal F}$ fibres, but we do not
do so here.) Therefore, fibres are of the following types:
\item{(1)} a closed circle;
\item{(2)} an open arc connecting two vertices;
\item{(3)} an open arc connecting a vertex to a saddle or a saddle
to itself.

\noindent Note that an open arc cannot connect a vertex to itself, other than
possibly a winding vertex. This is because the fibres
emanating from a non-winding vertex have distinct $\theta$-values.
Note also that a fibre cannot connect two distinct saddles,
because each fibre lies in a single page and each
page contains at most one saddle. A fibre that is incident to
a saddle is termed a {\sl separatrix}.

The complement of the vertices, the singular locus, the separatrices and the
boundary of $S$ has a special form. Each component of this complement
we term a {\sl tile}. This has a foliation induced by arcs and curves where $\theta$
is constant. It therefore admits a product structure. Hence,
each tile is an open annulus or an open disc, which we term an {\sl annular}
and {\sl disc} tile respectively. The discs have two vertices in their boundary,
and at most two saddles. (When the boundary of a disc tile runs over fewer than two saddles,
its closure contains an arc of $L$.) Note, however, that the boundary of a
tile may run over the same saddle more than once, as shown in Figure 11.
Hence, the closure of a disc or annular tile need not be a closed disc or annulus. 
There is a type of annular tile that is not shown in Figure 11, which has
boundary consisting of just two vertices. In this case, $S$ is a 2-sphere,
and if it has components of $L$ on both sides of it in $S^3$, then the arc presentation is
disconnected. We may therefore assume that there are no such tiles.

\vskip 18pt
\centerline{
\includegraphics[width=4.5in]{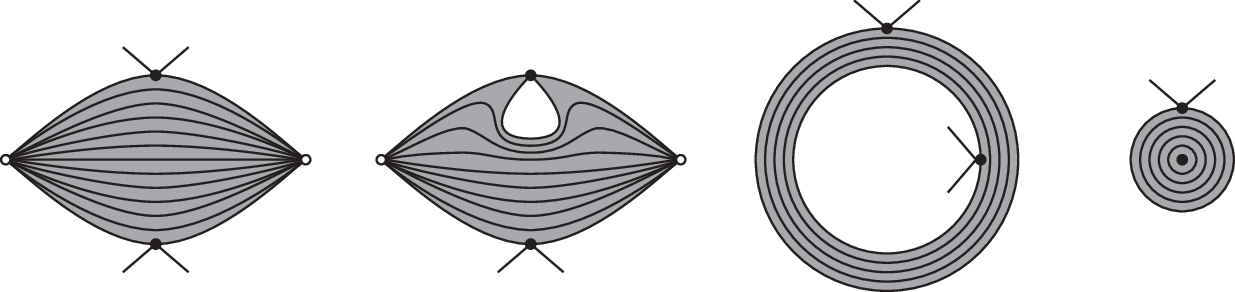}
}
\vskip 6pt
\centerline{Figure 11: Some tiles}

Note that if there are any poles, then there are necessarily
closed circle fibres near them. However, we will see shortly that
poles can be readily removed. Closed circle fibres also
arise near a separatrix that joins a saddle to itself. Note, however,
that in a small neighbourhood of each vertex of $S$, all the
fibres are intervals.

Dynnikov defines the complexity of the characteristic
surface $S$ in admissible form to be the number of
singularities of ${\cal F}$. We will use a slight variation
of this. We will consider the {\sl binding weight}
$w_\beta(S)$, which is the number of intersections between
$S$ and the binding curve $S^1_\phi$. In other words,
the binding weight of $S$ is the number of vertices of $S$,
as shown in Figures 10(a) and 10(d).

\vfill\eject
\noindent {\caps 3.3. Reducing the complexity of the characteristic surface}
\vskip 6pt

In [8], Dynnikov uses an Euler characteristic argument to show
that the singular foliation ${\cal F}$ must contain certain configurations.
In each case, he shows that one may either perform some exchange moves
and cyclic permutations followed by a destabilisation, or one may perform
some cyclic permutations, exchange moves and generalised exchange moves, 
after which one may reduce the complexity of the characteristic surface.
There are $8$ possible configurations that he considers. However,
in this paper, four of these play a particularly important role,
and we will focus initially on these.

For a vertex $s$ of ${\cal F}$, the closure of the union of all the fibres
of ${\cal F}$ approaching $s$ is called the {\sl star} of $s$. The {\sl
valence} of $s$ is the number of separatrices approaching $s$. 

Dynnikov defines an interior vertex $s$ as {\sl bad} if one of the
following cases arises:
\item{(1)} the star of $s$ contains at least two fibres
in distinct tiles that connect $s$ to boundary vertices;
\item{(2)} the star of $s$ contains a winding vertex.

\noindent If an interior vertex is not bad, it is {\sl good}.

The main cases that we consider now are:
\item{(1)} There is a pole.
\item{(2)} There is a good 2-valent interior vertex.
\item{(3)} There is a good 3-valent interior vertex.
\item{(4)} There is a 1-valent boundary vertex.

These are not the only possible cases, but they are the only ones that
we will be concerned with in this paper.
Note that Dynnikov explains in the proof of Lemma 5 in [8] that there can be no
1-valent interior vertex.

\vskip 6pt
\noindent {\caps 3.4. When there is a pole}
\vskip 6pt

In this case, there is a simple modification that can be performed
to the surface which reduces the number of singularities by $2$
without changing the binding weight. One considers
the tile incident to the pole. It has on its boundary
a saddle. One can isotope the surface so as to cancel the pole
and the saddle. This may move other parts of the surface, but
it does not introduce any other singularities. The link itself
does not need to be moved. In particular, no exchange moves, cyclic permutations
or destabilisations are performed at this step.

\vskip 6pt
\noindent {\caps 3.5. When there is a good 2-valent interior vertex}
\vskip 6pt

Suppose that the characteristic surface $S$ has a good 2-valent interior vertex $s$.
Then Dynnikov shows that there is a generalised exchange move
that can be applied to the arc presentation, which leaves
the complexity of $S$ unchanged, and then a further modification
to the surface which reduces its complexity.

\vfill\eject
Adjacent to $s$, there are two disc tiles, and hence the configuration
of ${\cal F}$ near $s$ is as shown in Figure 12.

\vskip 6pt
\centerline{
\includegraphics[width=2in]{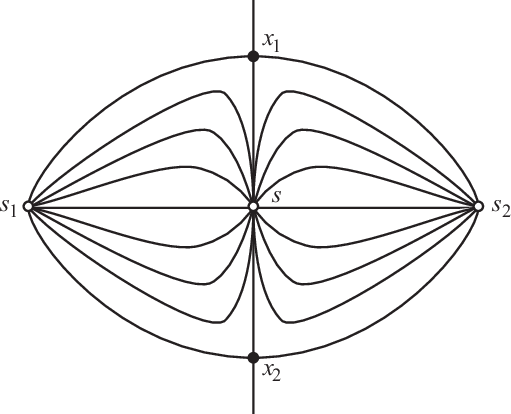}
}
\vskip 6pt
\centerline{Figure 12: A good 2-valent vertex}

The resulting arrangement of the characteristic surface is
shown in Figure 13. Dynnikov explains that, in this situation, one should perform
a generalised exchange move, exchanging the intervals 
$(s_1,s)$ and $(s,s_2)$. This has the effect of modifying the
foliation ${\cal F}$ without increasing its binding weight.
One can then perform an isotopy to $S$, which
reduces its binding weight by $2$.

\vskip 18pt
\centerline{
\includegraphics[width=3.5in]{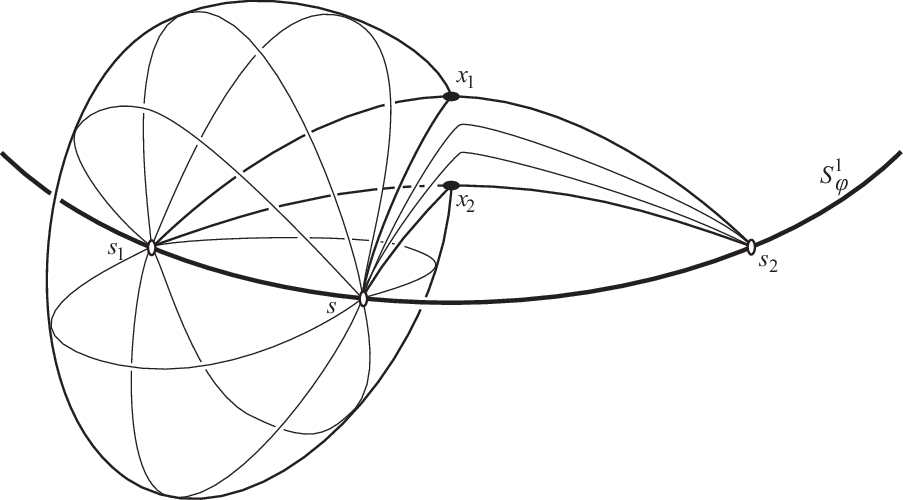}
}
\vskip 6pt
\centerline{Figure 13: The arrangement of the characteristic surface}

\vskip 18pt
\centerline{
\includegraphics[width=4.5in]{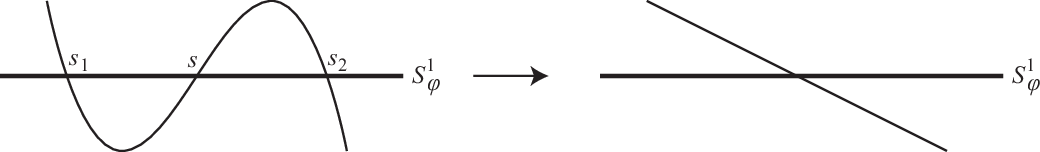}
}
\vskip 6pt
\centerline{Figure 14: The ambient isotopy of $S$}

This procedure does not change the foliation near $\partial S$. In particular,
no new winding vertices or boundary saddles are introduced. Moreover, in the
case where $L$ is the unknot, the decomposition of $L$ into `up' and `down'
arcs, plus possibly one extra arc, remains unchanged.

\vskip 6pt
\noindent {\caps 3.6. When there is a good 3-valent interior vertex}
\vskip 6pt

When there is a good 3-valent interior vertex, Dynnikov explains
how one can isotope $S$ without increasing its binding weight,
to create a good 2-valent interior vertex. This is admirably described in
the proof of Lemma 6 in [8], and so we only give a sketch here.

Let $s$ be the good 3-valent interior vertex. Let $s_2$, $s_3$ and $s_4$
be the three vertices in its star. Without loss of generality, suppose that
they are arranged around $S^1_\phi$ in the order $s$, $s_2$, $s_3$, $s_4$.
Let $x_1$ be the saddle that is connected by separatrices to $s$, $s_2$ and $s_3$,
and let $s_5$ be the other vertex connected to $x_1$ by a separatrix.
Let $x_2$ be the saddle that is connected by separatrices to $s$, $s_3$ and $s_4$,
and let $s_6$ be the other vertex connected to $x_2$ by a separatrix.
A picture of the foliation near $s$ is shown in the left of Figure 15.
Let $t_1 = \theta(x_1)$ and $t_2 = \theta(x_2)$. Suppose, without loss of
generality, that the fibres joining $s$ and $s_3$
have $\theta$ values lying in the interval $(t_1, t_2)$.

The first thing that one does is perform at most $n/2$ cyclic permutations,
so that $0 < t_1 < t_2 < 2 \pi$. Then Dynnikov explains that all events
in the interval $(t_1, t_2)$ need to moved out of this interval, where
an {\sl event} is the occurrence of a saddle or an arc of the link in some
page ${\cal D}_t$, where $t \in (t_1, t_2)$. This is done by moving
the events with endpoints in $(s, s_3)$ into the future, so that
they happen after $t_2$, and by moving
the events with endpoints in $(s_3, s)$ into the past, so that
they happen before $t_1$. In particular, the arcs of the link in these
intervals need to be moved past each other using exchange moves.
Suppose that there are $m$ such arcs with endpoints in the interval
$(s, s_3)$. Then there are at most $n-m$ arcs with endpoints in the interval
$(s_3, s)$ that need to be moved. So at most $m(n-m) \leq n^2/4$
exchange moves are required.

\vskip 18pt
\centerline{
\includegraphics[width=4.5in]{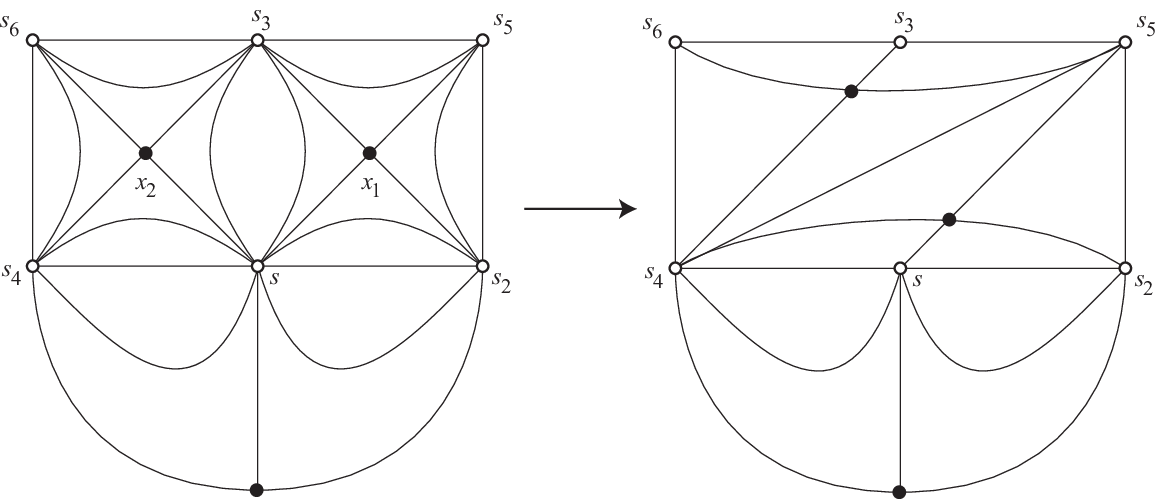}
}
\vskip 6pt
\centerline{Figure 15: A good 3-valent vertex}

Once this has been achieved, one then performs an isotopy,
which has the effect on the foliation as shown in Figure 15.
This turns $s$ into a good 2-valent interior vertex, and so
one then proceeds as in Section 3.5.

As in Section 3.5, this procedure does not change the foliation near
$\partial S$.

\vskip 6pt
\noindent {\caps 3.7. When there is a 1-valent boundary vertex}
\vskip 6pt

In this case, there are two possibilities for the configuration of ${\cal F}$
near the 1-valent boundary vertex $s$. These depend on
whether or not there is a boundary-saddle in the star of $s$.
They are shown in Figure 16.
In both cases, Dynnikov gives a modification to $L$ and $S$.
We concentrate on the case where the star of $s$ does not contain
a boundary saddle. The other case is similar. 

\vskip 18pt
\centerline{
\includegraphics[width=4.5in]{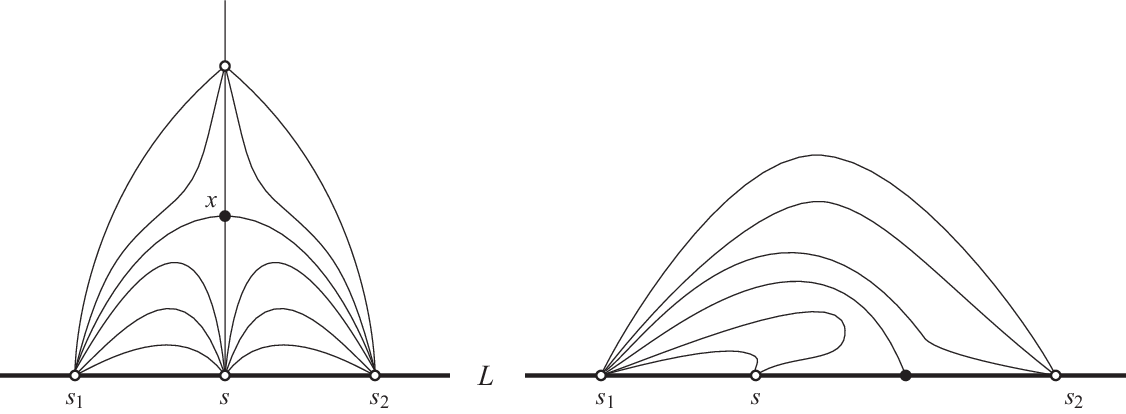}
}
\vskip 6pt
\centerline{Figure 16: A 1-valent vertex $s$}

Now, it cannot be the case that $s$ and $s_1$ are both winding vertices,
since $S$ is admissible.
Hence, the tile that is incident to both of them has total $\theta$-angle
less than $2 \pi$. There is therefore some page that is
disjoint from this tile. We first perform at most $n/2$ cyclic permutations
so that this page is at $\theta = 0$.
We then slide the arc of $L$ that joins $s$ and $s_1$ across this tile, maintaining
it in pages. This has the effect of performing
some exchange moves. At most $n$ of these are performed
in total, because the tile containing $s$ and $s_1$ is disjoint
from the page ${\cal D}_0$.

\vskip 18pt
\centerline{
\includegraphics[width=3.3in]{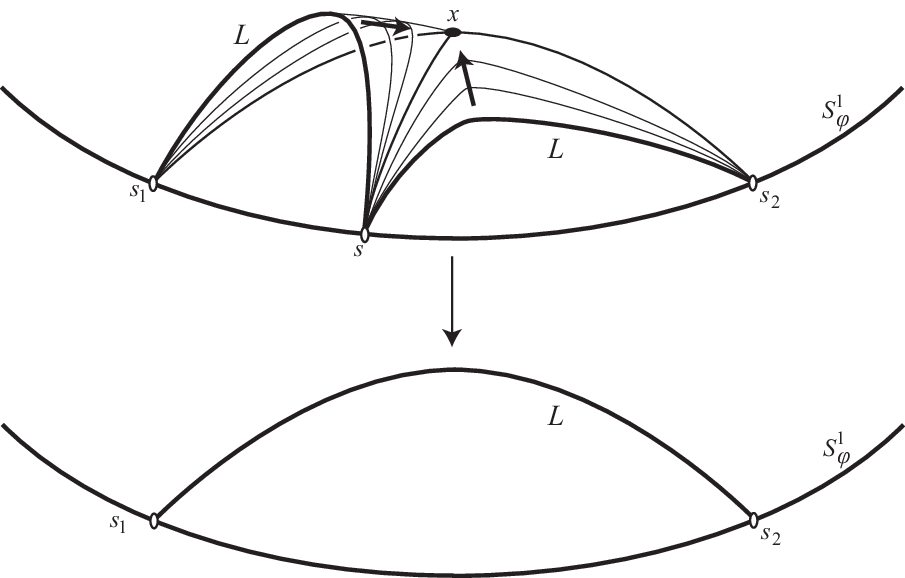}
}
\vskip 6pt
\centerline{Figure 17: The ambient isotopy of $S$}

We now consider the tile containing $s$ and $s_2$. 
We perform at most $n/2$ cyclic permutations so that the tile
misses the page ${\cal D}_0$. Then we slide the arc of $L$
that joins $s$ and $s_2$ across this tile. This process is stopped when the two arcs
of $L$  have adjacent $\theta$-values. Then a generalised destabilisation
is performed. By Lemma 2.5, this is a composition of
at most $n$ exchange moves, followed by a destabilisation.

This procedure does not introduce any winding vertices, since
the $\theta$-angle around each of the vertices $s_1$ and $s_2$
is reduced. However, the resulting surface need not be in admissible form, because
the saddle (labelled $x$ in the left of Figure 16) becomes a boundary-saddle.
If $S$ already has a boundary-saddle elsewhere, then a further isotopy is necessary 
if one wants the resulting surface to be in admissible form. In the next section,
we introduce a variation of admissible form, which we term alternative admissible form,
which is partly designed to get around this complication.

\vfill\eject
\centerline{\caps 4. Simplifying arc presentations of the unknot}
\vskip 6pt

In the previous section, we gave an outline of Dynnikov's argument,
which provides a sequence of exchange moves, cyclic permutations and destabilisations
taking an arc presentation of the unknot or split link to
a trivial or disconnected presentation. The argument relied
on destabilising the arc presentation or reducing the binding weight of the characteristic surface at
each stage. It is not very surprising that the number of
moves that are required can be bounded in terms of the
initial binding weight. In this section,
we prove a result along these lines. The main complication is
that it is not the case that, in Dynnikov's argument, a single
exchange move is used to reduce the binding weight by one.
Many moves may be needed, and these need to be
quantified. It is possible to do this by carefully analysing Dynnikov's proof, but the resulting upper bound on the
number of exchange moves and cyclic permutations is not optimal. Instead,
we present a variant of Dynnikov's theorem and proof, which leads to a better bound.
We are very grateful to Ivan Dynnikov for suggesting that a proof along these lines
would be possible.
This relies on a slightly modified version of admissible form, which is as follows.

Let $L$ be a link with a given arc presentation. Let $S$ be a compact surface embedded in
$S^3$ with interior disjoint from $L$ and with each component of $\partial S$ being a
component of $L$. Then $S$ is in {\sl alternative admissible form} if it satisfies (1),
(2), (3), (4), (5) and (8) in the definition of an admissible surface, together with the following:
\item{(9)} There are no winding vertices.
\item{(10)} Each arc of $L$ contains at most one boundary saddle of $S$.

This has some advantages and some disadvantages over admissible form. The main
disadvantage is that it might not be possible to isotope a given surface into alternative
admissible form, keeping the link fixed. But it is possible to do so after stabilising.

\noindent {\bf Lemma 4.1.} {\sl Let $D$ be an arc presentation of the unknot $L$ with
arc index $n$. Let $S$ be a spanning disc in admissible form. Suppose that it is not
in alternative admissible form, and hence has a winding vertex. 
Let its winding angle be at most $2 \pi m$ for some positive integer $m$.
Then, there is a sequence of $m-1$ stabilisations and at most $(m-1)(n+m)$ exchange moves, taking $D$ to a new
arc presentation $D'$, after which we may isotope $S$ to an alternative admissible surface, 
keeping $L$ fixed. The difference between the binding weight of $S'$ with respect to $D'$
and the binding weight of $S$ with respect to $D$ is $m -1$.}

\noindent {\sl Proof.}  When a stabilisation is performed on an arc presentation, it occurs
near a vertex $s$ of $L$. A new arc of $L$ is inserted into some page ${\cal D}_t$.
If we then perform at most $n$ exchange moves, we may take $t$ to be any value,
as long as this page contains no other arcs of $L$. We may also
suppose that ${\cal D}_t$ contains no singularities of the given admissible surface $S$.
If there is a fibre of the singular foliation on $S$ that is incident to $s$
and that lies in the page ${\cal D}_t$, then there is an obvious way of isotoping $S$ so that,
with respect to the new arc presentation, conditions (1), (2), (3), (4), (5), (6), and (8)
in the definition of admissibility hold. The effect of this on the singular
foliation near $s$ is shown in Figure 18. Away from this regular neighbourhood of $s$,
the singular foliation is unchanged.
We may do this $m - 1$ times at the winding vertex of $S$, so that the resulting surface $S'$
has no winding vertex. Note that each arc of $L$ ends up with at most one boundary saddle.
Hence, this surface is now in alternative admissible form. $\square$

\vskip 18pt
\centerline{
\includegraphics[width=3.5in]{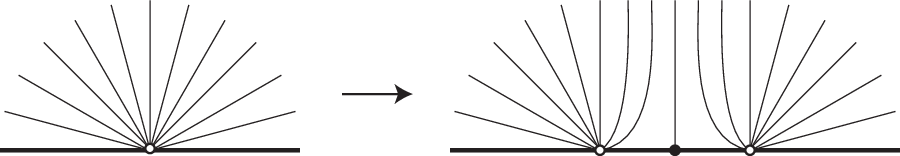}
}
\vskip 6pt
\centerline{Figure 18: Stabilising near a vertex}

We can now give an upper bound on the number of moves required to trivialise an arc presentation
of the unknot.

\noindent {\bf Theorem 4.2.} {\sl Let $D$ be an arc presentation for the unknot $L$ with arc index $n$.
Let $S$ be a spanning disc which is in alternative admissible form, with binding weight $w_\beta(S)$.
Then, there is a sequence of at most $4n^2 w_\beta(S)$ exchange moves, at most $n w_\beta(S)$ cyclic permutations,
at most $w_\beta(S)$ stabilisations and at most $w_\beta(S)$ destabilisations that takes $D$ to the trivial arc presentation.
Moreover, throughout this sequence, the arc index remains at most $n+1$.}

Stabilisations are used here and in Lemma 4.1, 
and so this is not `monotonic simplification' in the sense of Dynnikov [8].

Note that if $w_\beta(S)$ is bounded above by a polynomial function of the arc index $n$, then the number
of moves given by Theorem 4.2 is also bounded above by a polynomial in $n$.

\noindent {\sl Proof.} Because $S$ is in alternative admissible form, it inherits a singular foliation.
The language of admissible surfaces readily translates to this setting. However,
we modify the definition of good and bad vertices, as follows. An interior vertex of $S$
is now {\sl bad} if its star contains fibres $f_1$ and $f_2$ in distinct tiles, 
both of which are incident to boundary vertices, and such that both components of 
$S \backslash {\rm cl}(f_1 \cup f_2)$ contain at least one vertex of $S$.
We say that a boundary vertex is {\sl bad} if its star contains a fibre $f$ that is also incident to some
other boundary vertex, and such that both components of 
$S \backslash {\rm cl}(f)$ contain at least one vertex of $S$.
We say that a vertex is {\sl good} if it is not bad.

We may assume that $S$ has no poles, since if $S$ contains a pole, then 
there is a simple modification to $S$ which reduces its number of singularities without changing its
binding weight and without moving $L$.

For a vertex $s$ of $S$, define its {\sl interior valence} $d_i(s)$ and {\sl boundary valence} $d_b(s)$
to be the number of separatrices approaching $s$, that lie in the interior of $S$ and the boundary of $S$ respectively.
So, the sum of these two quantities is the valence of $s$.

\noindent {\sl Claim.} There is either a good interior vertex with valence $2$ or $3$, or a good boundary vertex $s$ 
such that $2 d_i(s) + d_b(s) \leq 3$.

In order to prove this, we will first construct a graph $G$ embedded in $S$. For each bad interior vertex and for each tile in its star that is
incident to a boundary vertex, pick a fibre in that tile and make it an edge of the graph. For each tile incident to
two boundary vertices and which does not contain an arc of $L$ in its closure, 
pick a fibre in that tile, which runs between these two vertices, and make it an edge of $G$. Take the
vertices of $G$ to be the endpoints of these edges. 

This graph divides
$S$ into discs. We will now pick one of these discs, $S'$, carefully.
If $G$ is empty, then set $S' = S$. So, suppose that $G$ is non-empty.
Let $N(G)$ be a thickening of $G$ away from $\partial S$. This is almost a regular
neighbourhood, except that $N(G) \cap \partial S = G \cap \partial S$.
Let $\alpha$ be $\partial N(G)$. Thus, $\alpha$ is a union of properly embedded
arcs, with disjoint interiors 
but which may intersect at their endpoints. Each arc of $\alpha$ runs parallel to one or two edges of $G$.
We say that an arc $\alpha'$ of $\alpha$ is {\sl trivial} if some
component of $S \backslash \alpha'$ contains no vertices of $S$.
In this case, the corresponding component of $S \backslash G$
contains a single separatrix running from a vertex of $G$ to a boundary
saddle. Let $\alpha_-$ be the resulting of removing all trivial arcs
from $\alpha$. Pick an arc of $\alpha_-$ that is outermost in the disc $S$.
This separates off a disc $S'$ with no arcs of $\alpha_-$
in its interior. Suppose first that $S'$ is not disjoint from $G$. Then
$S'$ contains at least two trivial arcs of $\alpha$, and so we
deduce that $S'$ contains a good boundary vertex $s$ with $d_b(s) = 2$
and $d_i(s) = 0$, as required by the claim. Thus, we may assume that
$S'$ is disjoint from $G$. It therefore corresponds to a component
of $S \backslash G$, which we will also call $S'$.

If $G$ is non-empty, then ${\rm cl}(S') \cap G$ is either a single edge joining two bad boundary vertices
or two edges joined at a bad interior vertex of $S$. Note that, by construction, 
$S'$ contains at least one vertex of $S$, which does not lie in $G$.

Now glue two copies of ${\rm cl}(S')$ along the two copies of ${\rm cl}(S') \cap \partial S$. Denote the resulting surface by $S_+$.
It is either a disc or sphere. This surface $S_+$ has a singular foliation. It has either zero, two or four vertices
in its boundary. In the latter case, at least one of these vertices has valence greater than one.
For if all four vertices in $\partial S_+$ had valence 1, then it is easy to check
that $S'$ contains no vertices, which is impossible.

Note that $S_+$ has no boundary saddles. Let $v_2^i$ and $v_3^i$ be the number of 
interior vertices of $S_+$ with valence $2$ and $3$ respectively. Let $v_1^b$ be the number
of boundary vertices of $S_+$ with valence $1$. Then $v_1^b < 4$.
Using the fact that $S_+$ has positive Euler characteristic,
Dynnikov's argument in the proof of Lemma 5 in [8] gives that $2v_2^i + v_3^i + v_1^b \geq 4$.
(See formula (8) in [8] for example.) Hence, $S_+$ contains in its interior a vertex with
valence at most $3$. This came from a good vertex $s$ of $S$.
When $s$ is in the interior of $S$, it is the vertex required by the claim.
(Note that a vertex in the interior of $S$ cannot have valence 1.)
So, suppose that $s$ lies in the boundary of $S$.
Each separatrix in the star of $s$ that lies in the interior of $S$ gives rise to two separatrices
in $S_+$. Each separatrix in the boundary of $S$ gives rise to just one separatrix of $S_+$.
So, we deduce that $2 d_i(s) + d_b(s) \leq 3$, which proves the claim.

When there is a good interior vertex in $S$ with valence $2$ or $3$, we would like to
apply the procedure described in
Sections 3.5 and 3.6. However, there is one minor complication.
We have modified the definition of a good interior vertex, and so an interior vertex $s$ that was
bad with the previous definition may now be good. In the star of such a vertex $s$, there
are two fibres $f_1$ and $f_2$ lying in distinct tiles, which are incident to
boundary vertices $s_1$ and $s_2$, say, and such that one component of 
$S \backslash {\rm cl}(f_1 \cup f_2)$ contains no vertex of $S$. 
We are concerned with the situation where $s$ has valence $2$ or $3$, and so we now consider these two cases.

Suppose first that $s$ has valence 3.
Then, the local picture near $s$ may not be quite as shown in Figure 15.
One or both of the saddles $x_1$ and $x_2$ may be boundary saddles, in which case
the vertices $s_5$ or $s_6$ might not be present. If $x_1$ and $x_2$ are both boundary
saddles, then we focus instead on $s_3$ which is a good boundary vertex with $d_b(s_3) = 2$ and $d_i(s_3) = 0$.
Such vertices are dealt with later in the argument. So, we may suppose that at most one of
$x_1$ and $x_2$ is a boundary saddle. If $x_2$ is a boundary saddle, the isotopy described in Section 3.6
may still be applied. When $x_1$ is a boundary saddle, we swap the roles of $x_1$ and $x_2$,
and so when we apply the isotopy described in Section 3.6, the resulting foliation is
the mirror image of that shown in the right in Figure 15 without the vertex $s_5$.
Therefore, in both cases, the valence of $s$ can be reduced to 2. It remains good.

So, suppose now that the valence of $s$ is 2. If $s$ is a good interior vertex that was
bad using the previous definition, then
the singular foliation near $s$ is shown in Figure 19. The arrangement of the
characteristic surface still is as shown in Figure 13, but now the arc in $S$ running from $s_1$
to $s_2$ via $x_1$ is actually an arc of $L$. It is clear that the generalised exchange move
and the isotopy of Figure 14
may still be applied, as along as they are combined with a generalised destabilisation of $L$ which
removes this arc.

These procedures reduce the binding weight by 2, and require at most $n$ cyclic
permutations, at most $n^2 + n$ exchange moves and at most one destabilisation.

\vskip 6pt
\centerline{
\includegraphics[width=2in]{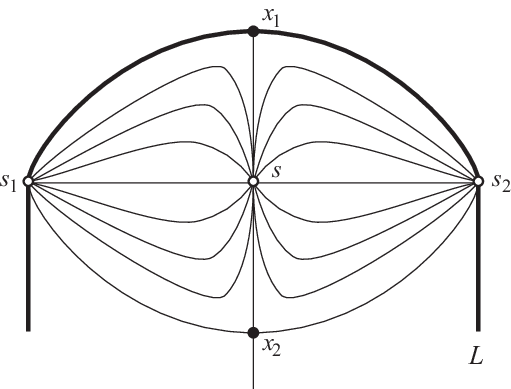}
}
\vskip 6pt
\centerline{Figure 19: A 2-valent interior vertex that is now good}

We now consider the case where there is a good boundary vertex $s$ such that $2 d_i(s) + d_b(s) \leq 3$.
Hence, we are in one of the following situations:
\item{(1)} $d_b(s) = 0$ and $d_i(s) = 0$;
\item{(2)} $d_b(s) = 0$ and $d_i(s) = 1$;
\item{(3)} $d_b(s) = 1$ and $d_i(s) = 0$;
\item{(4)} $d_b(s) = 1$ and $d_i(s) = 1$;
\item{(5)} $d_b(s) = 2$ and $d_i(s) = 0$.

\noindent Note that $d_b(s) \leq 2$, since at most two separatrices in the star of $s$ lie in the boundary of $S$.

We may assume that Case (1) does not arise, 
because a vertex cannot have zero valence, unless the arc presentation is
already is trivial.

Cases (2) and (3) are shown in Figure 16. As explained in Section 3.7, we may apply sequence of at most
$n$ cyclic permutations, at most $3n$ exchange moves and then a destabilisation. After this,
the spanning surface remains in alternative admissible form. Its binding weight has been decreased by $1$.
Note that in Case (2),
the saddle $x$ that is in the star of the vertex becomes a boundary saddle in the new spanning surface. 
The fact that boundary saddles
can be created in this way is one of the reasons why we use alternative admissible form.

In Cases (4) and (5), a new move is required. We will focus on Case (5), but Case (4) is similar.
A picture of the star of $s$ is shown in Figure 20. Note that $s_1$ lies in the interior of $S$,
because $s$ is good. We first perform a generalised stabilisation, which
replaces the arc of $L$ between $s$ and $s_2$ by two arcs, one running from $s$ to $s_1$,
the other running from $s_1$ to $s_2$. By Lemma 2.5, this is a composition of  a stabilisation and
at most $n$ exchange moves. These new arcs of $L$ follow fibres of the foliation of $S$ that lie
near the separatrices incident to $x_1$. The unknot $L$ with this new arc presentation inherits a spanning disc, which is
a subset of $S$, in alternative admissible form. This is shown in the right of Figure 20. 
With respect to this new surface, $d_b(s) = 1$ and $d_i(s) = 0$.
So, we are in Case (3), and therefore a sequence of at most $n$ cyclic permutations, at most $3n$ exchange moves
and then a destabilisation can be performed. Note that, although a stabilisation has been
performed, it is followed by a destabilisation, and so the arc index remains at most $n$
after this process.

\vskip 12pt
\centerline{
\includegraphics[width=4.9in]{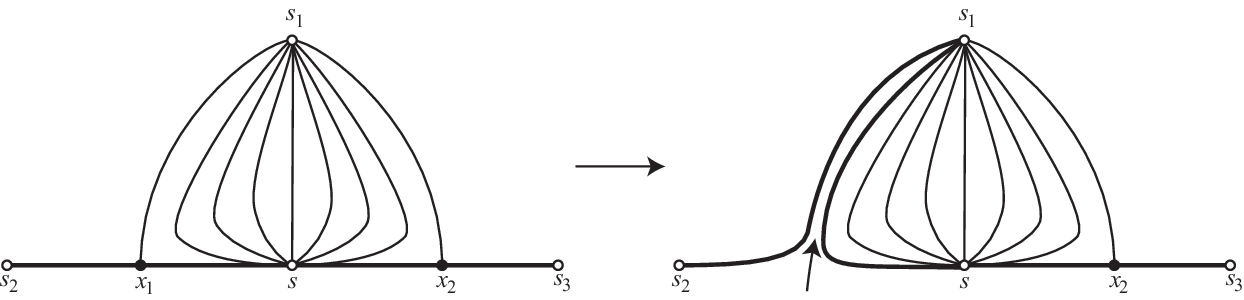}
}
\vskip 6pt
\centerline{Figure 20: Case (5) in the proof}

Since the binding weight has decreased by at least $1$ at each stage, and we have bounded the number
of exchange moves, cyclic permutations, stabilisations and destabilisations at each stage, the theorem
follows immediately. $\square$

\vskip 18pt
\centerline{\caps 5. Normal surfaces}
\vskip 6pt

In this section, we recall some key aspects of normal surface
theory. We also extend the theory a little, by introducing the new
concept of a boundary-vertex normal surface.

\vskip 6pt
\noindent {\caps 5.1. Definitions}
\vskip 6pt

Let $P$ be a compact 3-dimensional polyhedron. Then a disc properly embedded
in $P$ is said to be {\sl an elementary normal disc} if 
\item{(1)} it is disjoint from the vertices and intersects the edges transversely;
\item{(2)} it intersects each face in a collection of properly embedded arcs; and
\item{(3)} it intersects each edge at most once.

When $P$ is a tetrahedron, an elementary normal disc necessarily intersects the 1-skeleton
in three or four points. Normal discs of this form are called {\sl triangles}
and {\sl squares}. Examples are shown in Figure 22.

Let $M$ be a compact 3-manifold with a polyhedral decomposition ${\cal P}$.
Then a surface properly embedded in $M$ is {\sl normal} if
it intersects each polyhedron in a disjoint union of elementary normal discs.

Note that this is a variation on the usual notion of normality. Many authors
require that elementary normal discs satisfy an extra condition:
for each arc of intersection with an interior face,
the endpoints of the arc do not lie on adjacent edges, one
of which is in $\partial M$, while the other is not.
We do not make this requirement here. Our notion of normality
is very close to that used by Jaco and Oertel in [19].

We say that an arc properly embedded in a 2-dimensional polygon 
is {\sl normal} if it is disjoint from the vertices and has endpoints
in distinct edges. When $M$ has a polyhedral decomposition, its
boundary $\partial M$ also inherits a polyhedral structure. We say
that a collection of disjoint simple closed curves in $\partial M$ is
{\sl normal} if its intersection with each face in $\partial M$ is a collection of
normal arcs.

One of the key tenets of normal surface theory is that many topologically
relevant surfaces may be placed in normal form. This is usually proved
by showing that, when a properly embedded surface is not normal,
then there is a modification that can be made to it which reduces
the number of intersections with the 1-skeleton. Hence, a surface
with minimal number of intersections with the 1-skeleton (among a suitable
collection of surfaces) is typically
normal. In fact, when the surface is closed, these modifications do not increase the number
of intersections with any edge. (See Theorem 3.3.21 in [22] for example.) We may therefore obtain a version
of this result which uses a variation of the usual notion of complexity,
which is defined as follows.

Let $M$ be a compact 3-manifold with a polyhedral decomposition ${\cal P}$.
Fix a subcomplex $\beta$ of the 1-skeleton. For a surface $S$ properly
embedded in $M$ in general position with respect to the 1-skeleton of ${\cal P}$,
define the {\sl weight} of $S$, denoted $w(S)$, to be the number of intersection points
between $S$ and the 1-skeleton of ${\cal P}$. Define the {\sl $\beta$-weight}
of $S$ to be the number of intersection points between $S$ and $\beta$,
denoted $w_\beta(S)$. We will consider the pair $(w_\beta(S), w(S))$
and order these pairs lexicographically. Thus, $(w_\beta(S), w(S))$
is less than $(w_\beta(S'), w(S'))$ if and only if {\sl either} $w_\beta(S) < w_\beta(S')$, 
{\sl or} $w_\beta(S) = w_\beta(S')$ and $w(S) < w(S')$.

Note that the terminology $w_\beta(S)$ is already being used to denote the
binding weight of an admissible surface $S$. This is intentional,
because later in the paper, we will choose ${\cal P}$ and $\beta$
so that these quantities coincide.

A straightforward modification to the proof of Theorem 3.3.21 in [22]
gives the following result.

\noindent {\bf Theorem 5.1.} {\sl Let $M$ be a compact orientable
3-manifold with a polyhedral decomposition ${\cal P}$ that has a
subcomplex $\beta$ in its 1-skeleton. Suppose that $M$ is reducible. 
Then there is a reducing sphere $S$ in normal form, such that
$(w_\beta(S), w(S))$ is minimal among all reducing spheres that
are in general position with respect to the 1-skeleton.}

We will also need to work with normal surfaces with boundary.
In this case, the usual normalisation procedure may need to move
the boundary of a surface. With the strong notion of normality
that is used by many authors, this movement of the boundary of the
surface is hard to avoid. However, with the weaker version of normality
we are using in this paper, it is possible to ensure that the
boundary of the surface does not need to be moved, under a fairly
mild hypothesis. The main modification occurs
when there is an arc of intersection between the surface $S$
and an interior face of the polyhedral decomposition with endpoints on the
same edge, and with this edge lying in $\partial M$.
Then, usually one performs a boundary compression to simplify
the surface. If $S$ is orientable, then its boundary inherits
an orientation and we see that, in this situation, the
boundary of the surface intersects this edge in two points
of opposite sign. Thus, if we ensure that this does not arise,
then this modification is not required. We therefore obtain the
following result.

\vfill\eject
\noindent {\bf Theorem 5.2.} {\sl Let $M$ be a compact orientable
3-manifold with a polyhedral decomposition ${\cal P}$ that has a
subcomplex $\beta$ in its 1-skeleton. Suppose that $M$
has compressible boundary. Let $C$ be a normal simple closed curve in $\partial M$
that bounds a disc in $M$. Suppose that, for each edge in $\partial M$,
all points of intersection between $C$ and that edge have the same sign. Then
there is a compression disc $S$ in normal form, with $\partial S = C$,
such that $(w_\beta(S), w(S))$ is minimal among all compression discs that
are in general position with respect to the 1-skeleton and that have boundary
equal to $C$.}

\vskip 6pt
\noindent {\caps 5.2. The normal surface equations}
\vskip 6pt

Let $M$ be a compact 3-manifold with a polyhedral decomposition ${\cal P}$.
Suppose that there are $k$ types of elementary normal discs in ${\cal P}$.
Then each properly embedded normal surface $S$ in $M$ 
determines a sequence of non-negative integers $(x_1, \dots, x_k)$.
Each $x_i$ is the number of elementary normal discs of a fixed type,
and is called the {\sl co-ordinate} of this disc type.
This sequence is known as the {\sl normal surface vector} for $S$, and we denote it by $[S]$.

This vector satisfies a system of linear equations called the
{\sl matching equations}. There is a set of equations for each face $F$ of ${\cal P}$
with polyhedra on both sides. When $S$ is a normal surface properly embedded in $M$,
the elementary discs in the polyhedra adjacent to $F$ intersect $F$
in a collection of normal arcs. For each type of normal arc in $F$,
there must be the same number of arcs of this type from the polyhedra
on both sides. These conditions are the matching equations.

Some elementary normal disc types in a polyhedron
necessarily intersect. We call two discs of this type {\sl incompatible}.
Thus, incompatible elementary discs cannot occur in
a properly embedded normal surface. For example,
in the case of a tetrahedron, two squares of different types necessarily intersect.
Therefore the vector for a normal surface satisfies
the constraints which, for each pair of incompatible disc types,
force the co-ordinate of at least one of them to be zero.
These conditions are called the {\sl compatibility conditions}.

The following key result is one of the cornerstones of normal
surface theory (see Section 1 in [19] for example).

\noindent {\bf Theorem 5.3.} {\sl There is a one-one correspondence
between properly embedded normal surfaces, up to normal isotopy,
and solutions to the matching equations by non-negative integers
that satisfy the compatibility conditions.}

Because of this strong relationship between normal surfaces and
solutions to certain equations, it is useful to take advantage of tools from linear algebra.

The {\sl normal surface solution space} ${\cal N}$ is the set of vectors in ${\Bbb R}^k$ with non-negative
{\sl real} co-ordinates that satisfy the matching equations and the compatibility
conditions. Thus, the points of ${\cal N} \cap {\Bbb Z}^k$ correspond to
properly embedded normal surfaces.

It is easy to see that the normal surface solution space has a polyhedral
structure, in the sense that it is a union of convex polytopes glued
along certain faces. More specifically, suppose that we pick a subset
$Z$ of the co-ordinates, with the property that when two elementary
normal discs are incompatible, at least one of their co-ordinates lies
in $Z$. Consider the set of vectors with real non-negative entries
that satisfy the matching equations, and that satisfy the extra condition
that whenever a co-ordinate lies in $Z$, it is forced to be zero.
We denote this set by ${\cal N}_Z$.
Then ${\cal N}_Z$ is simply the intersection of a subspace of ${\Bbb R}^k$
with the non-negative quadrant $\{ (x_1, \dots, x_k) : x_i \geq 0 \ \forall i \}$.
Hence, it is a cone on a compact polytope. This polytope is just the
intersection of this set with the hyperplane $\{ (x_1, \dots, x_k): x_1 + \dots + x_k = 1\}$.
We denote it by $P_Z$. Note that ${\cal N}$ is the union of ${\cal N}_Z$, over all possible
subsets $Z$.

Let $S$, $S_1$ and $S_2$ be properly embedded normal surfaces. Then
$S$ is said to be the {\sl sum} of $S_1$ and $S_2$ if $[S] = [S_1] + [S_2]$.
We often write $S = S_1 + S_2$.
The sum of $n$ parallel copies of $S$ is denoted by $nS$.
Now, the Euler characteristic of $S$ is a linear function of the number
of elementary normal discs of each type. Hence, when $S = S_1 + S_2$, then
$\chi(S) = \chi(S_1) + \chi(S_2)$.

The normal surface $S$ is a {\sl vertex surface} if it is connected, and whenever $nS$ is the
sum of $S_1$ and $S_2$ for some positive integer $n$, then each of $S_1$ and $S_2$ is a
multiple of $S$.


%


\vskip 6pt
\noindent {\caps 5.3. Realising certain surfaces as vertex surfaces}
\vskip 6pt

Jaco and Tollefson [20] proved that many topologically relevant surfaces
may in fact be realised as vertex surfaces. One of their results is as follows
(see Lemma 5.1 in [20]).

\noindent {\bf Theorem 5.4.} {\sl Let $M$ be a compact orientable
3-manifold with a triangulation $T$. Suppose that $M$ is reducible. 
Then there is a vertex normal surface $S$ that is a reducing sphere,
such that $w(S)$ is minimal among all reducing spheres that are
in general position with respect to the 1-skeleton.}

We will need variation on this result, which differs from it in two ways. 
Firstly, we will not be dealing with a triangulation. Instead, we will 
start with a triangulation ${\cal T}$ (of the 3-sphere) in which the link $L$
is simplicial, and we will remove a small
regular neighbourhood of $L$, forming
a polyhedral structure ${\cal P}$. Now, many of Jaco and Tollefson's
arguments do not extend from triangulations to polyhedral
structures. However, any {\sl closed} normal surface in ${\cal P}$ is
also normal in ${\cal T}$. The arguments of Jaco and Tollefson
do work in this setting. Secondly, we will use a slightly more refined version of
complexity, as in Theorem 5.1.  We therefore obtain the 
following result. The proof of this precisely follows that of Lemmas 5.1 and 4.8 in [20], and is omitted.

\noindent {\bf Theorem 5.5.} {\sl Let ${\cal T}$ be a triangulation
of a compact orientable 3-manifold. Let $M$ be the compact
3-manifold that results from removing a small open neighbourhood
of a subcomplex $L$ of the 1-skeleton. Let ${\cal P}$ be the
resulting polyhedral structure. Let $\beta$ be a subcomplex
of the 1-skeleton of ${\cal P}$. Suppose that $M$ is reducible. Then
there is a reducing sphere that is a vertex normal surface with respect
to ${\cal T}$, and such $(w_\beta(S), w(S))$ is minimal among all reducing spheres that are
in general position with respect to the 1-skeleton.}

\vskip 6pt
\noindent {\caps 5.4. Boundary-vertex surfaces}
\vskip 6pt

When dealing with vertex surfaces, one
loses some control over their boundary behaviour. In order
to get around this, we introduce a new notion.

Let $M$ be a compact orientable 3-manifold with a polyhedral
decomposition ${\cal P}$. Let $S$ be a properly embedded normal surface in $M$.
Then $S$ is a {\sl boundary-vertex surface} if $S$ is connected and whenever $nS$ is the
sum of normal surfaces $S_1$ and $S_2$, where
$\partial S_1$ and $\partial S_2$ are both multiples of
$\partial S$, then each of $S_1$ and $S_2$ is a
multiple of $S$.

\vfill\eject
Boundary-vertex surfaces will play an important role in the proof of our
theorems, in the case of the unknot. We will therefore explore them in
some detail now.

Fix a collection of disjoint simple closed curves $C$ in $\partial M$ that are normal.
The {\sl $C$-normal surface solution space} ${\cal N}^C$ is the set of vectors in 
the normal solution space
${\cal N}$ with boundary that is a multiple of $C$.

As in the case of the usual normal surface solution space, ${\cal N}^C$
is a union of convex polytopes glued along certain faces. This is because
a vector in ${\cal N}$ lies in ${\cal N}^C$ if and only if satisfies a collection
of extra linear equations. Consider two different arc types of normal
arcs in the 2-cells of $\partial M$. Let $c_i$ and $c_j$ be the number
of arcs of $C$ of these two types. For a normal surface $S$, the number of
arcs in $\partial S$ of these two types are linear functions
$\phi_i$ and $\phi_j$ of $[S]$. So, to lie in ${\cal N}^C$,
$[S]$ must satisfy the linear equation $c_j \phi_i [S] = c_i \phi_j [S]$.
These equations, as we run over all pairs of arc types in $\partial M$,
give the extra conditions required
to determine ${\cal N}^C$. Now, just as ${\cal N}$ is a union of the polytopes
${\cal N}_Z$, we may form similar polytopes ${\cal N}_Z^C$
with the above extra linear constraints. So, ${\cal N}_Z^C = {\cal N}^C \cap {\cal N}_Z$.
Then ${\cal N}^C$ is the union of  ${\cal N}_Z^C$ over all possible $Z$.
Note that ${\cal N}_Z^C$ is a cone over a compact polytope,
where the compact polytope is again the intersection with
$\{ (x_1, \dots, x_k): x_1 + \dots + x_k = 1\}$. We denote
this compact polytope by $P_Z^C$.

The condition that $S$ is a boundary-vertex surface is precisely
that $[S]$ is a multiple of a vertex of some $P_Z^{\partial S}$
and that $S$ is connected. The reason for this is as follows. Suppose
that $[S]$ is a multiple of a vertex of some $P_Z^{\partial S}$ and that $nS = S_1 + S_2$
where $\partial S_1$ and $\partial S_2$ are both multiples of
$\partial S$. Then, for each co-ordinate of $S$ that is zero, the
corresponding co-ordinates of $S_1$ and $S_2$ are zero. So, 
$S_1$ and $S_2$ both lie in ${\cal N}_Z^{\partial S}$, and so some multiples
of these surfaces lie in $P_Z^{\partial S}$. However, since 
$S$ is a multiple of a vertex of $P_Z^{\partial S}$, we deduce that both
$S_1$ and $S_2$ are multiples of $S$. Conversely, suppose that 
$[S]$ is not a multiple of any vertex of any $P_Z^{\partial S}$.
Let $Z$ be the set of zero co-ordinates of $S$. Then
a multiple $k[S]$ lies in $P_Z^{\partial S}$ for some positive real $k$. It can therefore be
expressed as an affine linear combination $\lambda_1 v_1 + \dots + \lambda_n v_n$
of the vertices of $P_Z^{\partial S}$, where $\lambda_1 + \dots + \lambda_n = 1$
and each $\lambda_i$ is non-negative. Choose such an expression where
as many of the $\lambda_i$ as possible are zero. After re-ordering, we express
$k[S]$ as $\lambda_1 v_1 + \dots + \lambda_m v_m$ where each $\lambda_i$ is
positive. Since $m$ is minimal, the coefficients $\lambda_1, \dots, \lambda_m$
are uniquely determined. Hence, they are the unique solution to a system of linear
equations with rational coefficients, and therefore they are rational. Rescaling,
we obtain $S$ as a non-trivial sum of surfaces, each with boundary a multiple of $\partial S$,
none of which is a multiple of $S$. Thus, $S$ is not
a boundary-vertex surface.

We will need to realise compression discs as boundary-vertex surfaces.
The precise result, which is an analogue of Theorem 5.5, is as follows.

\noindent {\bf Theorem 5.6.} {\sl Let $M$ be a compact orientable irreducible 3-manifold
with a polyhedral decomposition ${\cal P}$, and a subcomplex $\beta$ in its 1-skeleton. 
Suppose that $\partial M$ is
compressible, and let $C$ be an essential normal simple closed curve in $\partial M$ that bounds a disc
in $M$. Suppose that, for each edge in $\partial M$, all points of intersection between
$C$ and that edge have the same sign. Then there exists a normal disc $S$
bounded by $C$, such that
\item{(1)} $S$ is a boundary-vertex surface, and
\item{(2)} $(w_\beta(S), w(S))$ is minimal among all normal discs with boundary equal to $C$.

}

\vfill\eject
We will now embark upon a proof of this. As mentioned above, the arguments
of Jaco and Tollefson in [20] do not readily translate to the polyhedral setting.
We therefore provide a more direct argument.

We need the following lemma. This is proved in exactly the same way
as Lemma 2.1 in Jaco and Oertel [19], to which we refer the reader for a proof.

\noindent {\bf Lemma 5.7.} {\sl Let $M$ be a compact orientable irreducible
3-manifold with a polyhedral decomposition ${\cal P}$,
with a subcomplex $\beta$ in its 1-skeleton. Let $S$ be a properly embedded,
incompressible, normal surface
such that $(w_\beta(S), w(S))$ is minimal among all surfaces
isotopic to $S$ via an isotopy that keeps $\partial S$ fixed.
Suppose that $S = S_1 + S_2$, and that the number of components of $S_1 \cap S_2$ is
minimal among all normal surfaces $S_1'$ and $S_2'$ such that
$\partial S_1' = \partial S_1$, $\partial S_2' = \partial S_2$,
$S'_1$ and $S'_2$ are isotopic to $S_1$ and $S_2$ keeping
their boundaries fixed and $S = S'_1 + S'_2$. Then no component
of $S_1\cap S_2$ is a simple closed curve bounding a disc in
$S_1$ or $S_2$.}

\noindent {\bf Corollary 5.8.} {\sl Let $M$ be a compact orientable irreducible
3-manifold $M$ with a polyhedral decomposition ${\cal P}$
with a subcomplex $\beta$ in its 1-skeleton. Let $S$ be a properly embedded, incompressible,
normal surface such that $(w_\beta(S), w(S))$ is minimal among all surfaces
isotopic to $S$ via an isotopy that keeps $\partial S$ fixed.
Then $S$ cannot be written as $S_1 + S_2$, where
$S_2$ is a 2-sphere.}

\noindent {\sl Proof.} We may assume that $S_1 \cap S_2$ is
minimal among all normal surfaces $S_1'$ and $S_2'$ such that
$\partial S_1' = \partial S_1$, $\partial S_2' = \partial S_2$,
$S'_1$ and $S'_2$ are isotopic to $S_1$ and $S_2$ keeping
their boundaries fixed and $S = S'_1 + S'_2$. Since $S_2$
is a 2-sphere, each component of $S_1 \cap S_2$ bounds a
disc in $S_2$, which contradicts Lemma 5.7. $\square$

\noindent {\sl Proof of Theorem 5.6.} By Theorem 5.2,
there is a compression disc $S$ in normal form, with $\partial S = C$,
such that $(w_\beta(S), w(S))$ is minimal among all compression discs that
are in general position with respect to the 1-skeleton and that have boundary
equal to $C$.

Note first that this implies that, for each positive integer $n$, $(w_\beta(nS), w(nS))$
is minimal among all collections of $n$ disjoint discs with boundary equal to $nC$.
For if there was a collection of $n$ such discs with smaller complexity, then
one of these discs would have to have complexity less than that of $S$,
which is a contradiction.

Now, $[S]$ lies in the $C$-normal solution solution space. It therefore
lies in some polytope ${\cal N}^C_Z$. This is a cone
on the compact polytope $P^C_Z$. Let $\lambda$ be the
unique real number such that $\lambda [S] \in P^C_Z$.
Now, $P^C_Z$ is the affine hull of its vertices $v_1, \dots, v_m$.
Hence, there are non-negative real numbers $\lambda_1, \dots, \lambda_m$
which sum to $1$ such that $\lambda_1 v_1 + \dots + \lambda_m v_m = \lambda [S]$.
Suppose that as many of the $\lambda_i$ as possible are zero.
We may assume that the first $k$ of them, say, are non-zero
and the remainder are zero. So,
$\lambda_1 v_1 + \dots + \lambda_k v_k = \lambda [S]$.
Divide by $\lambda$ to get an expression
$\mu_1 v_1 + \dots + \mu_k v_k = [S]$.
By our minimality assumption, these real numbers
$\mu_1, \dots, \mu_k$ are unique.
Now each $v_i$ has rational co-ordinates and so
because of the uniqueness of the $\mu_i$s, each
$\mu_i$ is therefore rational. Hence, clearing
denominators, we get an expression
$$n_1 [S_1] + \dots + n_k [S_k] = nS.$$
Here, each $S_i$ is a connected $C$-normal surface, which is a boundary-vertex surface.
Also, $n$ and each $n_i$ is a positive integer.
Hence,
$$n_1 \chi(S_1) + \dots + n_k \chi(S_k) = n \chi(S).$$
Since each $S_i$ is $C$-normal, its boundary consists
of multiples of $C$. So,
$$n_1 |\partial S_1| + \dots + n_k |\partial S_k| = n |\partial S| = n.$$
Therefore,
$$n_1(\chi(S_1) - |\partial S_1|) + \dots + n_k (\chi(S_k) - |\partial S_k|) = 0.$$
There are therefore two cases:
\item{(1)} For some $i$, $\chi(S_i) > |\partial S_i|$.
\item{(2)} For each $i$, $\chi(S_i) = |\partial S_i|$.

Let us consider Case 1 first. Let $\hat S_i$ be the
result of attaching a disc to each boundary component
of $S_i$. Then $\chi(\hat S_i) = \chi(S_i) + |\partial S_i| > 2 |\partial S_i|$.
But $\hat S_i$ is a closed connected surface, and so its Euler characteristic
is at most 2. We deduce that $|\partial S_i| = 0$. Thus, $S_i$ is a 2-sphere
or projective plane. Now, $S_i$ cannot be a projective plane,
for a regular neighbourhood would be a punctured ${\Bbb RP}^3$,
which would force $M$ to be reducible, and this is contrary to assumption.
Therefore, $S_i$ is a 2-sphere. 
We hence get an expression $nS = S_i + W$, for some
normal surface $W$. By Corollary 5.8, this is impossible.

Let us now consider Case 2. Then each $S_i$ is a disc, torus or Klein bottle.
We claim that, in fact, no $S_i$ is a torus or Klein bottle. Suppose it were.
Write $nS = S' + S_i$. Then $S'$ has the same boundary and the same
Euler characteristic as $nS$. It cannot have any 2-sphere or projective plane components,
for this would contradict Corollary 5.8 or irreducibility. Hence, it consists of
$n$ discs, plus possibly some tori and Klein bottles. Let $S''$ be the
union of the disc components of $S'$. Then the total complexity
of $S''$ is strictly less than that of $nS$. Therefore, some component 
of $S''$ has strictly smaller complexity than $S$. This is a contradiction.

We deduce that each $S_i$ must be a disc. So,
$$n_1 + \dots + n_k = n_1 \chi(S_1) + \dots + n_k \chi(S_k) = n \chi(S) = n.$$
Now, 
$$n_1 w_\beta(S_1) + \dots + n_k w_\beta(S_k) = n w_\beta(S).$$
Since each $S_i$ is a disc with boundary equal to $C$, the minimality
assumption on $(w_\beta(S), w(S))$ implies that  $w_\beta(S_i) \geq w_\beta(S)$.
Hence,
$$n w_\beta(S) = (n_1 + \dots + n_k) w_\beta(S)
\geq n_1 w_\beta(S_1) + \dots + n_k w_\beta(S_k) = n w_\beta(S).$$
We deduce that, for each $i$, $w_\beta(S_i) = w_\beta(S)$.
Applying the same argument, we also deduce that $w(S_i) = w(S)$.
Hence, each $S_i$ is a normal disc with boundary $C$
and with minimal complexity. Any of these is our required
boundary-vertex surface. $\square$

\vskip 6pt
\noindent {\caps 5.5. Estimating the size of normal surfaces}
\vskip 6pt

The following is due to Hass, Lagarias and Pippenger (Lemma 6.1 in [14]).

\noindent {\bf Theorem 5.9.} {\sl Let $M$ be a compact 3-manifold
with a triangulation having $t$ tetrahedra. Then, each vertex normal surface $S$,
where $[S] = (x_1, \dots, x_{7t})$, satisfies
$$\max_{1 \leq i \leq 7t} |x_i| \leq 2^{7t - 1}.$$}

\vskip -12pt
We will need the following version of this for compressing discs in polyhedral decompositions.

\noindent {\bf Theorem 5.10.} {\sl Let $M$ be a compact orientable 3-manifold
with a polyhedral decomposition. Let $c$ be an upper bound for the number of elementary normal
disc types in each polyhedron, and let $k$ be the number of elementary
disc types in total. Let $S$ be a compression disc
for $\partial M$  which is a normal boundary-vertex surface.
Let $(x_1, \dots, x_k)$ be the vector $[S]$. Let $y_1, \dots, y_\ell$ denote the weights
of the edges in $\partial M$. Then
$$\max_{1 \leq i \leq k} |x_i| \leq (2c)^{k-1} \left ( \sum_{i=1}^\ell |y_i| \right).$$}

\noindent {\sl Proof.} Consider the following set of linear equations:
\item{(1)} The matching equations.
\item{(2)} The equation $x_i = 0$, for each co-ordinate where $[S]_i$ is zero.
\item{(3)} The equations that specify that $x_i = [S]_i$ for all edges in $\partial M$.

These can be expressed as
$A x = y$, where $A$ is a matrix, $x = (x_1, \dots, x_{k})^T$ and $y$ is
a column vector with the first set of entries being zero, and the remaining entries being
the co-ordinates $y_1, \dots, y_\ell$ of $[\partial S]$. Now, since
$S$ is a boundary-vertex surface, the only solution to
these equations is $[S]$. Hence, $A$ has zero kernel. So, its
rank equals the number of columns. Hence, we may find
a square submatrix $B$ with the same number of columns
and with non-zero determinant. The equations corresponding
to the rows of $B$ become $B x = y'$ for a submatrix $y'$ of $y$.
Inverting, we get $x = B^{-1} y'$. Now, the rows of $B$ have entries
that are $0$, $1$ and $-1$, and there are at most $2c$ non-zero
entries in each row. Also, $B^{-1}$ equals ${\rm adj}(B)/{\rm det}(B)$,
where ${\rm adj}(B)$ is the adjugate matrix. Since $B$ has integral
entries and non-zero determinant, $|{\rm det}(B)| \geq 1$.
Each entry of ${\rm adj}(B)$ is a determinant of a minor
of $B$ and so has modulus at most $(2c)^{k-1}$. 
The required bound on the modulus of each co-ordinate of $x$ immediately
follows. $\square$

\vskip 6pt
\noindent {\caps 5.6. Normally parallel surfaces}
\vskip 6pt

Another useful feature of normal surfaces is that it is possible to speak of
parts of the surface as being normally parallel. The formal definition of this is as follows.

Let $M$ be a compact 3-manifold with a polyhedral decomposition ${\cal P}$. Let $S$ be
a (possibly disconnected) surface properly embedded in $M$ that is in normal form
with respect to ${\cal P}$. Then two subsurfaces $S_0$ and $S_1$ of $S$ are said to be
{\sl normally parallel} if there are subsurfaces $S'_0$ and $S'_1$ of $S$, each of
which is a union of elementary normal discs, and satisfying $S'_0 \supseteq S_0$
and $S'_1 \supseteq S_1$, and an embedding $H \colon S'_0 \times [0,1] \rightarrow M$ such that
the following hold:
\item{(1)} For each elementary normal disc $D$ of $S'_0$ and each $t \in [0,1]$, 
$H(D, t)$ is an elementary normal disc.
\item{(2)} $H(S'_0 \times \{ i \}) = S'_i$ for $i = 0$ and $1$.
\item{(3)} $H(S_0 \times \{ i \}) = S_i$ for $i = 0$ and $1$.

\vfill\eject
\centerline{\caps 6. Triangulations and arc presentations}
\vskip 6pt

\noindent {\caps 6.1. Dynnikov's triangulation}
\vskip 6pt

Dynnikov gave a triangulation of the 3-sphere associated with an
arc presentation of a link $L$. In this subsection, we describe this
triangulation.

As in Section 2, the 3-sphere is viewed as a join $S^1_\theta \ast S^1_\phi$.
Let $n$ be the arc index of the arc presentation.
Then $L$ intersects the binding circle $S^1_\phi$ in $n$ points. The intersection between each page ${\cal D}_t$
and $L$ is either empty or a single open arc. In the latter case, we may assume that
this arc is a concatenation of two arcs which are joined at $S^1_\theta$. We may take
each of these arcs to be $(\phi, \tau, \theta)$, for fixed $\theta$ and $\phi$,
and with $\tau$ varying between $0$ and $1$.

With $L$ in this form, we now define the triangulation of $S^3$,
in which $L$ is simplicial. If $s_1 < \dots < s_n$ are the
vertices $L \cap S^1_\phi$, and $t_1 < \dots < t_n$ are the
points $L \cap S^1_\theta$, we subdivide $S^1_\phi$ and
$S^1_\theta$ at these points. We choose the parametrisation
of $\theta$ and $\phi$ so that these points are equally spaced
around $S^1_\phi$ and $S^1_\theta$. Thus, each circle has been
subdivided into $n$ 1-simplices. We give $S^3$ the triangulation
that is the join of these two triangulations of $S^1_\phi$
and $S^1_\theta$. A typical 3-simplex is therefore of the form $[s_i, s_{i+1}] \ast [t_j, t_{j+1}]$,
for 1-simplices $[s_i, s_{i+1}] \subset S^1_\phi$ and $[t_j, t_{j+1}] \subset S^1_\theta$,
where the indexing is mod $n$.

This triangulation ${\cal T}$ will be of crucial importance in this paper.
In the case where $L$ is a split link, we will arrange that a splitting 2-sphere
is normal with respect to ${\cal T}$. However, when $L$ is the unknot, 
the characteristic surface is a spanning disc, which cannot be made
normal with respect to ${\cal T}$, since $L$ is a subset of the 1-skeleton. 
It is therefore necessary to work with a modified version of the triangulation,
which we define in the next subsection.

\vskip 6pt
\noindent {\caps 6.2. A modification of the triangulation}
\vskip 6pt

The first thing that we do is replace each 1-simplex in $S^1_\theta$
and $S^1_\phi$ by two 1-simplices. We again work with the triangulation
of the 3-sphere that is the join of these triangulations. We denote this also by ${\cal T}$.
This has $4n^2$ tetrahedra. The purpose of doing this is so that, for each
tetrahedron $\Delta$, $L \cap \Delta \cap S^1_\phi$ is at most
one point, and similarly $L \cap \Delta \cap S^1_\theta$ is at most one point.
Hence, for each
tetrahedron $\Delta$, the intersection $\Delta \cap L$ is now
at most two isolated points or a single edge.

We now remove a regular neighbourhood of $L$. 
The effect of this on each tetrahedron is to truncate
some vertices, or slice off an edge. This is shown in Figure 21.
This converts each tetrahedron into a polyhedron.
Let ${\cal P}$ denote the resulting polyhedral decomposition
of the exterior of $L$.

\vskip 6pt
\noindent {\caps 6.3. The number of elementary disc types}
\vskip 6pt

In this subsection, we provide the following crude upper bound on the
number of elementary normal disc types in each polyhedron of ${\cal P}$.

\vskip 6pt
\centerline{
\includegraphics[width=3.4in]{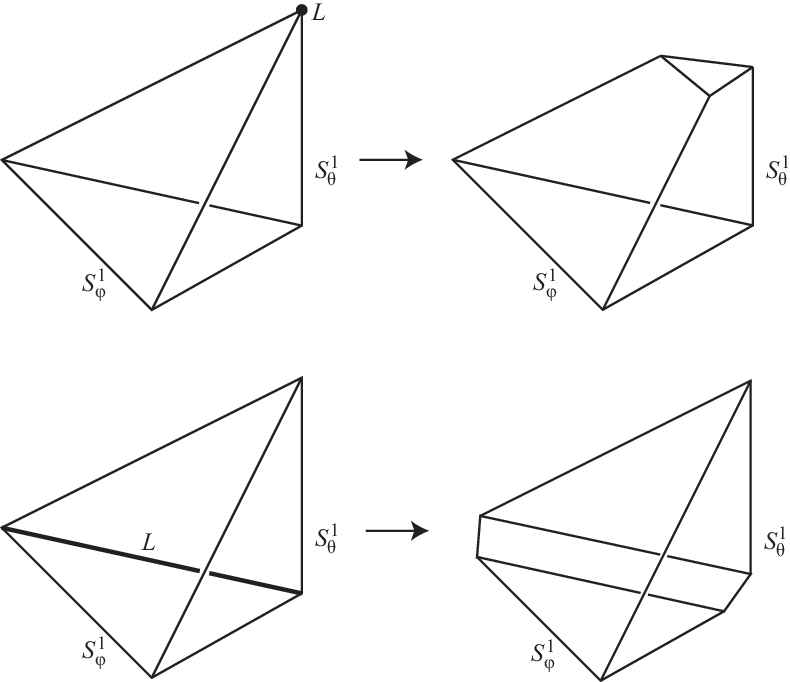}
}
\vskip 6pt
\centerline{Figure 21: Truncating the tetrahedra}
\vskip 6pt

\noindent {\bf Lemma 6.1.} {\sl The number of elementary normal disc types in each polyhedron of 
${\cal P}$ is at most $10^{6}$.}

\noindent {\sl Proof.} Each face of ${\cal T}$ is a triangle. When truncated to form
${\cal P}$, this face ends up with at most $5$ sides. Since an elementary normal
disc can intersect each edge at most once, it therefore intersects this face
in at most 2 normal arcs. There are at most 10 ways of inserting two normal arcs
into the face (since these avoid at most one of the edges, and this avoided edge determines
the normal arcs up to one further choice). There are at most 10 ways of inserting one
normal arc. Hence, there are at most $21$ possible configurations for the intersection
between the face and an elementary normal disc. The normal disc is almost
determined by its intersection with these four faces. The one ambiguity is
when $L$ intersects the tetrahedron in an edge, which is sliced
off to form a rectangular face. Then when the elementary normal disc
intersects all four edges of this rectangular face, there are two possible
ways that it can intersect this face. So, there are at most $2 \times 21^4 < 10^6$
possible elementary normal disc types in each polyhedron. $\square$

\vskip 6pt
\noindent {\caps 6.4. The specified longitude}
\vskip 6pt

The boundary $\partial N(L)$ of this polyhedral structure
inherits a cell structure. Each truncated vertex of a tetrahedron
gives rise to a triangular 2-cell. Each sliced-off edge gives rise
to a rectangular 2-cell.

In the case where $L$ has a single component,
we will now pick a normal, simple closed curve $C$ in $\partial N(L)$
which has winding number one along $N(L)$ and zero linking number with $L$. We will
term this curve the {\sl specified longitude}.
We first create a normal curve $C'$ in $\partial N(L)$, which has winding number
one along $N(L)$, but not necessarily zero linking number with $L$.

Now $L$ is a union of arcs, each of which is the closure of the intersection with
some page. When the arc index $n$ is even, we label these arcs alternately
as {\sl up} and {\sl down} arcs. When the arc index $n$ is odd, this is not
possible, and so we label the arcs as alternately {\sl up} and {\sl down},
with the exception of one arc which is unlabelled. We also pick an orientation
on $L$.

It is also the case that $L$ is a union of edges of the triangulation ${\cal T}$.
For each edge, there are four rectangular 2-cells in $\partial N(L)$ that encircle it.
Two of these rectangles have slightly greater $\theta$-values than the
arc of $L$; two have slightly smaller $\theta$-values. Similar statements
are true for the $\phi$-values. 
We now label the edges of $L$ in this triangulation. If the edge lies in a labelled
arc, we give it the same label. If the edge lies in an unlabelled arc, then we consider the
labelled arc to which it is incident, and give it the opposite label. 
Now arrange $C'$ in the neighbouring rectangles according to the following recipe:
\item{(1)} If the edge of $L$ is labelled `up'  and runs from
$S^1_\phi$ to $S^1_\theta$, then choose $C'$ in this neighbourhood to have slightly greater 
$\theta$-value and slightly greater $\phi$-value.
\item{(2)} If the edge of $L$ is labelled `up' and runs from
$S^1_\theta$ to $S^1_\phi$, then choose $C'$ in this neighbourhood to have slightly greater 
$\theta$-value and slightly smaller $\phi$-value.
\item{(3)} If the edge of $L$ is labelled `down' and runs from
$S^1_\phi$ to $S^1_\theta$, then choose $C'$ in this neighbourhood to have slightly smaller
$\theta$-value and slightly smaller $\phi$-value.
\item{(4)} If the edge of $L$ is labelled `down' and runs from
$S^1_\theta$ to $S^1_\phi$, then choose $C'$ in this neighbourhood to have slightly smaller
$\theta$-value and slightly greater $\phi$-value.

At each point of $L \cap S^1_\theta$ or $L \cap S^1_\phi$, 
there is a collection of triangles of $\partial N(L)$. Coming into these,
there are the endpoints of two arcs of $C'$ lying in rectangular
2-cells. Join these by a path of normal arcs in the triangles
which is as short as possible. (At the point of $L \cap S^1_\theta$ in the middle
of the unlabelled arc, $C'$ will also need to cross some rectangular 2-cells.)
The result is the simple closed curve $C'$.

Suppose that the rectangular diagram associated with this arc presentation
has writhe $k$. Then we claim that the modulus
of the linking number between $C'$ and $L$ is at most $|k|+n+1$. 
We see that $C'$ runs parallel to each vertical and horizontal edge of the rectangular diagram, 
except possibly at the midpoint of just one edge, where it may jump
from one side of the edge to the other. This exceptional case will correspond
to the arc of $L$ containing a boundary saddle. When vertical and horizontal
arcs of the diagram meet at their endpoints, a crossing between $C'$
and $L$ can occur. We deduce that ${\rm lk}(C',L)$ differs from 
the writhe of the rectangular diagram by at most $n+1$. Therefore, $|{\rm lk}(C',L)| \leq |k| + n+1$,
as claimed. Note that $|k|$ is at most the number of crossings of the rectangular diagram, which is at most
$(n-1)^2$, and so we also deduce that $|{\rm lk}(C',L)| < n^2$.

To obtain $C$, we perform some Dehn twists to $C'$,
the twisting curve being a meridian that encircles $L$ half-way
along an edge of $L$. If $L$ has an unlabelled arc, then choose the twisting
curve to be a meridian of one of its edges. We perform enough Dehn twists so
that ${\rm lk}(C,L) = 0$.

We say that the number of Dehn twists that we performed is the {\sl twisting number}
of $C$. Hence, this number is at most $|k|+n+1$.

It is a consequence of the construction that $C$ is normal in $\partial N(L)$ and that, for each edge of
the cell structure of $\partial N(L)$, $C$ intersects that edge in points of the same sign.

We now estimate the weight of $C$, which is the number of points of
intersection with the 1-skeleton of ${\cal P}$. The number of triangles of $\partial N(L)$
at each 0-cell of $L$ is $4n-8$. By construction, $C'$ runs
through at most $2n$ of these. Therefore, the weight of $C'$ is at most $4n^2$.
The creation of $C$ from $C'$ introduces at most $4n^2$ points of intersection.
So, the weight of $C$ is at most $8n^2$.

\vfill\eject
\noindent {\caps 6.5. Making the elementary normal discs piecewise-linear}
\vskip 6pt

Each 3-simplex of ${\cal T}$ may be identified with a Euclidean tetrahedron, 
since it is the join of two edges (one lying in $S^1_\theta$, the
other lying in $S^1_\phi$), which we may take to be Euclidean straight lines.
Each polyhedron in ${\cal P}$ is a subset of a tetrahedron in ${\cal T}$,
which we may choose to be convex. We may also choose the gluing maps
between the faces of adjacent polyhedra to be isometries.

Our goal in this subsection is to realise each elementary normal disc of a normal surface
as piecewise-linear, with respect to the Euclidean structure on the polyhedron
that contains it. So, consider a surface $S$ properly embedded in the
exterior of $L$ that is normal with respect to ${\cal P}$.

We first arrange the points of $S \cap S^1_\theta$ in a certain way.
We have arranged that each 1-simplex in $S^1_\theta$
and $S^1_\phi$ has equal length. We first ensure that each point of $S^1_\theta \cap S$ lies in
the middle half of the 1-simplex in $S^1_\theta$ that contains it. In other words,
it lies closer to the midpoint of this 1-simplex than to either of its
endpoints. This will be technically convenient later in the argument.

We next arrange for $S$ to intersect each face of ${\cal P}$
in straight arcs, without moving their endpoints.
We then arrange for $S$ to lie inside each polyhedron $P$ of
${\cal P}$ in a certain way. The boundary of the elementary normal discs is a union of normal arcs in $\partial P$,
which we have taken to be straight in the Euclidean structure. The elementary normal
discs that are triangles can then be realised as flat. The elementary normal
squares can each be realised as two flat triangles, joined along a straight line. We call these
two triangles {\sl half-squares}. When the
square intersects $S^1_\theta$ and $S^1_\phi$, we choose this straight line so that
it runs between $S^1_\theta$ and $S^1_\phi$. We can choose the straight lines in the remaining
squares so that the union of the squares is embedded.

When a normal surface $S$ is closed, its piecewise-linear structure is now
completely determined. However, when $S$ has non-empty boundary, there
are many more types of elementary disc to consider. We realise these as piecewise-linear in
the following way.

Cut the polyhedron $P$ along a thin regular neighbourhood of the triangles and squares in $S \cap P$, 
creating a union of (possibly non-convex) polyhedra. Each such polyhedron $P'$
is star-shaped, centred at some point $v$, say, in its interior. Create a collection of copies of $\partial P'$
by performing dilations based at $v$ with dilation factor smaller than 1. We create
as many copies as there are components of $S \cap {\rm int}(P')$. The curves $S \cap \partial P'$
are simple closed curves in the sphere $\partial P'$. Hence, there is one, $\alpha$, that is
innermost in $\partial P'$. Attach to $\alpha$ an annulus, which runs to the outermost
dilated copy of $\partial P'$. Take this annulus to be a subset of a cone on $\alpha$
with cone point $v$. Now attach, to the other boundary component of the annulus,
the disc in the dilated copy of $\partial P'$ that it bounds. The resulting disc is the required piecewise-linear
elementary normal disc spanned by $\alpha$. Repeat this procedure with a curve in
$(S \cap \partial P') - \alpha$ that is innermost in $\partial P'$, but this time using
the second-outermost dilated copy of $P'$. Continuing in this fashion, we realise
all of $S \cap P'$ as piecewise-linear.

\vskip 6pt
\noindent {\caps 6.6. PL-admissible form}
\vskip 6pt

Let $S$ be a surface properly embedded in the exterior of $L$
that is normal with respect to ${\cal P}$, and that is piecewise-linear.
As in the case of admissible form, 
this surface $S$ inherits a singular foliation ${\cal F}$ on 
$S - S^1_\phi$ defined by $d \theta = 0$. (See Figure 22 for example.)

\vskip 18pt
\centerline{
\includegraphics[width=1.8in]{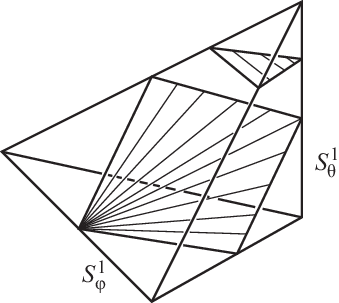}
}
\vskip 6pt
\centerline{Figure 22: Some elementary normal discs with their foliation}

The normal surface is not admissible, for many reasons. It is piecewise-linear,
not smooth. We have yet to make sense of a `singularity' for a such a piecewise-linear
surface, but with any reasonable definition, its singularities cannot be said
to be of `Morse type'. Finally, it need not have the correct behaviour near $\partial S$.
In this subsection, we introduce the notion of PL-admissible form;
the normal surface will have this structure.

Consider a piecewise-linear surface $S$ embedded in ${\Bbb R}^3$ with height function $h$ given
by the final co-ordinate. Suppose that no 1-cell of $S$ is horizontal with respect to $h$.
A point $p$ in $S$ is {\sl non-singular} (with respect to $h$) if
it has a disc neighbourhood $N$ in $S$ such that $\{ x \in N : h(x) = h(p) \}$
is a properly embedded arc in $N$, running through $p$. Otherwise
$p$ is {\sl singular}. We say that a singular point $p$ is a {\sl pole}
if it has a neighbourhood $N$ such that $\{ x \in N : h(x) = h(p) \}$
is just $\{ p \}$.

Note that the singular points are isolated.
The ones that are not poles are {\sl generalised saddles}. An example
of the singular foliation near a generalised saddle is shown in Figure 23.
A generalised saddle $p$ has a disc neighbourhood $N$ such that
$\{ x \in N : h(x) = h(p) \}$ is a star-shaped graph with central vertex $p$.
When $p  \in S - \partial S$, the number of edges of this graph coming out
of $p$ is an even integer at least $4$. When this integer is $4$, we say that
$p$ is a {\sl saddle}. When $p \in \partial S$, the number of edges coming out
of $p$ is an integer at least 2. When this is $2$, we say that $p$ is a {\sl boundary-saddle}
(see Figure 24).

\vskip 6pt
\centerline{
\includegraphics[width=1.3in]{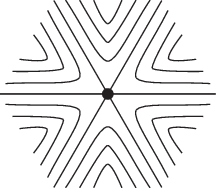}
}
\vskip 6pt
\centerline{Figure 23: A generalised saddle}

A surface $S$ properly embedded in the exterior of $L$
is {\sl PL-admissible} if the following hold:
\item{(1)} It is piecewise-linear in each polyhedron of ${\cal P}$.
\item{(2)} It intersects the binding circle transversely at finitely many points.
\item{(3)} With respect to the function $\theta$ on $S^3 - S^1_\phi$,
$S$ has no horizontal 1-cells and finitely many singularities.
\item{(4)} Each page contains at most one arc of $L$ and
at most one singularity of $S$, but not both.

We can translate the terminology of admissible surfaces to this setting.
A {\sl vertex} of $S$ is a point of $S \cap S^1_\phi$. When the singularities
are removed from the singular foliation ${\cal F}$ on $S - S^1_\phi$,
the result is a genuine foliation. Each leaf is a {\sl fibre}. A fibre that is
incident to a generalised saddle is a {\sl separatrix}. Each component
of the complement of the vertices, the singular locus and the separatrices is a {\sl tile}.
For a vertex $s$ of ${\cal F}$, the closure of the union of all the fibres approaching
$s$ is the {\sl star} of $s$. The {\sl valence} of $s$ is the number of separatrices
approaching $s$.

Note that PL-admissible surfaces have quite different behaviour near
the boundary than in the case of admissible surfaces. This is for several
reasons. In the case of admissible surfaces, their boundary is $L$, which is
a union of arcs in pages. On the other hand, PL-admissible surfaces lie
in the exterior of $L$, and hence have boundary on $\partial N(L)$.
Their boundary curves need not be union of horizontal arcs. In fact, they have
no horizontal arcs in their boundary, because of the assumption that no 1-cell is
horizontal. An example of the singular
foliation near $\partial S$ is shown in Figure 24.

\vskip 18pt
\centerline{
\includegraphics[width=2.8in]{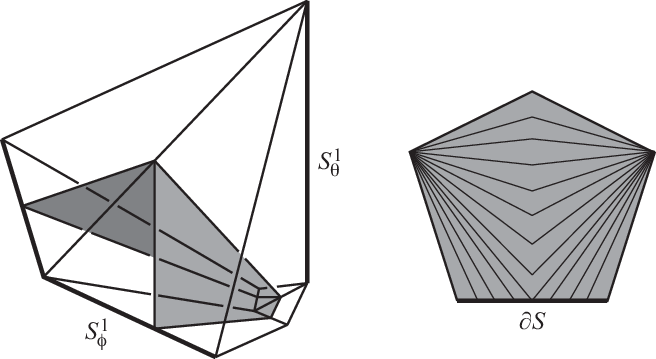}
}
\vskip 6pt
\centerline{Figure 24: Foliation near the boundary}

Let $S$ be a normal surface which has been made
piecewise-linear as described in Section 6.5.

\noindent {\bf Lemma 6.2.} {\sl  In the interior of each elementary normal disc of $S$,
there are at most $24$ singularities. Of these, at most $12$ are generalised
saddles, and all of these are saddles.}

\noindent {\sl Proof.} Note first that a singularity in the interior of a piecewise-linear surface only occurs
when more than two flat discs meet at a point. Moreover, at least four
such discs have to meet a point for this to be a generalised saddle.
At least six such discs have to meet for the point to be a generalised saddle
that is not a saddle.

Let $D$ be an elementary normal disc. When $D$ is a triangle or square, it has no
singularities in its interior. So suppose that $D$ is not of this form. Then
in our construction, $D$ consists
of two parts: an annulus $A$ which runs between $\partial D$ and a dilated
copy of $\partial D$, and a disc which is a subset of a
dilated polyhedron. Singularities that lie in $A - \partial D$ must lie in $\partial A - \partial D$,
and these have at most 4 flat discs incident to them. So, any such
singular points must be poles or saddles. The number of such singularities
is at most the number of points of intersection between $\partial D$ and
the 1-skeleton of ${\cal P}$, which is at most 12, since this
is the maximal number of edges of a polyhedron in ${\cal P}$.
The vertices in the interior of the disc part of $D$ correspond to vertices of the
polyhedron. It is easy to check that there are at most 12 of these.
None of these can be a generalised saddle because they all have
three flat discs incident to them. $\square$

\noindent {\bf Lemma 6.3.} {\sl Each singular point in $\partial S$
has at most two fibres incident to it, and so is a boundary-saddle or a pole.}

\noindent {\sl Proof.} Each elementary normal disc is flat near $\partial S$.
So the only way that a singularity can appear on $\partial S$ is at the
intersection between two elementary normal discs. Since just two
flat discs meet here, this implies that the singular point has precisely
two or zero fibres incident to it. (See Figure 24.) $\square$

\vskip 6pt
\noindent {\caps 6.7. Exceptional and typical separatrices}
\vskip 6pt

An example of a separatrix is shown in Figure 22, running
in the elementary normal square from $S^1_\theta$ to $S^1_\phi$.
We say that a separatrix that lies entirely in an elementary normal square is {\sl typical}.
Otherwise, it is {\sl exceptional}.

\noindent {\bf Lemma 6.4.} {\sl There are at most $408n^2$ exceptional
separatrices of $S$.}

\noindent {\sl Proof.} Each exceptional separatrix emanates from a generalised saddle.
There are three possible locations for a generalised saddle: on the boundary of $S$,
in the interior of an elementary normal disc and on $S^1_\theta$. We consider
these generalised saddles in turn.

The weight of $C$ is at most $8n^2$. Hence,
the number of singularities on $\partial S$ is at most $8n^2$. By Lemma 6.3,
each gives rise to at most two exceptional separatrices.

When a generalised saddle lies in the interior of an elementary normal disc,
this normal disc cannot be a triangle or square, and so it must intersect
$\partial N(L)$. There are at most $8n^2$ such discs. Each contains
at most $12$ generalised saddles in its interior, all of which are saddles,
by Lemma 6.2. So, these give rise to at most $384n^2$ exceptional separatrices.

The remaining separatrices are incident to $S^1_\theta$. To be an exceptional
separatrix, it must start in an elementary normal disc that is not a triangle
or square. There are at most $8n^2$ of these, and each such disc intersects
$S^1_\theta$ at most once. So, we obtain at most $8n^2$ exceptional separatrices of this
form. This gives a total at most $408n^2$ exceptional separatrices. $\square$

\vskip 6pt
\noindent {\caps 6.8. Ordinary tiles and deep vertices}
\vskip 6pt

We say that a tile of $S$ is {\sl ordinary} if it satisfies the following conditions:
\item{(1)} Its closure is a disc disjoint from $\partial S$.
\item{(2)} Its boundary is a union of typical separatrices.
\item{(3)} It lies in the union of the elementary squares and triangles.

\noindent {\bf Lemma 6.5.} {\sl The number of disc tiles that are not ordinary
is at most $1644n^2$. Moreover if $S$ is closed, then every disc tile is ordinary.}

In order to prove this, we will need to introduce the following definition.

For a disc tile $T$ of $S$, we define its {\sl $\theta$-width} as follows.
Let $s_1$ and $s_2$ be its vertices.
Pick a properly embedded arc $\gamma$ in $T$ with endpoints in distinct
components of $\partial T - \{ s_1, s_2 \}$. Then the
$\theta$-width of the tile is $\left | \int_\gamma d\theta \right |$.

\noindent {\bf Lemma 6.6.} {\sl For all disc tiles with at most $816n^2$ exceptions,
the tile has $\theta$-width at least $\pi /2n $.}

\vfill\eject
\noindent {\sl Proof.} By Lemma 6.4, there are at most $408n^2$ exceptional separatrices, 
and these lie in the boundary of at most $816 n^2$ tiles.
Therefore, consider a tile which has no exceptional separatrix
in its boundary. Consider a vertex $s$ of the tile.
The two singular fibres in the boundary of the tile emanating from $s$ lie
in elementary normal squares. These lie in distinct
tetrahedra of ${\cal T}$. There are at most $2n$
tetrahedra arranged around the 1-simplex containing
$s$. They each account for $\theta$-angle around that 1-simplex
of at least $2\pi / 2n$. Since we have arranged
that each point of intersection between $S$ and $S^1_\theta$ lies in
the middle half of the 1-simplex that contains it,
we deduce that the difference in $\theta$ value
between these two singular fibres is at least $2 \pi / 4n$.
$\square$

\noindent {\sl Proof of Lemma 6.5.} If a disc tile is not ordinary, either its closure
intersects $\partial S$ or it contains an exceptional separatrix in its boundary.  

We say that a tile is a {\sl boundary-tile} if it intersects
$\partial S$ in an arc. Not all tiles with closure that intersects $\partial S$ need be
boundary-tiles. This is because the closure of a tile can intersect $\partial S$ at isolated
points, which are boundary saddles. But there are at most $8n^2$
of these. 

We claim that there are at most $820n^2$ boundary-tiles. The total
$\theta$-width of the boundary tiles equals the total $\theta$-angle that
$C$ runs through, in other words, $\int_C |d \theta|$. Now, $C$ is a union
of normal arcs in $\partial N(L)$. As $C$ runs along rectangular faces of $\partial N(L)$,
its $\theta$-angle barely changes, except near the endpoints of the rectangle
that lie near $S^1_\phi$.
At these endpoints, it then runs through triangular faces
of $\partial N(L)$. As it does so, its change in $\theta$-angle is at most $2 \pi$.
So, the total $\theta$-angle that $C$ runs through is at most $2 \pi$
times the number of vertices of $L$, in other words, $2 \pi n$.
Now, for all  but at most $816n^2$ tiles, the $\theta$-width of the tile
is at least $\pi/2n$, by Lemma 6.6. So, the number of boundary-tiles is
at most $816n^2 + 4n^2 = 820 n^2$, as claimed. 

Finally, each exceptional separatrix lies in the boundary of two tiles.
So, by Lemma 6.4, this gives rise to at most $816n^2$ tiles that
are not ordinary. $\square$

We say that two vertices of $S \cap S^1_\phi$
are of the same {\sl type} if their stars are normally
parallel. We say that a vertex is {\sl deep} if its star is disjoint from
$\partial S$ and every separatrix in the boundary of this star
is typical.

\noindent {\bf Lemma 6.7.} {\sl The number of deep vertex types
is at most $48n^2$.}

\noindent {\sl Proof.} Consider a deep vertex and all the vertices
that are of the same type. Their stars are normally parallel. Consider the outermost stars
in this collection. Transversely orient these so that they are both pointing away
from the other stars of the same type. (If there is just one star in the
collection, we consider it twice, with the two different transverse orientations.)
By the definition of a deep vertex, these stars have only typical separatrices in
their boundary. So, each such star is a union of elementary normal triangles,
squares and half-squares. Since each of these stars is not normally parallel 
in the specified transverse direction to 
another star of the same type, we deduce that it contains an elementary normal triangle,
square or half-square that is not parallel
to another elementary normal triangle, square or half-square in the specified transverse
direction. There are at
most $4$ triangle types and at most one square type in $S$ in each truncated tetrahedron,
and there are at most $4n^2$ truncated tetrahedra. Hence, there are
at most $48n^2$ outermost triangles, squares or half-squares in the specified transverse
direction. Each one of these outermost
triangles, squares or half-squares that is in the star of a deep vertex lies in a single tile, and therefore lies in
the star of at most two vertices. This gives the upper bound.
$\square$

\noindent {\bf Lemma 6.8.} {\sl The number of vertices of $S$ that are not deep
is at most $3288n^2$.}

\noindent {\sl Proof.} Each vertex that is not deep lies in the boundary of a disc 
tile that is not ordinary. (Note that annular tiles are not incident to any vertices.)
By Lemma 6.5, at most $1644n^2$ disc tiles are not ordinary.
Each gives rise to two vertices that are not deep. $\square$

\vskip 6pt
\noindent {\caps 6.9. Poles}
\vskip 6pt

\noindent {\bf Lemma 6.9.} {\sl Let $S$ be the characteristic surface in
normal PL-admissible form. Let $w_\beta(S)$ be the number of points
of intersection between $S$ and $S^1_\phi$.
Suppose that $(w_\beta(S), w(S))$ is minimal among all
normal characteristic surfaces with boundary equal to $\partial S$.
Then, $S$ has at most $208n^2$ poles. Moreover, if $S$ is closed, then in fact
it contains no poles.}

\noindent {\sl Proof.} Let $p$ be a pole of the foliation. We claim that $p$
has non-empty intersection with an elementary normal disc that is not a triangle or square.
Now, the interior of each elementary triangle and square contains no poles.
So, a pole that is only incident to triangles and squares must lie in the
2-skeleton of ${\cal P}$. By construction, it does not lie in the interior
of a face of ${\cal T}$. In fact, near edges of ${\cal T} - (S^1_\theta \cup S^1_\phi)$,
the foliation also has no singularities. Since only vertices of $S$ lie
on $S^1_\phi$, we deduce that the pole $p$ lies on $S^1_\theta$.
Let $e$ be the edge of the triangulation
containing $p$. Let $B$ be the union of the tetrahedra
incident to $e$. If a square is incident to $p$, then it contains a
fibre ending on $p$, and so $p$ is not then a pole. Thus, $p$
is only incident to triangles. The union of these triangles
is a disc $D$ properly embedded in $B$. It forms the link in $B$
of one of the endpoints $x$ of $e$. Let $D'$ be the remainder of the link of $x$ in ${\cal T}$.
Note that $D$ and $D'$ have the same number of triangles, by the way that
${\cal T}$ is constructed. Remove $D$ from $S$,
replace it with $D'$. Then perform a further small
isotopy which makes the surface transverse to the 1-skeleton of ${\cal P}$.
This leaves $\partial S$ unchanged, and it also
does not change $w_\beta(S)$. But it has decreased $w(S)$.
By Theorem 5.1 or 5.2, there is a normal characteristic surface
with the same boundary as $S$ but smaller complexity.
This is contrary to hypothesis, proving the claim.

So, consider an elementary normal disc that is not a triangle or square.
There are at most $8n^2$ of these.
By Lemma 6.2, it contains at most $24$ poles in its interior. Any pole
on its boundary lies on $\partial S$ or $S^1_\theta$. The elementary normal disc intersects $S^1_\theta$
at most once. So, the number of poles not lying on $\partial S$ is at most
$200 n^2$. There are at most $8n^2$ poles lying on $\partial S$.

Note that the claim also implies that, when $S$ is closed, it contains no poles. This is
because $S$ then consists only of triangles and squares.
$\square$

\vskip 6pt
\noindent {\caps 6.10. Moves on PL-admissible surfaces}
\vskip 6pt

Dynnikov's argument, described in Section 3, dealt with admissible surfaces.
But many of these arguments work just as well with PL-admissible surfaces.
For example, we have the following result.

\noindent {\bf Proposition 6.10.} {\sl Let $D$ be an arc presentation of a link $L$
with arc index $n$. Let $S$ be a PL-admissible surface properly embedded
in the polyhedral decomposition ${\cal P}$.
\item{(1)} Suppose that $S$ contains a deep 2-valent vertex. Then there is a generalised
exchange move on the link, followed by an ambient isotopy of the link complement,
taking $S$ to a surface $S'$ such that $w_\beta(S') \leq w_\beta(S) - 2$.
\item{(2)} Suppose that $S$ contains a deep 3-valent vertex. Then there is a sequence of
at most $n/2$ cyclic permutations, at most $n^2/4$ exchange moves, a generalised
exchange move and some ambient isotopies on the link complement,
taking $S$ to a surface $S'$ such that $w_\beta(S') \leq w_\beta(S) - 2$.

}

\noindent {\sl Proof.} (1) A picture of the star of a deep 2-valent vertex is shown in Figure 12.
The saddles $x_1$ and $x_2$ shown there may be generalised saddles,
and so may have many separatrices emanating from them, but this does
not affect the argument. The arrangement of the characteristic surface
is shown in Figure 13, and one may make the same ambient isotopy
which reduces the number of intersections with the binding circle by 2.

(2) A picture of the star of a deep 3-valent vertex is shown in Figure 15,
but in the case where all generalised saddles are actual saddles.
When $x_1$ and $x_2$ are generalised saddles,  then there
may be several vertices in Figure 15 between $s_5$ and $s_2$, and between
$s_6$ and $s_4$, which are joined by separatrices to $x_1$ and $x_2$
respectively. But we can nevertheless perform the modification described in
Figure 15 without involving these vertices. This requires at most 
$n/2$ cyclic permutations and at most $n^2/4$ exchange moves.
It converts the deep 3-valent vertex into a deep 2-valent one. 
We then proceed as in (1). $\square$

However, we require a stronger version of this, which involves many vertices
at a time. This is absolutely central to this paper.

\noindent {\bf Proposition 6.11.} {\sl Let $D$ be an arc presentation of a link $L$
with arc index $n$. Let $S$ be a PL-admissible surface properly embedded
in the polyhedral decomposition ${\cal P}$.
\item{(1)} Suppose that $S$ contains $m$ deep 2-valent vertices, all with normally
parallel stars. Then there is a generalised
exchange move followed by an ambient isotopy, 
taking $S$ to a surface $S'$ such that $w_\beta(S') \leq w_\beta(S) - 2m$.
\item{(2)} Suppose that $S$ contains $m$ deep 3-valent vertices, all with normally
parallel stars. Then there is a sequence of at most $n/2$ 
cyclic permutations, at most $n^2/4$ exchange moves, a generalised
exchange move and some ambient isotopies,
taking $S$ to a surface $S'$ such that $w_\beta(S') \leq w_\beta(S) - 2m$.

}

\noindent {\sl Proof.} This follows the above argument. However, $m$ copies
of the surface shown in Figure 12 or 15 are used, all of which are parallel.
Thus, in (1), once the generalised exchange move is performed,
the ambient isotopy shown in Figure 14 can be applied, which reduces
the number of intersections with $S^1_\phi$ by $2m$. The argument in (2)
is similar. $\square$

\vskip 6pt
\noindent{\caps 6.11. Relating admissible and normal surfaces}
\vskip 6pt

In this paper, we are considering four types of surface: admissible surfaces, alternative admissible surfaces,
PL-admissible surfaces, and normal surfaces. It will be crucial to be able to pass between these different types of
surface, as each will play an important role. In this subsection, we explain
how to do this in one direction, while maintaining control of the complexity
of the surfaces.


\noindent {\bf Proposition 6.12.} {\sl Let $D$ be an arc presentation
for the unknot $L$ with arc index $n$. Let $S$ be a compression disc for
$\partial N(L)$ in $S^3 - {\rm int}(N(L))$ which is in normal PL-admissible form with respect to ${\cal P}$. Suppose that $\partial S$ is
equal to the specified longitude, and that its twisting number is $t$. Then there is a characteristic surface $S'$ for 
$L$ which is in admissible form such that $w_\beta(S') \leq w_\beta(S) + n$.
Moreover, if $S'$ contains a winding vertex, then its winding angle is at most $2 \pi t$.
}

\noindent {\sl Proof.} Let $N(L)$ be the regular
neighbourhood of $L$ which is removed when forming ${\cal P}$. 
Let $N_-(L)$ be a much smaller regular
neighbourhood of $L$. We initially set $S'$ to equal $S$ in $S^3 - {\rm int}(N(L))$.

When Dynnikov shows in [8] how a characteristic surface may be placed in
admissible form, the first thing that he does is to arrange it near $L$
so that it has the correct boundary behaviour. (See Section 3.1.) We do the same here,
so that the characteristic surface $S'$ lies in $N_-(L)$ in this specified way. It therefore picks up
$n$ intersection points with $S^1_\phi$, which are precisely
the vertices of the arc presentation.

We now need to explain how to arrange $S'$ in $N(L) - N_-(L)$.
Note that in this region, there lie the arcs $S^1_\phi \cap (N(L) - N_-(L))$,
which are vertical in its product structure. We need to ensure that
$S'$ has no intersection points with these arcs. Then the binding weight
of $S'$, which is just the number of intersection points with $S^1_\phi$, is
$w_\beta(S) + n$.

Now, the two curves $S \cap \partial N(L)$ and $S' \cap \partial N_-(L)$
are already fixed. Using the product structure on ${\rm cl}(N(L) - N_-(L))$,
we may identify $\partial N(L)$ and $\partial N_-(L)$, and
therefore view these two curves as lying on the same torus.
The former curve is equal to the specified longitude,
and the latter is arranged according to the recipe given by Dynnikov,
as described in Section 3.1.
But the specified longitude is defined precisely so that these
are equal, up to an ambient isotopy in the complement of
$S^1_\phi$. Thus, there is a way of inserting $S'$ into this
product region, so that it is an annulus interpolating between
these two curves, and without introducing any new intersection points
with $S^1_\phi$. 

Note that $S'$, as constructed, is piecewise-linear, not smooth. Also, its
singularities are poles and generalised saddles. But a small ambient isotopy,
supported away from $N_-(L)$, makes $S'$ smooth with Morse-type singularities.
This does not change its binding weight, and it turns $S'$ into an admissible
surface. 

In the definition of the specified longitude in Section 6.4, a normal curve $C'$ was first
defined. The specified longitude was obtained from $C'$ by performing
$t$ Dehn twists along a meridian of $L$. This is the location for a winding
vertex of $S'$ (if it has one). By construction, its winding angle is therefore at most
$2 \pi t$. $\square$


%
%

\vskip 18pt
\centerline {\caps 7. The Euler characteristic argument}
\vskip 6pt

\noindent {\bf Theorem 7.1.} {\sl Let $L$ be the unknot or a split link. Fix
an arc presentation of $L$ with arc index $n$ that is not disconnected. Let $S$ be a characteristic
surface in PL-admissible normal form with respect
to the polyhedral decomposition ${\cal P}$, as described in Section 6.5.
Suppose that $S$ is a boundary-vertex surface. In the case
where $L$ is the unknot, suppose also that $\partial S$ is the
specified longitude. Let $w_\beta(S)$ be the number of points 
in $S \cap S^1_\phi$, and let $w(S)$ be the weight of $S$. Suppose that $(w_\beta(S), w(S))$ is minimal
among all characteristic surfaces with the same boundary as $S$. Then, the number of deep
2-valent and 3-valent vertices is at least
$${w_\beta(S) \over 2 \times 10^9 n^4} - 4833n^2.$$
Moreover, if $L$ is a split link and hence $S$ is closed, the number of such vertices is at least
$${w_\beta(S) \over 2 \times 10^9 n^4}.$$
}

We now define a Euclidean subsurface of $S$, which will play a 
key role in the proof. The {\sl designated Euclidean subsurface} $E$ of $S$ is obtained as follows.
It includes the interiors of the ordinary tiles. If two such tiles are adjacent along a separatrix, 
add the interior of this separatrix. If a vertex is 4-valent and is completely
surrounded by ordinary tiles, add it in. Similarly, if a saddle
is completely surrounded by ordinary tiles, add it in. 

We now give $E$ a Riemannian metric that is locally isometric to the
Euclidean plane. Each ordinary tile has, by definition, only typical separatrices in
its boundary, each of which runs between a vertex of $S$ and a generalised saddle.
So it has precisely 4 separatrices in its boundary (as in the
left of Figure 11). We realise it as the interior of a Euclidean square with side length 1.
When the interiors of edges are added, they are realised as
Euclidean geodesics with length 1. The Euclidean metric extends
over the vertices and saddles that are added to form $E$, in a natural way.

Denote the {\sl combinatorial  length} $\ell(\partial E)$ of $\partial E$ 
to be the number of separatrices in $\partial E$ plus the number of components of
$\partial E$ that are isolated points. 


\noindent {\bf Lemma 7.2.} {\sl Suppose that each
point of $E$ is at a distance at most $R$ from $\partial E$. Then, the
area of $E$ is at most $\pi (R+1)^2 \ \ell(\partial E)$.}

\noindent {\sl Proof.} For each point $y$ in $E$, there is a shortest
path from $y$ to $\partial E$, which is a Euclidean geodesic.
Let $x$ be the endpoint of this geodesic in $\partial E$.
Then $y$ lies in the image of the exponential map based at $x$.
Call this map $\exp_x$. It is defined on a star-shaped subset
of $T_x E$ centred at the origin, which we denote by ${\rm dom}(\exp_x)$.
In fact, if we set $S(R,x)$ to be $\exp_x(B(R,0) \cap {\rm dom}(\exp_x))$,
then $y$ lies in $S(R,x)$. Thus, we have shown that $E$ equals
$\bigcup_{x \in \partial E} S(R,x)$.

We now show in fact that $E$ equals the union of $S(R+1, x)$, as $x$
runs over all 0-cells in $\partial E$. By a 0-cell, we mean a corner
of one of the tiles, which may be a generalised saddle or vertex of ${\cal F}$. For suppose that
$\alpha$ is a shortest geodesic joining $y$ to $\partial E$
and that its endpoint $x$ is in the interior of a side of one of the tiles.
Then $\alpha$ is orthogonal to this side. So,
if we slide $x$ to one of the endpoints $x'$ of this side, keeping $\alpha$ a geodesic,
then it remains in the same set of tiles. In particular, it remains in $E$.
This process increases the length of $\alpha$ by at most $1$.
So, $y$ lies in $S(R+1, x')$.

Now, $\exp_x$ is a local isometry from $B(R+1,0) \cap {\rm dom}(\exp_x)$
onto $S(R+1,x)$. Hence, the area of $S(R+1,x)$ is at most 
$\pi(R+1)^2$. So, the area of $E$ is
at most $\ell(\partial E)$ times the maximal
area of $S(R+1,x)$, which gives the required bound. $\square$

The proof of the following key result will take up the entirety
of Section 9.

\noindent {\bf Theorem 7.3.} {\sl Let $L$, $n$ and $S$ be as in Theorem 7.1.
Let $E$ be the designated Euclidean subsurface
of the characteristic surface $S$. Then each point of $E$ has distance
at most $8000n^2$ from $\partial E$.}

\noindent {\sl Proof of Theorem 7.1.} Let $\Gamma$ be the following 1-complex
embedded in $S$.  Its 1-cells are the separatrices. Its 0-cells are 
the endpoints of these separatrices, plus the poles.
This includes the vertices of $S$, the generalised saddles and the 
endpoints of separatrices on $\partial S$.

Let $S_+$ be two copies of $S$ glued
along $\partial S$ via the identity map. (So, when $S$ is a sphere, $S_+$ is two 2-spheres.) So,
$\chi(S_+) \geq 2$. Let $\Gamma_+$ be the union of the copies of $\Gamma$
in $S_+$. Then, $S_+ - \Gamma_+$ is a collection of open annuli
and discs. Let $S_-$ be the result of removing the open
annuli from $S_+$. Then $\chi(S_-) = \chi(S_+) \geq 2$.

Now, $S_-$ inherits a cell structure. When a separatrix in $S$ ends at
a non-singular point on $\partial S$, then combine the two copies
of this separatrix in $S_-$ into a single 1-cell. So, each 0-cell of $S_-$
comes from a vertex, pole or generalised saddle of ${\cal F}$.
For $i=0$, $1$ and $2$, let $S_-^i$ denote the $i$-cells of $S_-$. For each 0-cell $v$, let $d(v)$ denote
its valence.

Each 2-cell of $S_-$ has precisely four 1-cells
in its boundary that are not loops. Hence,
$2|S_-^1| \geq 4|S_-^2|$. Therefore,
$$2 \leq \chi(S_-) = |S^0_-| - |S^1_-| + |S^2_-| \leq |S_-^0| - |S_-^1|/2  = 
 \sum_{v \in S_-^0} (1 - d(v)/4).$$
Let $V$ denote the set of vertices of $S$, let $P$ denote the set of poles of $S$ and let $X$ denote the set of 
generalised saddles in the interior of $S$. The boundary-saddles
give rise to 0-cells of $S_-$ with valence 4, and so they do not contribute
to the above summation. Therefore,
$$2 \leq \sum_{v \in S_-^0} (1 - d(v)/4) = 2|P| + 2\sum_{v \in V} (1 - d(v)/4) + 2 \sum_{x \in X} (1 - d(x)/4).$$
For $k \geq 2$, let $v_k$ denote the number of vertices in $S$ with valence $k$.
Note that there are no vertices of valence 1 in the interior of $S$, as explained in the proof of Lemma 5 in [8].
Note also that $|P| \leq 208 n^2$, by Lemma 6.9. 
Moreover, $|P|$ is zero when $S$ is closed. So, $|P| \leq 208 n^2 |\partial S|$. So,
$$2v_2 + v_3 = \sum_{\scriptstyle v \in V \atop \scriptstyle  d(v) < 4} (4 - d(v))
\geq 4 + \sum_{\scriptstyle v \in V \atop \scriptstyle  d(v) > 4} (d(v) - 4)
+ \sum_{x  \in X} (d(x) - 4) - 832n^2|\partial S|.$$
Note that
$$\sum_{\scriptstyle  v \in V \atop \scriptstyle  d(v) > 4} (d(v) - 4) = 
\sum_{k > 4} \sum_{\scriptstyle  v \in V \atop \scriptstyle d(v) = k} (d(v) - 4)
= \sum_{k > 4} v_k (k-4) \geq \sum_{k > 4} v_k.$$
Similarly, because each generalised saddle in the interior of $S$ has even valence at least 4,
we deduce that
$$\sum_{x \in X} (d(x) - 4) \geq {1 \over 3} \sum_{\scriptstyle  x \in X \atop \scriptstyle d(x) \not= 4} d(x).$$
So,
$$2v_2 + v_3 
> {1 \over 3} \left ( 3 \sum_{k > 4} v_k 
+ \sum_{\scriptstyle x  \in X \atop \scriptstyle  d(x) \not= 4} d(x) \right ) - 832n^2 |\partial S|. \eqno{(1)}$$

Let $v_4^E$ be the vertices lying in the interior of $E$,
each of which is 4-valent by construction. Let $v_4^{NE}$
denote the number of remaining 4-valent vertices.
Each of the vertices in the interior of $E$ contributes $1$ to
the area of $E$. So, by Lemma 7.2 and Theorem 7.3,
$$v_4^E \leq \pi (8000n^2 + 1)^2 \ \ell(\partial E).$$ 
So, $$v_4 = v_4^E + v_4^{NE} \leq 
\pi (8000n^2 + 1)^2 \ \ell(\partial E) + v_4^{NE}
\leq \pi (8000n^2 + 1)^2 (\ell(\partial E) + v_4^{NE}).
$$
Each fibre in $\partial E$ is adjacent to a disc tile that is not ordinary.
The number of such tiles is at most $1644n^2|\partial S|$, by Lemma 6.5. 
Each contributes at most $4$ to $\ell(\partial E)$. Each isolated point
of $\partial E$ is a generalised saddle with valence not equal to
$4$ or a vertex with valence not equal to $4$. Each 4-valent vertex
not in the interior of $E$ is adjacent to a disc tile that is not ordinary. This
tile contributes at most $2$ to $v_4^{NE}$. So, we
deduce that
$$\ell(\partial E) + v_4^{NE} \leq 9864 n^2 |\partial S| + \sum_{k \not= 4} v_k 
+ \sum_{\scriptstyle  x \in X \atop \scriptstyle  d(x) \not= 4} d(x).$$
Therefore,
$$\eqalign{
v_2 + v_3 &\geq \ell(\partial E) + v_4^{NE}
- \sum_{k > 4} v_k - 
\sum_{\scriptstyle  x \in X \atop \scriptstyle  d(x) \not= 4} d(x) - 9864 n^2 |\partial S|\cr
&\geq {v_4 \over \pi (8000n^2 + 1)^2} 
- \sum_{k > 4} v_k - 
\sum_{\scriptstyle  x \in X \atop \scriptstyle  d(x) \not= 4} d(x) - 9864 n^2 |\partial S|.} \eqno{(2)}$$
Adding 3 times (1) to (2), we deduce that
$$7 v_2 + 4v_3 >
{v_4 \over \pi(8000n^2 + 1)^2} + 2 \sum_{k > 4} v_k  - 12360 n^2 |\partial S|.$$
Therefore,
$$8 v_2 + 5 v_3 > 
{v_4 \over \pi (8000n^2 + 1)^2} + \sum_{k \not= 4} v_k  - 12360 n^2 |\partial S|
> {|V| \over \pi (8000n^2 + 1)^2} - 12360 n^2 |\partial S|.$$

The next stage is to discard vertices that are not deep.
By Lemma 6.8, the number of these is at most $3288 n^2 |\partial S|$.
Therefore, the number of deep 2-valent and 3-valent vertices in
$S$ is at least
$$\eqalign{
v_2 + v_3 - 3288 n^2 |\partial S| &\geq {8 v_2 + 5 v_3 \over 8} - (3288 n^2 |\partial S|) \cr
&\geq {w_\beta(S) \over 8 \pi (8000 n^2 + 1)^2} - 4833 n^2 |\partial S|\cr
&\geq {w_\beta(S) \over 2 \times 10^9 n^4} -  4833 n^2 |\partial S|,}$$
as required. $\square$

We are now in a position to prove Theorems 1.4 and 1.3, assuming Theorem 7.3.

\noindent {\sl Proof of Theorem 1.4.} Let $D$ be an arc presentation of
a split link $L$. Suppose that $D$ is not disconnected. Let $n$ be its arc index.
Let ${\cal T}$ denote Dynnikov's triangulation of
$S^3$, given in Section 6.1. By Theorem 5.1, 
there is a splitting 2-sphere for $S^3 - L$ which is in normal form
with respect to ${\cal T}$. 
By Theorem 5.5, there is such a sphere $S$ which is a vertex surface with respect to ${\cal T}$,
for which $(w_\beta(S), w(S))$ is minimal.
By Theorem 5.9, the binding weight $w_\beta(S)$ of this surface is at most $n 2^{7n^2}$,
which is less than $2^{8n^2}$.
The surface $S$ inherits a singular foliation.
Then by Theorem 7.1, the number
of 2-valent and 3-valent vertices in $S$ is at least
$${w_\beta(S) \over 2 \times 10^9 n^4}.$$
By Lemma 6.7, these vertices come in at most $48n^2$ types.
So, there is a collection of at least $w_\beta(S) / (10^{11} n^6)$ 2-valent
or 3-valent vertices, all of the same type. Applying Proposition 6.11, there is a sequence of
at most $n/2$ cyclic permutations, at most $n^2/4$ exchange moves,
a generalised exchange move and some ambient isotopies, which reduces
the binding weight of the surface by at least $2w_\beta(S) / (10^{11} n^6)$.

After we have performed these moves, the result is
a new arc presentation of $L$. This gives a new triangulation, which we
will call ${\cal T}'$. Now, the new surface $S'$ need not be normal, but by Theorem 5.1, there is
reducing 2-sphere which is normal with respect to ${\cal T}'$ and
with no greater binding weight.
We may therefore repeat the above argument with this new arc
presentation, triangulation and splitting sphere.

Let $x$ be the number of these steps required to reduce the complexity
down to less than 1, by which time we must have reach a disconnected arc presentation,
as required. Then
$$2^{8 n^2}\left ( 1 - {2 \over 10^{11} n^6} \right )^{x-1} \geq 1,$$
because after $x-1$ steps, the binding weight is still at least $1$, by the definition of $x$.
Taking logs, we obtain
$$\left ({x}-1 \right ) \log \left ( 1 - {2 \over 10^{11} n^6} \right ) + 8 n^2 \log 2
\geq 0.$$
Now, $\log (1 - y) \leq - y$ for any $y$ between $0$ and $1$, and so
$$(x - 1) \leq (8 n^2 \log 2) (10^{11} n^6/2).$$
Therefore $x \leq 3 \times 10^{11} n^8$, which proves the theorem. $\square$

\noindent {\sl Proof of Theorem 1.3.}
We now consider the case where $L$ is the unknot. The argument is similar
to that of Theorem 1.4, but it is made more complicated by the
presence of boundary. A flowchart for the proof is shown in Figure 25.

\vskip 18pt
\centerline{
\includegraphics[width=2.8in]{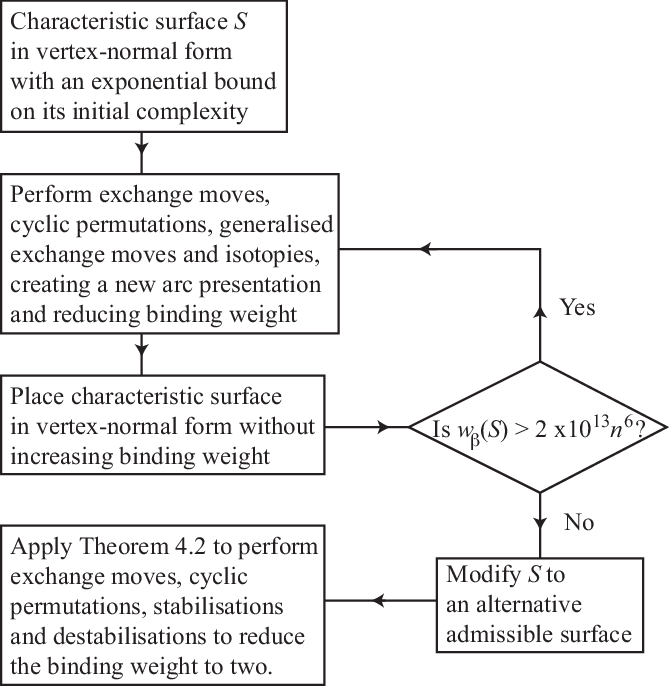}
}
\vskip 6pt
\centerline{Figure 25: Flowchart for the proof of Theorem 1.3}

Let $M$ be the exterior of $L$, and let ${\cal P}$ be the polyhedral structure 
for $M$ defined in Section 6.2. Let $C$ be the specified longitude. 
Now apply Theorem 5.2 to find a compression disc $S$ for $\partial M$
which is normal and with boundary equal to $C$.
By Theorem 5.6, we may choose $S$ so that it is a boundary-vertex surface
and so that $(w_\beta(S), w(S))$ is minimal.
So, by Theorem 5.10 and Lemma 6.1, its binding weight $w_\beta(S)$ is at most 
$$(2n) (2 \times 10^6)^{4 \times 10^6 n^2}(8n^2) < 2^{10^8 n^2}.$$

Now suppose that $w_\beta (S) > 2 \times 10^{13} n^6$. Then by Theorem 7.1, the number
of deep 2-valent and 3-valent vertices in $S$ is at least
$${w_\beta (S) \over 2 \times 10^9 n^4} - 4833n^2 \geq {w_\beta (S) \over 4 \times 10^{9} n^4}.$$
By Lemma 6.7, these vertices come in at most $48n^2$ types.
So, there is a collection of at least $w_\beta(S) / (2 \times 10^{11} n^6)$ deep 2-valent
or 3-valent vertices, all of the same type. By Proposition 6.11, we may perform at most
$n/2$ cyclic permutations, at most $n^2/4$ exchange moves, a generalised
exchange move and some ambient isotopies, to reduce the binding weight
by at least $w_\beta(S) / (10^{11} n^6)$. 

As in the proof of Theorem 1.4, this creates a new arc presentation of $L$.
This then gives a new polyhedral decomposition ${\cal P}'$. Let $S'$ be the result of $S$
after making this modification. Then, $S'$ is a surface properly embedded in the
exterior of the new copy of $L$. Its boundary remains a longitude on $\partial N(L)$. 
Moreover, the decomposition of $L$ into `up' and `down' arcs, plus possibly
one extra arc, is preserved. And near these arcs, $S'$ continues to lie in the
up and down directions from $L$. We may therefore
isotope $\partial S'$, taking it to the new specified longitude for $L$, without
changing its binding weight. By Theorem 5.6, there exists a normal disc $S''$ in ${\cal P}'$ with boundary
this specified longitude that is a boundary-vertex surface,
such that $w_\beta(S'') \leq w_\beta(S')$. We choose
$S''$ so that $(w_\beta(S''), w(S''))$ is minimal. We then repeat the above argument.

Let $x$ be the number of steps required to reduce the complexity
down to at most $2 \times 10^{13} n^6$. By the above argument, 
$$x \leq (10^8n^2 \log 2) (10^{11} n^6) + 1 \leq 10^{19} \, n^8.$$
This is at most $3 \times 10^{18} n^{10}$ exchange moves, at most $5 \times 10^{18} n^{9}$ 
cyclic permutations and at most $10^{19} n^8$
generalised exchange moves.
Once we have reduced the binding weight below $2 \times 10^{13} n^6$, we apply Proposition 6.12
to create a characteristic surface in admissible form with binding weight at most
$2 \times 10^{13}n^6 + n$. If it has a winding vertex, its winding angle is at most 
$2 \pi$ times the twisting number of the specified longitude. Now, the exchange moves,
cyclic permutations and generalised exchange moves that we have performed so far
do not affect the writhe of the rectangular diagram, which therefore remains $k$. So, 
as explained in Section 6.4, the twisting number of the specified longitude is at most $|k|+n+1$.
So, by Lemma 4.1, there is a sequence of at most $|k|+n < n^2$ stabilisations, less than $n^2(n+n^2)$ 
exchange moves and an ambient
isotopy of the knot complement taking the characteristic surface into alternative
admissible form with binding weight at most $2 \times 10^{13}n^6 + n + n^2$.
Then, by Theorem 4.2, there is a sequence of at most 
$4 (n + n^2)^2 (2 \times 10^{13}n^6 + n + n^2)$ exchange moves, at most 
$(n + n^2) (2 \times 10^{13}n^6 + n + n^2)$ cyclic permutations, at most 
$(2 \times 10^{13}n^6 + n + n^2)$ stabilisations and at most 
$(2 \times 10^{13}n^6 + n + n^2)$ destabilisations
taking $D$ to the trivial arc presentation.
So, the total number of exchange moves is at most
$$(3 \times 10^{18}) n^{10} + n^2(n + n^2) +
4 (n + n^2)^2 (2 \times 10^{13}n^6 + n + n^2) \leq 4 \times 10^{18} n^{10}.$$
The total number of cyclic permutations is at most
$$(5 \times 10^{18}) n^{9} + (n + n^2) (2 \times 10^{13}n^6 + n + n^2) \leq 6 \times 10^{18} n^9.$$
The number of stabilisations and destabilisations are each at most
$$n^2 + (2 \times 10^{13}n^6 + n + n^2) \leq 3 \times 10^{13} n^6,$$
as required. $\square$

\vskip 18pt
\centerline{\caps 8. Branched surfaces}
\vskip 6pt

The remainder of this paper is devoted to the proof of Theorem 7.3.
The proof will be given in Section 9, but it requires some background
theory on branched surfaces, which we recall in this section.
This is mostly standard material, which can be found in [9], for example.

\vskip 6pt
\noindent {\caps 8.1. Definitions}
\vskip 6pt

A {\sl branched surface} is a compact 2-complex $B$ smoothly embedded in
a 3-manifold $M$, with the following properties. At each point $x$ of $B$,
there is a specified tangent plane in $T_x(M)$ and all the 1-cells and 2-cells
that contain $x$ have tangent spaces at $x$ that lie in this tangent plane.
This tangent plane is denoted by $T_x(B)$.
Thus, at each point $x$ in the interior of a 1-cell of $B$, $T_x(B)$
is divided into two half-planes by the tangent space of
the 1-cell. We term these the two {\sl sides} at $x$.
We require that, at each such point $x$, either there are 2-cells on both sides of $x$
or the 1-cell is incident to a single 2-cell. The closure of the union of the points $x$ of the
former type is the {\sl branching locus} of $B$. The closure of the union of the points $x$ of
the latter type is the {\sl boundary} of $B$, which we denote by $\partial B$.
(Note that we do not require $\partial B$ to lie in $\partial M$.)
The 2-cells of $B$ are called the {\sl patches} of $B$.

Note that this definition is somewhat more general than the one that is
frequently used, for example in [9]. There, a branched surface is defined via
its possible local models. In our definition, it is not the case that there
are only finitely many local models. An example is shown in Figure 26, but
this is not the general situation.

\vskip 18pt
\centerline{
\includegraphics[width=2in]{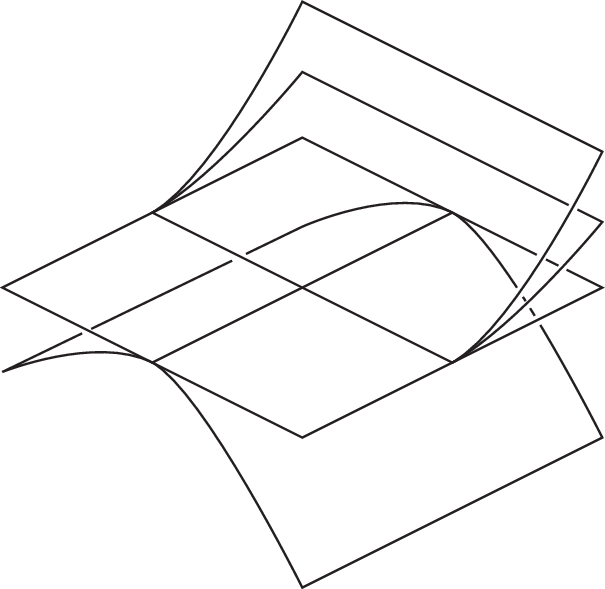}
}
\vskip 6pt
\centerline{Figure 26: A branched surface}

\vfill\eject
A thickening $N(B)$ of $B$ has a decomposition as a union of 
{\sl fibres}, each of which is homeomorphic to an interval. (This thickening is almost
a regular neighbourhood, except that $\partial B \subset \partial N(B)$.) Away from a small
regular neighbourhood of the 1-skeleton of $B$, this is just an $I$-bundle. There
is a map $\pi \colon N(B) \rightarrow B$ which collapses
each fibre to a point. For each $x \in B$, the fibre through $x$
is required to have tangent space that is complementary to $T_x(B)$.
Also, each fibre is required to intersect $\partial N(B)$ in its endpoints,
plus possibly a finite collection of closed intervals (see Figure 27).
The {\sl horizontal boundary} $\partial_h N(B)$ is the union of the endpoints
of these fibres. The {\sl vertical boundary} $\partial_v N(B)$ is
${\rm cl}(\partial N(B) - \partial_h N(B))$. Each component of
${\rm cl}(\partial_vN(B) - \pi^{-1}(\partial B))$ is termed a {\sl cusp}.

\vskip 18pt
\centerline{
\includegraphics[width=1.7in]{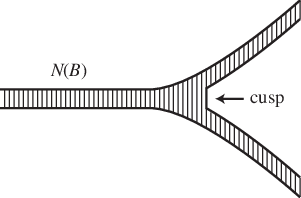}
}
\vskip 6pt
\centerline{Figure 27: The fibred neighbourhood of $B$}

Due to the potential interaction between $\partial B$ and the branching locus,
these cusps can have slightly complicated topology. However,
in the following case, they are rather simple.

\noindent {\bf Lemma 8.1.} {\sl Let $B$ be a branched surface in an orientable 3-manifold $M$.
Suppose that $B \cap \partial M = \partial B$ and that $\pi^{-1}(\partial B) = N(B) \cap \partial M$.
Suppose also that at each $x \in \partial B$, the tangent plane $T_x(B)$ does not equal
$T_x(\partial M)$.
Then each cusp either is an annulus or is a disc $D$ such that $D \cap \partial M$ is two
disjoint arcs in $\partial D$.}

\noindent {\sl Proof.} The cusps lie in a regular neighbourhood of the 1-skeleton of $B$.
Near each 1-cell of $B$, they have a simple form. They have the structure of $I$-bundles
over this 1-cell, where each $I$-fibre is the intersection between the cusp and a fibre
of $N(B)$. We need to analyse how these $I$-bundles join together near the 0-cells of $B$. 
We claim they patch together to give each cusp the structure of an $I$-bundle.

Let $v$ be a 0-cell of $B$ in the interior of $M$. Pick a small smoothly embedded disc $P$ running
through $v$ with tangent plane at $v$ equal to $T_v(B)$. Let $N$ be a
thickening of $P$, which is an $I$-bundle over $P$, and let $A$ be the
$I$-bundle over $\partial P$, which is an annulus in $\partial N$. Then we may
arrange that $B \cap \partial N$ lies in the interior of $A$, and that it is transverse
to the $I$-fibres of $A$. It is a branched 1-manifold
in $A$. By our hypothesis on $\partial B$, this branched 1-manifold has no boundary.
Hence, it divides $A$ into a collection of bigons, together with at least two annuli with smooth boundary.
At each such bigon, two cusps of $N(B)$ enter and are joined together.
We deduce that the $I$-bundle structures do indeed patch together correctly here.

A similar analysis applies near each 0-cell $v$ of $B$ that lies in $\partial M$.
Again, pick a small smoothly embedded disc $P$ that contains $v$,
and with tangent plane at $v$ that equals $T_v(B)$. Then $P$ may be
chosen so that $P \cap \partial M$ is a single arc in $\partial P$,
which contains $v$ in its interior. Thicken $P$ to an $I$-bundle $N$ over $\partial P$.
Let $W$ be the $I$-bundle over ${\rm cl}(\partial P - \partial M)$,
which is a disc. Again, we may assume that $B \cap W$ is transverse to the fibres.
It therefore divides $W$ into a collection of bigons, at least two smooth discs which include
collars on the horizontal boundary of $W$, and some discs, each
of which has a single cusp in its boundary and a single arc of intersection
with $\partial W$. At this latter type of disc, a cusp of $N(B)$ hits
$\partial M$ and terminates. At each bigon, two bits of cusp of $N(B)$ join together.
Thus, again, we deduce that $I$-bundle structures patch together as required,
proving the claim.

Since each cusp lies in the
boundary of the orientable 3-manifold $N(B)$, we deduce that it is orientable and
hence an annulus or disc, as required. $\square$

\vskip 6pt
\noindent {\caps 8.2. Surfaces carried by branched surfaces}
\vskip 6pt

A compact surface $S$ is {\sl carried} by a branched surface $B$ if $S$ is embedded in $N(B)$,
it is transverse to the fibres and $S \cap \pi^{-1}(\partial B) = \partial S$.
These conditions ensure that, for each patch of $B$, the cardinality of
$\pi^{-1}(x) \cap S$ is constant for all $x$ in the interior of that patch. This cardinality
is termed the {\sl weight} of $S$ in that patch.
These weights form a collection of non-negative
integers, which is known as the {\sl vector} associated with $S$,
and is denoted by $[S]$.
The weights satisfy a system of linear equations, which are known 
as the {\sl matching equations}. These specify that, at each 1-cell
in the branching locus of $B$, the total weight of the patches on one side
is equal to the total weight of the patches on the other.
Conversely, given a solution to these matching equations by non-negative
integers, one can form a compact surface carried by $B$ with these
weights.

\vskip 6pt
\noindent {\caps 8.3. Summation of surfaces}
\vskip 6pt

Let $S$, $S_1$ and $S_2$ be surfaces carried by $B$. Then
$S$ is said to be the {\sl sum} of $S_1$ and $S_2$ if $[S] = [S_1] + [S_2]$.
We say that $S_1$ and $S_2$ are {\sl summands} of $S$.

There is an alternative way of viewing summands of a surface.

\noindent {\bf Lemma 8.2.} {\sl Let $S$ and $S_1$ be surfaces carried
by a branched surface $B$. Then $S_1$ is a summand of $S$ if and only if,
in every patch of $B$, the weight of $S$ is at least the weight
of $S_1$.}

\noindent {\sl Proof.} Suppose that $S$ and $S_1$ satisfy this weight condition.
Consider the vector $[S] - [S_1]$. Since $[S]$
and $[S_1]$ satisfy the matching equations, so does $[S] - [S_1]$.
By assumption, each of its co-ordinates is non-negative.
Hence, it corresponds to a surface $S_2$ carried by $B$,
and $S$ is the sum of $S_1$ and $S_2$.

Conversely, if $S$ is the sum of $S_1$ and $S_2$, then clearly,
the weight of $S$ is at least the weight of $S_1$ in each patch.
$\square$

\vskip 6pt
\noindent {\caps 8.4. Branched surfaces associated to normal surfaces}
\vskip 6pt

Let ${\cal P}$ be a polyhedral decomposition of a compact 3-manifold $M$. Associated to
any normal, properly embedded surface $S$, there is a branched surface $B_S$,
which we term a {\sl normal} branched surface. It carries $S$.

It is constructed as follows. For each type of elementary normal disc in $S$,
we take one such disc. We arrange for these discs to be smoothly embedded.
They form the patches of $B_S$.
For each face $F$ of ${\cal P}$ with polyhedra on both sides
and for each arc type of $F \cap S$, we glue all patches of $B_S$
which contain this arc type along this arc.

This is a branched surface, because the tangent planes to
$B_S$ can be defined as follows. For each $x \in B_S$ lying in a 1-cell of ${\cal P}$,
pick a tangent plane $T_x(B_S)$ that does not contain the tangent plane of
the 1-cell. For each $x \in B_S$ lying in the interior of a face of ${\cal P}$,
pick a tangent plane not equal to the tangent plane of the face. We can do this
compatibly with the choices for the points in the 1-cells of ${\cal P}$. For each point $x$ inside
the interior of a patch, we define $T_x(B_S)$ to be the tangent
plane of the elementary normal disc in which it lies. Although no Riemannian metric has been
specified, one should still think of the tangent planes of $B_S$ at the
1-cells and 2-cells of ${\cal P}$ as being `orthogonal' to those cells.

Note that $S$ is carried by $B_S$. For we may take a regular
neighbourhood $N(S)$ of $S$, such that $N(S)$ intersects each
polyhedron in a union of elementary normal discs. Then, when
two elementary normal discs of $S$ are normally parallel, we attach
the space between them to $N(S)$. Also, when two arcs of $S$ in a face
of ${\cal P}$ are normally parallel, we attach a slight thickening of the
space between them to $N(S)$. The resulting space is
a 3-dimensional subset of the 3-manifold, which we term $N(B_S)$.
It is composed of a collection of regions, each of which
is the product of an elementary normal disc type of $S$ with an interval.
There is therefore a map $\pi \colon N(B_S) \rightarrow B_S$
which collapses these intervals to points. It is
clear that $N(B_S)$ is a fibred regular neighbourhood of $B_S$.
By construction, $S$ is a subset of $N(B_S)$ that is transverse
to the fibres. The boundary of $B_S$ is precisely $B_S \cap \partial M$.
So, $\pi^{-1}(\partial B_S) \cap S = S \cap \partial M = \partial S$.

When a surface $S'$ is carried by $B_S$, it is normal with respect to
${\cal P}$. Moreover, the vector for $S'$ as a surface carried by $B_S$
is equal to its normal surface vector. Hence, summation of surfaces in the branched surface $B_S$
corresponds to the summation of normal surfaces. More precisely,
suppose that $S'$, $S_1$ and $S_2$ are surfaces carried
by $B_S$ such that $S'$ is the sum of $S_1$ and $S_2$.
Then these are normal and $S' = S_1 + S_2$ as normal surfaces.
Conversely, if $S$, $S_1$ and $S_2$ are normal surfaces satisfying
$S = S_1 + S_2$, then, when one forms the branched surface $B_S$ starting
from the normal surface $S$, then $S$, $S_1$ and $S_2$ are
all carried by $B_S$ and $S$ is the sum of $S_1$ and $S_2$ in
$B_S$.

\vskip 6pt
\noindent {\caps 8.5. Branched surfaces carried by branched surfaces}
\vskip 6pt

We say that a branched surface $B_1$ is {\sl carried} by a branched surface
$B_2$ if 
\item{(1)} $B_1$ is smoothly embedded in $N(B_2)$, and
\item{(2)} for each point $x$ in $B_1$, $T_xB_1$ is transverse to the
fibres of $N(B_2)$.

We do not require that $\partial B_1$ lies in $\pi^{-1}(\partial B_2)$,
where $\pi \colon N(B_2) \rightarrow B_2$ is the collapsing map for $B_2$.
In fact, $\pi(\partial B_1)$ is permitted to run through the interior
of patches of $B_2$.

\noindent {\bf Lemma 8.3.} {\sl If $B_1$ is carried by $B_2$, then any closed surface
carried by $B_1$ is also carried by $B_2$.}

\noindent {\sl Proof.} Let $S$ be a closed surface carried by $B_1$.
We may assume that, at each point of $x$ of $S$, $T_x(S)$ is arbitrarily
close to $T_{\pi(x)}(B_1)$, where $\pi \colon N(B_1) \rightarrow B_1$ is
the collapsing map for $B_1$. Since $T_{\pi(x)}(B_1)$ is transverse to
the fibre of $N(B_2)$ through $\pi(x)$, we can therefore arrange
that $T_x(S)$ is also transverse to the fibre of $N(B_2)$ through $x$.
For a closed surface, this is the definition of $S$ being carried by $B_2$.
$\square$

\vfill\eject
\centerline {\caps 9. Euclidean subsurfaces of the characteristic surface}
\vskip 6pt

This section is devoted to the proof of Theorem 7.3. 

\noindent {\bf Theorem 7.3.} {\sl Let $L$ be the unknot or a split link. Fix
an arc presentation of $L$ with arc index $n$.
Let $S$ be a characteristic surface in PL-admissible normal form
with respect to the polyhedral decomposition ${\cal P}$,
and that is a boundary-vertex surface. In the case where $L$ is the unknot,
suppose that $\partial S$ is the specified longitude. Suppose also that $(w_\beta(S), w(S))$ is minimal
among all characteristic surfaces for $L$ with the same boundary as $S$.
Let $E$ be the designated Euclidean subsurface
of $S$. Then each point of $E$ has distance at most $8000 n^2$ from $\partial E$.}

For a subset $F$ of a metric space, and a positive real number $r$, let $N_r(F)$ denote
the set of points with distance at most $r$ from $F$.
Thus, in our situation, Theorem 7.3 asserts that in the metric space $E$, $N_{8000 n^2}(\partial E)$
is all of $E$.

\vskip 6pt
\noindent {\caps 9.1. Overview of the proof}
\vskip 6pt

The strategy for the proof is as follows. Suppose that there is a point
in $E$ with distance more than $8000n^2$ from $\partial E$. Then,
around this point, there is a large Euclidean region. We will show
that this implies that there is a normal torus which is a summand
for some multiple of $S$. This will imply that $S$ is not a boundary-vertex
surface, which is contrary to hypothesis.

Throughout this section, $S$ will be a characteristic surface for $L$,
which is in normal form with respect to the polyhedral decomposition,
and that is a boundary-vertex surface. Also, $(w_\beta(S), w(S))$ is minimal
among all characteristic surfaces for $L$ with the same boundary as $S$.
As above, $E$ will denote the designated Euclidean subsurface of $S$.
We will prove Theorem 7.3 by contradiction, and therefore suppose
that there is some point $z$ in $E$ with distance more than $8000n^2$
from $\partial E$. Let $E'$ denote the component of $E$
containing $z$.

\vskip 6pt
\noindent {\caps 9.2. A branched surface carrying the Euclidean subsurface}
\vskip 6pt

Starting with the designated Euclidean subsurface $E$, we can form a branched
surface $B$ as follows. By construction, $E$ is a union of
square-shaped tiles. We say that two tiles of $E$
are {\sl normally parallel} if they are normally parallel in ${\cal P}$.
We first form a 2-complex $\overline B$, where
each 2-cell of $\overline B$ arises from a normal equivalence
class of tiles of $E$. Each 2-cell therefore has the shape of a square tile.
We call each of these 2-cells a {\sl Euclidean patch}. When  two tiles of $E$
are incident along a separatrix, then we glue the associated patches
of $\overline B$ along the corresponding edges.

However, $\overline B$ is not quite a branched surface, because there may be 1-cells
of $\overline B$ with more than one 2-cell on one side but no 2-cells on the other.
To remedy this, we attach to $\overline B$ some extra 2-cells, as follows.
For each normal equivalence class of tile $T$ of $S$ that is incident
to $E$ but not a subset of $E$, we attach a thin neighbourhood of $T \cap \partial E$.
This is a collection of thin discs. The boundary of each such disc
consists of two long arcs and two short arcs. When two such discs are incident
because they share a common isotopy class of short arc in their boundary,
we glue these discs together along this arc, forming new branching locus on their
boundary. The result is the branched surface $B$.

Note that we make no identifications on the long arcs that do not lie in $\partial E$.
Instead, they become part of the boundary of $B$.

Note also there is a retraction map $B \rightarrow \overline B$, which 
collapses the thin discs attached to $\overline B$. This is a homotopy
equivalence.

\noindent {\bf Lemma 9.1.} {\sl The number of Euclidean patches of $B$ is
at most $24n^2$.}

\noindent {\sl Proof.} Each Euclidean patch of $B$ corresponds to a normal
isotopy class of ordinary tiles. These tiles are all normally parallel,
and so there are two that are outermost, $T_1$ and $T_2$
say. We claim that each of $T_1$ and $T_2$ must contain an elementary normal
triangle, square or half-square that is outermost in $N(B_S)$. Suppose that
this is not the case, for $T_1$, say. Then, adjacent to each elementary normal
triangle, square or half-square in $T_1$ on both sides of $T_1$, there is another elementary
normal triangle, square or half-square
and the union of these forms two tiles which are normally parallel
to $T_1$ on both sides of $T_1$. These are ordinary tiles, which
contradicts the assumption that $T_1$ is outermost in $N(B)$. This proves the claim.

As a consequence of the claim, the number of Euclidean patches of $B$ is
at most the number of normal isotopy classes of triangles, squares and 
half-squares in $S$. There are at most 6 of these in each truncated tetrahedron.
There are at most $4n^2$ truncated tetrahedra in the polyhedral decomposition.
This proves the lemma. $\square$

\noindent {\bf Lemma 9.2.} {\sl $B$ is carried by $B_S$.}

\noindent {\sl Proof.} We need to find an embedding $B \rightarrow N(B_S)$.
Each Euclidean patch of $B$ is a normal equivalence
class of Euclidean tiles. The remaining patches are subsets of tiles.
Each tile is made up pieces of elementary normal discs.
Embed $B$ into $N(B_S)$ by including each such piece into
a regular neighbourhood of the relevant patch of $B_S$.
It is easy to see that this inclusion map has the right properties. $\square$

Let $B'$ be the component of $B$ such that $N(B')$ contains $E'$. Let $\overline B'$
be the component of $\overline B$ such that $N(\overline B')$ contains $E'$. Then the
retraction map $B' \rightarrow \overline B'$ is a homotopy equivalence,
using which we may identify $\pi_1(B')$ and $\pi_1(\overline B')$.

Note that $E'$ is not carried by $B'$, because we added non-Euclidean
patches to $\overline B'$. But we can take a small regular neighbourhood
of $E'$ in $S$, denoted $\hat E$, so that $\hat E$ is carried by $B'$.

\vskip 6pt
\noindent {\caps 9.3. Reducing to the case of trivial monodromy}
\vskip 6pt

We now define a homomorphism $\mu \colon \pi_1(B') \rightarrow O(2)$,
where $O(2)$ is the group of orthogonal transformations of ${\Bbb R}^2$.
We term this the {\sl monodromy} of the branched surface $B'$.

It is convenient to subdivide the cell structure on $\overline B'$,
introducing a new vertex into the midpoint of each 1-cell of $\overline B'$,
and introducing a new vertex in the centre of each 2-cell and coning
off from this vertex. Let $b$ be a basepoint for $\overline B'$,
which is a vertex at the centre of one of the original 2-cells.

Around each vertex of this new cell structure,
we pick a Euclidean disc of radius $1/4$, say,
that lies in $B'$. Given two vertices which are the endpoints
of a 1-cell of $\overline B'$, there is a canonical isometry taking one disc to the
other, which is Euclidean translation along the 1-cell. If one
follows a loop that encircles a 2-cell of $\overline B'$, the composition of these
Euclidean isometries is the identity. Thus, one may define
$\mu \colon \pi_1(\overline B') \rightarrow O(2)$ as follows. Given a
cellular loop $\ell$ in $\overline B'$ based at $b$, it is a composition of paths along
1-cells, and this then gives a composition of Euclidean isometries.
This composition is a Euclidean isometry which takes
the disc neighbourhood of $b$ to itself. It is therefore
an element of the orthogonal group $O(2)$. Because the monodromy
around each 2-cell is trivial, this gives a well defined homomorphism
$\mu \colon \pi_1(\overline B') \rightarrow O(2)$, and hence a homomorphism
$\mu \colon \pi_1(B') \rightarrow O(2)$. Note that $\mu(\ell)$
is an isometry that preserves the tile containing $b$, and hence
the image of $\mu$ lies in a subgroup of $O(2)$ of order 8.

We now define a finite-sheeted cover $\tilde B$ of $B'$, as follows.
We let $\tilde B$ be the
covering space of $B'$ corresponding to 
kernel of the monodromy homomorphism $\mu \colon \pi_1(B') \rightarrow O(2)$.

We record some properties of $\tilde B$.

\noindent {\bf Property 9.3.} {\sl $\tilde B$ is a branched surface.}

\noindent {\sl Proof.} There is an inclusion of $B'$ into the 3-manifold $N(B')$,
which  is a homotopy equivalence. Hence, associated with the kernel
of $\mu \colon \pi_1(B') \rightarrow O(2)$, there is a covering
space of $N(B')$, which we denote by $N(\tilde B)$. This
is a regular neighbourhood of $\tilde B$, and hence
is the required 3-manifold. Note that there is a
collapsing map $\pi \colon N(\tilde B) \rightarrow \tilde B$.
$\square$

\noindent {\bf Property 9.4.} {\sl $\tilde B$ has trivial monodromy.}

\noindent {\sl Proof.} Implicit in this statement is the assertion that one can define a monodromy
homomorphism $\tilde \mu \colon \pi_1(\tilde B) \rightarrow O(2)$. 
But the method of doing this is by direct analogy with
the case of $B'$. By construction, the monodromy homomorphism
of $\tilde B$ has trivial image. $\square$

\noindent {\bf Property 9.5.} {\sl $\tilde B$ is transversely
orientable.}

\noindent {\sl Proof.} The obstruction to finding a transverse orientation
to a branched surface is the existence of a closed loop $\ell$ in the branched surface, 
so that as one travels
around this loop, and one keeps track of a transverse orientation,
this is reversed by the time one returns to the starting point.
This is evident in the monodromy homomorphism $\tilde \mu$.
For then $\tilde \mu (\ell)$ has non-trivial image after composing
with the determinant homomorphism $O(2) \rightarrow
\{ \pm 1 \}$. This contradicts the fact that $\tilde \mu$ has trivial
image. $\square$

\noindent {\bf Property 9.6.} {\sl The number of Euclidean patches of $\tilde B$ is at most $192 n^2$.}

\noindent {\sl Proof.} This follows from the fact that the number of Euclidean patches of $B$ is at most $24n^2$.
$\square$

\noindent {\bf Property 9.7.} {\sl The total length
of the intersection between the singular locus of $\tilde B$
and the Euclidean patches of $\tilde B$ is at most $768n^2$.}

\noindent {\sl Proof.} The singular locus is a subset of the
1-skeleton of $\tilde B$. Since each Euclidean patch of $\tilde B$ is isometric
to a Euclidean square of side length 1, the total length
of the singular locus incident to the Euclidean patches
is at most 4 times the number of Euclidean patches.
$\square$

There is an inclusion $i \colon \hat E \rightarrow N(B')$
and a collapsing map $\pi \colon N(B') \rightarrow B'$. The kernel of
$\mu \pi_\ast i_\ast$ is a finite index subgroup of $\pi_1(\hat E)$.
Let $\tilde E$ be the corresponding covering space of $\hat E$.
Then $\tilde E$ is carried by $\tilde B$.

The actual result we will prove in this section is as follows.

\noindent {\bf Proposition 9.8.} {\sl Suppose that there is a point
$x$ in $\tilde E$ with distance more than $8000n^2$
from $\partial \tilde E$. Then $\tilde E$ has a torus summand,
when viewed as a surface carried by the branched surface $\tilde B$.}

We now show how Theorem 7.3 follows from this.

We are supposing that there is a point $z$ in $E'$ with distance more than $8000 n^2$ from
$\partial E$. Let $\tilde z$ be a point in the inverse image of $z$ in $\tilde E$.
Then this  has distance more than $8000 n^2$ from $\partial \tilde E$. 
This is because a path from $\tilde z$ to $\partial \tilde E$ projects to a path from
$z$ to $\partial E$ with the same length. So, applying Proposition 9.8,
we deduce that $\tilde E$ has a torus summand $T$. 

Now, the covering map $\tilde B \rightarrow B'$ sends $[\tilde E]$
to a non-zero multiple $m[\hat E]$ of $[\hat E]$. The covering map sends $[T]$ to a vector satisfying the
matching equations for $B'$ and with zero boundary. This corresponds to a closed surface $T'$
carried by $B'$. By Lemma 8.2, 
$T'$ is a summand of $m [\hat E]$. Since $T'$ is a closed
surface, Lemmas 8.3 and 9.2 imply that
$T'$ is also carried by $B_S$. So, $[T']$ is a summand of
$m [S]$. We deduce that $S$ is not a boundary-vertex surface. But this is
contrary to the hypothesis of Theorem 7.3.

Thus, Proposition 9.8 implies Theorem 7.3. We therefore now work almost
exclusively with $\tilde E$ and $\tilde B$.

We fix a transverse orientation of $\tilde B$. This induces a transverse orientation of $\tilde E$.

\vskip 6pt
\noindent {\caps 9.4. Grids and annuli}
\vskip 6pt

The proof now divides into two cases. Either there is a closed
geodesic in $\tilde E - N_{1000 n^2}(\partial \tilde E)$ with length at most 
$12000n^2$, or there is not. 

Suppose first that there is such a closed geodesic. Since $\tilde E$ has trivial monodromy,
this closed geodesic is a multiple of a simple closed geodesic $\alpha$.

Now, $\alpha$ represents a non-trivial element of $\pi_1(\tilde E)$.
Let $\tilde E_\infty$ be the universal cover of $\tilde E$, and let $\tilde \alpha$
be one component of the inverse image of $\alpha$ in $\tilde E_\infty$. 
Then corresponding to $\alpha$, there is a covering transformation $\tau$ of $\tilde E_\infty$.
For each point $\tilde x_\infty$ on $\tilde \alpha$, $\tau$ acts on $\tilde x_\infty$ by translation
along $\tilde \alpha$. This is also true of points close to $\tilde x_\infty$.
Now if $N_r(\tilde \alpha)$ is disjoint from $\partial \tilde E_\infty$, for some $r > 0$, then $N_r(\tilde \alpha)$
is isometric to $[-r,r] \times \tilde \alpha$. Hence, $N_r(\tilde \alpha) / \langle \tau \rangle$
is isometric to $[-r,r] \times \alpha$. The covering map $\tilde E_\infty \rightarrow \tilde E$
sends $N_r(\tilde \alpha)$ onto $N_r(\alpha)$. If there are two points of
$N_r(\tilde \alpha)$ that do not differ by an element of $\langle \tau \rangle$
but which are sent to the same point in $N_r(\alpha)$, then $\tilde E$
is a torus. This is impossible because $\partial \tilde E$ is non-empty. We therefore deduce that
$N_r(\alpha)$ is isometric to a Euclidean annulus. Summarising,
we have proved the following.

\noindent {\bf Lemma 9.9.} {\sl Suppose that $\alpha$ is a simple closed
geodesic in $\tilde E - N_{1000 n^2}(\partial \tilde E)$ with length at most
$12000n^2$. Then, for all $r \leq 1000 n^2$, $N_r(\alpha)$ is a Euclidean
annulus with core curve $\alpha$.}

Suppose now there is no closed
geodesic in $\tilde E - N_{1000 n^2}(\partial \tilde E)$ with length at most 
$12000n^2$. Then, the exponential 
map based at any $x \in \tilde E - N_{7000 n^2}(\partial \tilde E)$ defines an isometry between a
Euclidean disc of radius $6000n^2$ and $N_{6000 n^2}(x)$.
Let $x$ lie at the centre of a tile.
Hence, centred at $x$, there is a {\sl grid}, which is a union of square tiles  that is
isometric to Euclidean square. (See Figure 28.) We may find such grids with
any odd integer side length less than $6000 \sqrt{2} n^2$.
Note that $8000 n^2 +1 \leq 6000 \sqrt{2} n^2$.

For a positive integer $r \leq 4000 n^2$ and any $x \in \tilde E - N_{7000n^2}(\partial\tilde E)$, 
we let $D(x,r)$ denote a grid centred at a tile containing $x$
with side length $2r+1$. Note that when $x$ lies in more than
one tile, this is slightly ambiguous, but this ambiguity will not cause
any problems.

\vskip 18pt
\centerline{
\includegraphics[width=1.3in]{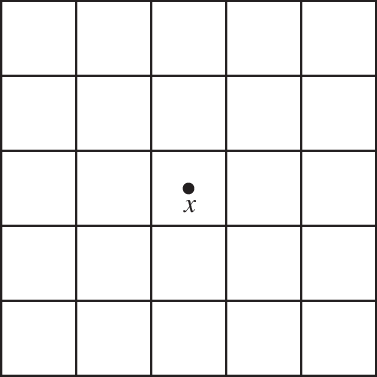}
}
\vskip 6pt
\centerline{Figure 28: A grid centred at $x$}

The proof now divides into these two cases. We focus first on the case
where each $x$ in $\tilde E - N_{7000 n^2}(\partial \tilde E)$ lies in the central tile
of a grid with side length $8000 n^2 + 1$.

\vskip 6pt
\noindent {\caps 9.5. First-return maps}
\vskip 6pt

For a transversely oriented subsurface $F$ of $\tilde E$, we now define the {\sl first-return map}.
This is function (which need not be continuous) $r_F \colon {\rm dom}(r_F)
\rightarrow F$, where the domain of definition ${\rm dom}(r_F)$ is a
subsurface of $F$.
For each point $x$ in $F$, there is a fibre $I_x$ in $N(\tilde B)$ through $x$. This fibre
is divided into two by $x$. Let $\alpha_x$ be the component of $I_x - \{ x \}$
into which $F$ points at $x$. Define ${\rm dom}(r_F)$ to be those $x \in F$
such that $\alpha_x \cap F \not= \emptyset$. For $x \in {\rm dom}(r_F)$,
define $r_F(x)$ to be the point of $\alpha_x \cap F$ that is closest to 
$x$ in $\alpha_x$.

\vskip 18pt
\centerline{
\includegraphics[width=2.9in]{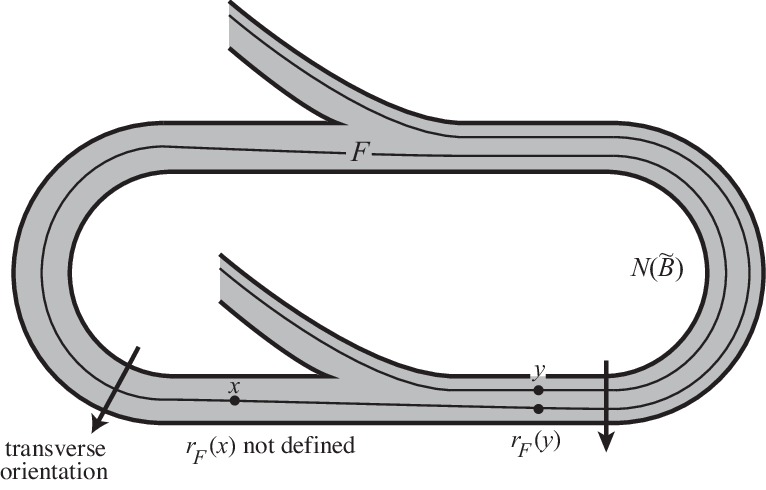}
}
\vskip 6pt
\centerline{Figure 29: First-return map}

\noindent {\bf Lemma 9.10.} {\sl If $F \subseteq F'$, then ${\rm dom}(r_F) \subseteq
{\rm dom}(r_{F'})$.}

\noindent {\sl Proof.} If $x \in {\rm dom}(r_F)$, then $\alpha_x \cap F \not=\emptyset$,
and so $\alpha_x \cap F' \not= \emptyset$. Therefore $x \in {\rm dom}(r_{F'})$.
$\square$

\noindent {\bf Lemma 9.11.} {\sl $r_F$ does not have a fixed point.}

\noindent {\sl Proof.} This is because, for each $x \in {\rm dom}(r_F)$,
$x$ and $r_F(x)$ are
distinct points in the fibre $I_x$. $\square$

\noindent {\bf Lemma 9.12.} {\sl If $F$ is a connected subsurface of $\tilde E$
that is a union of Euclidean tiles, then $F - {\rm dom}(r_F)$ consists of at
most $192n^2$ tiles.}

\noindent {\sl Proof.} For each Euclidean patch of $\tilde B$, the tiles of $F$ in the fibred
neighbourhood of this patch
are parallel. Since $F$ is connected and $\tilde B$ is transversely
orientable, these tiles of $F$ are all coherently oriented.
Hence, on all but one of these tiles of $F$, $r_F$ is defined.
There are at most $192n^2$ Euclidean patches of $\tilde B$, by Property 9.6, which establishes
the lemma. $\square$

\vfill\eject
\noindent {\caps 9.6. The first-return map for large grids}
\vskip 6pt

The following is an easy observation.

\noindent {\bf Lemma 9.13.} {\sl Let $D$ be a grid with side length
at least $14n$. Then $r_D$ is defined at some point of $D$.}

\noindent {\sl Proof.} The branched surface $\tilde B$ contains at most $192 n^2$
patches. Since $D$ has at least $196n^2$ tiles, $D$ must run over
some patch of $\tilde B$ at least twice. Hence, $r_D$ is defined
on one of these patches in $D$. $\square$

However, we need the following rather stronger statement.

\noindent {\bf Proposition 9.14.} {\sl Every point $x \in \tilde E - N_{7000n^2}(\partial \tilde E)$
lies in ${\rm dom}(r_{D(x,4000n^2)})$.}

The key step in the proof of this is the assertion that, for a large grid $D$,
the points in $D$ where $r_{\tilde E}$ fails to be defined lie close to $\partial D$.
In fact, it is convenient to work with several grids simultaneously, as
follows.

\noindent {\bf Proposition 9.15.} {\sl Let $D_1, \dots, D_m \subseteq \tilde E$ be a collection of
disjoint grids, each with side length at least $1500n^2$. For each
$D_i$, let $d_i$ be 
$\sup \{d(y, \partial D_i) : y \in D_i - {\rm dom}(r_{\tilde E}) \}$.
Then $\sum_i d_i \leq 384n^2$.}

\noindent {\sl Proof.} Let $D$ denote the union of the grids
$D_1, \dots, D_m$. We now form a union of annuli and discs $C$ in $N(\tilde B)$,
such that $\tilde E \cap C \subseteq \partial C$, as follows.
Start with the cusps of $N(\tilde B)$. In a regular neighbourhood
of the Euclidean patches of $\tilde B$, the cusps of $\tilde B$
are annuli and discs. This follows from Lemma 8.1, setting $M$ to be
this regular neighbourhood of the Euclidean patches and considering
the branched surface $\tilde B \cap M$. Hence, we may extend each such cusp
vertically into the interior of $N(\tilde B)$ until it just touches $\tilde E$. (See Figure 30.)
Let $C$ be the result. Note that $C \cap \tilde E$ is a collection of simple closed curves
and properly embedded arcs in $\tilde E$. Divide $C \cap \tilde E$ into $\partial_- C$ and $\partial_+ C$, 
where the transverse orientation on $\tilde E$ points into $C$ at  $\partial_- C$,
and out of $C$ at $\partial_+ C$.
Then $\partial_- C \cap D$ forms the intersection between $D \cap {\rm dom}(r_{\tilde E})$
and ${\rm cl}(D - {\rm dom}(r_{\tilde E}))$.  It is a collection of simple closed curves and properly
embedded arcs in $D$.

\noindent {\sl Claim 1.} For each grid $D_i$, any point on $\partial_- C \cap D_i$
that is furthest from $\partial D_i$ lies on an arc component of $\partial_- C \cap D_i$.

\vskip 12pt
\centerline{
\includegraphics[width=2.5in]{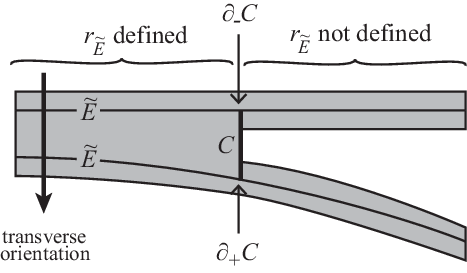}
}
\vskip 6pt
\centerline{Figure 30: The annuli $C$}

Let us assume the claim for the moment.
Now, the total length of the cusps of $\tilde B$ is at most $768n^2$,
by Property 9.7,
and so $\partial_- C \cap D$ also has length at most $768 n^2$. Consider
a point $p$ on $\partial_- C \cap D_i$ which has maximal distance from $\partial D_i$.
By the claim, $p$ can be at a distance at most $384n^2$ from $\partial D_i$. Let $D'_i$ be the
grid with the same centre as $D_i$, but with $p$ on its boundary.
Hence, the interior of $D'_i$ either lies entirely in ${\rm dom}(r_{\tilde E})$
or is entirely disjoint from ${\rm dom}(r_{\tilde E})$. But, by Lemma 9.12,
${\rm dom}(r_{\tilde E})$ is defined for all but at most $192n^2$ tiles.
We are assuming that $D_i$ has side length at least $1500 n^2$.
Hence, the side length of $D'_i$ is at least $1500n^2 - (2 \times 384n^2) = 732n^2$. Therefore there
are at least $(732 n^2)^2$ tiles in $D'_i$ which is more than $192n^2$.
So, we deduce that $D'_i$ lies in ${\rm dom}(r_{\tilde E})$.
Therefore, a point in ${\rm cl}(D_i - {\rm dom}(r_{\tilde E}))$ at maximal
distance from $\partial D_i$ must lie in $\partial_- C \cap D_i$.
So, $d_i$ is at most half the length of $\partial_- C \cap D_i$.
Therefore, $\sum d_i$ is at most half the length of $\partial_- C \cap D$,
which is at most $384n^2$. This proves the proposition.

We still need to prove Claim 1. 

We give $\partial_- C \cap D$ a transverse orientation in $D$, pointing it towards
${\rm dom}(r_{\tilde E})$. Thus, it points `into' $\tilde B$ and away from the cusps.

\noindent {\sl Claim 2.} Each simple closed curve of $\partial_- C \cap D$
points into the disc in $D$ that it bounds.

Claim 1 is a consequence of Claim 2, as follows. Cut $D_i$ along the arc
components of $\partial_- C \cap D_i$, and let $D''_i$ be the disc containing
the centre of $D_i$. Since the arc components have length at most 
$768 n^2$, $D''_i$ contains the grid with the same centre as $D_i$
and with side length $732n^2 - 1$. So, $D''_i$ contains at least $(732n^2 - 1)^2$ tiles.
We will rule out the possibility that there are any simple closed curves
of $\partial_- C$ in $D''_i$. Let $\gamma$ be the union of those components
of $\partial_- C \cap {\rm int}(D''_i)$ that are outermost, in other words, that do not lie
within another component of $\partial_- C \cap {\rm int}(D''_i)$. The total length of $\gamma$
is at most $768n^2$, and so the total number of tiles that it can bound
is at most $(768n^2)^2/ 4 = (384n^2)^2$. But by Claim 2, ${\rm dom}(r_{\tilde E}) \cap D''_i$
lies within $\gamma$. So, at least $(732n^2 - 1)^2 - (384n^2)^2$ tiles do not
lie in ${\rm dom}(r_{\tilde E})$. However, we have already seen in Lemma 9.12 that $r_{\tilde E}$
is defined on all but at most $192n^2$ tiles of $\tilde E$. This is a contradiction,
proving Claim 1.

We now must prove Claim 2. Suppose that there is a simple closed
curve component $\beta$ of $\partial_- C \cap D$ that points out of the
disc that it bounds in $D$. We therefore get a configuration as shown in Figure 31.

\vskip 12pt
\centerline{
\includegraphics[width=3in]{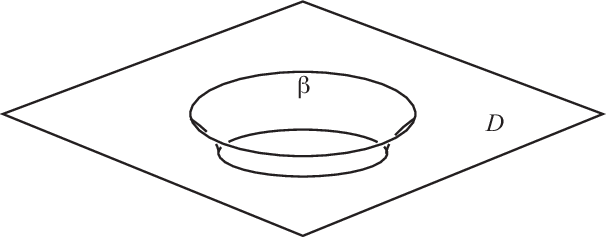}
}
\vskip 6pt
\centerline{Figure 31: An outward-pointing cusp}

Let $\tilde C'$ be the component of $C$ containing $\beta$. 
Let $C'$ be the image of $\tilde C'$ in $N(B')$ under the covering
map $N(\tilde B) \rightarrow N(B')$. This is a cusp of $N(B')$
that has been extended to $E'$. Let $2S$ be two normally parallel copies
of the characteristic surface $S$. Shrink the annulus $C'$ a little so that
its boundary lies in $2S$. Let $C''$ be the resulting annulus.
It is embedded. (Note that $C'$ might not have been embedded since
its two boundary components might have intersected each other.)
The two curves $\partial C''$ bound discs $W_1$ and $W_2$ in $2S$. One
of these discs is parallel to the image in $E'$ of the disc in $D$
bounded by $\beta$. Hence, $W_1$ and $W_2$ are not normally
parallel, because the cusp $C''$ lies between them.

Now, $W_1$ and $W_2$ are disjoint. For if they were nested, say $W_1 \subset W_2$,
then we could remove $W_2$ from $S$ and replace it by $W_1$, thereby
create a normal characteristic surface $S'$ with the same boundary as $S$, but
with $(w_\beta(S'), w(S')) < (w_\beta(S), w(S))$. 

We now form two new normal characteristic surfaces $S_1$ and $S_2$ with the same boundary as
$2S$. The first of these is obtained from $2S$ by removing $W_2$ and replacing
it with a normally parallel copy of $W_1$. Similarly, $S_2$ is obtained from
$2S$ by removing $W_1$ and inserting a normally parallel copy
of $W_2$. Then $S_1$ and $S_2$ are both distinct from $2S$, up to normal
isotopy, because $W_1$ and $W_2$ are not normally parallel.
Note that, as normal surfaces, $4S = S_1 + S_2$. Hence, we deduce that
$S$ is not a boundary-vertex surface. This contradiction proves Claim 2. $\square$

\noindent {\sl Proof of Proposition 9.14.} Let $x$ be a point in $\tilde E - N_{7000n^2}(\partial \tilde E)$.
We will define two increasing sequences of non-negative integers
$m_i$ and $k_i$ and a collection of maps $D(x, 2000n^2 - k_i) \times [0, m_i]
\rightarrow N(\tilde B)$ with the following properties:
\item{(1)} The map is an embedding on 
$D(x, 2000n^2 - k_i) \times [0, m_i)$.
\item{(2)} $D(x, 2000n^2 - k_i) \times \{ 0 \} = D(x, 2000n^2 - k_i) \subset \tilde E$.
\item{(3)} The transverse orientation on 
$D(x, 2000n^2 - k_i)$, which is inherited from that of $\tilde E$,
points into $D(x, 2000n^2 - k_i) \times [0, m_i)$
\item{(4)} For each point $\{ \ast \}$ in $D(x, 2000n^2 - k_i)$, 
$\{ \ast \} \times [0, m_i]$ is a subset of a fibre in $N(\tilde B)$.
\item{(5)} The intersection between $D(x, 2000n^2 - k_i) \times [0, m_i]$ and 
$\tilde E$ is $D(x, 2000n^2 - k_i) \times ([0,m_i] \cap {\Bbb Z})$.

\noindent This sequence will continue until $D(x, 2000n^2 - k_i) \times \{m_i \} \subseteq \tilde E$
has non-empty intersection with $D(x, 2000n^2) \subseteq \tilde E$. 

We start with $D(x, 2000n^2) \times \{ 0 \}$. Set $i$, $m_0$ and $k_0$ to be $0$.
We now apply the following procedure.

\item{(1)} Suppose $D(x, 2000n^2 - k_i) \times [0, m_i]$ has been defined and that
$m_i > 0$.
If $D(x, 2000n^2 - k_i) \times \{m_i\}$ is disjoint from $D(x, 2000n^2)$, then
add it to this product region. Increase $m_i$ by $1$. Pass to step 2.
If $D(x, 2000n^2 - k_i) \times \{m_i\}$ intersects $D(x, 2000n^2)$, then
the procedure terminates.

\item{(2)} Is $r_{\tilde E}(y)$ defined for all $y \in D(x, 2000n^2 - k_i) \times \{ m_i \}$?
If not, then pass to step 3. Otherwise, remain on this step. This 
means that below $D(x, 2000n^2 - k_i) \times \{ m_i \}$, there is
another part of $\tilde E$. Define this to be $D(x, 2000n^2 - k_i) \times \{ m_i +1 \}$.
Between these two surfaces, there is a product region, which we take to be
$D(x, 2000n^2 - k_i) \times (m_i, m_{i}+1)$. Return to step 1.

\item{(3)} In this situation, $r_{\tilde E}(y)$ is not defined for some $y \in D(x, 2000n^2 - k_i ) \times \{ m_i \}$.
This means that there is at least one cusp of $N(\tilde B)$ directly below some
part of $D(x, 2000n^2 - k_i) \times \{ m_i \}$. Let $d_i$ be the maximal distance of such a cusp
from the boundary of $D(x, 2000n^2 - k_i) \times \{ m_i \}$. Applying Proposition 9.15 to the discs
$D(x, 2000n^2 - k_0) \times \{ m_0\}, \dots ,D(x, 2000n^2 - k_i) \times \{ m_i\}$
gives that $\sum_{j=1}^i d_j$ is at most $384n^2$. Set $k_{i+1} = \sum_{j = 1}^i d_j$.
Therefore $D(x, 2000n^2 - k_{i+1})$ is a grid of side length at least 
$2 \times (2000n^2 - 384n^2) \geq 1500n^2$. Let $m_{i+1} = m_i + 1$.
Increase $i$ by 1, and pass to step 1.

When this process terminates, we deduce the existence of points $y \in
D(x, 2000n^2 - k_i)$ and $y' \in D(x, 2000n^2)$ such that 
$y \times \{ m_i \} = y' \times \{ 0 \}$. Thus,
$r_{D(x, 2000n^2)}(y) = y'$.

We now extend $D(x, 2000n^2)$ to the grid $D(x, 4000n^2)$. Now, $x \times \{ m_i \}$
lies within the grid $D(y \times \{ m_i \}, 2000n^2)$. 
Hence, we deduce that $x \times \{ m_i \}$ lies in $D(x, 4000n^2)$. There may be other points of
$D(x, 4000n^2)$ on the fibre between $x \times \{ 0 \}$ and $x \times \{ m_i \}$.
But we deduce that the first-return map for $D(x, 4000 n^2)$ is defined at $x$.
$\square$

\vskip 18pt
\centerline{
\includegraphics[width=5in]{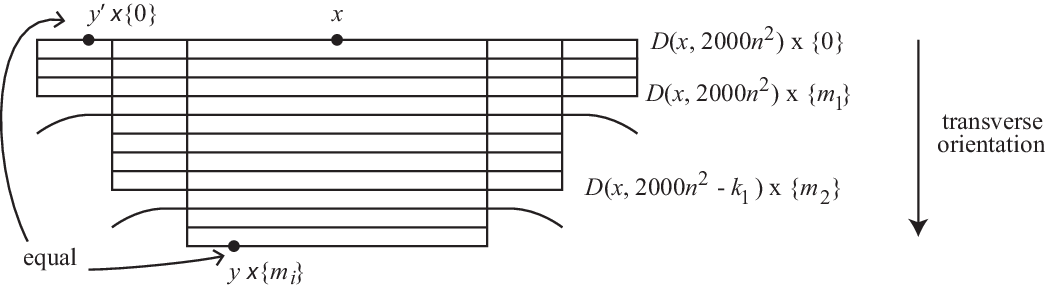}
}
\vskip 6pt
\centerline{Figure 32: Schematic picture of the product regions}

\vskip 18pt
\noindent {\caps 9.7. Translation invariance of first-return maps}
\vskip 6pt

Let $x$ be a point in $\tilde E - N_{7000n^2}(\partial \tilde E)$. 
The points $x$ and $r_{D(x, 4000n^2)}(x)$ both lie in the disc
$D(x, 4000n^2)$ and so there is a well-defined Euclidean translation
vector $v_x$ taking $x$ to $r_{D(x, 4000n^2)}(x)$.
This vector lies in the tangent space $T_x\tilde E$. 

\noindent {\bf Proposition 9.16.} {\sl The vector field $\{ v_x : x \in 
\tilde E - N_{7000 n^2}(\partial \tilde E) \}$ is covariant constant.}

In other words, this vector field on a component of 
$\tilde E - N_{7000 n^2}(\partial \tilde E)$ is the same as the one obtained
by starting with the vector $v_x$ for some fixed $x$ in that
component, and translating using Euclidean parallel translation.

\noindent {\sl Proof.} Clearly the vector field is covariant constant on each tile,
since $\tilde B$ has trivial monodromy.
So suppose that $x$ and $x'$ lie at the centres of adjacent tiles
$t$ and $t'$ of $\tilde E$. Let $\tau$ be the Euclidean translation of length 1 taking
$x$ to $x'$. Then, when passing from $D(x, 4000n^2)$ to $D(x', 4000n^2)$,
the translation $\tau$ is performed. Since $\tilde B$ has trivial monodromy,
the tile containing $r_{D(x, 4000n^2)}(x)$ is also translated by $\tau$.
Hence, it lies in the same patch of $N(\tilde B)$ as $t'$. We claim
that this is the tile containing $r_{D(x', 4000n^2)}(x')$. For otherwise,
there is a tile of $D(x', 4000n^2)$ lying between it and $x'$.
But, then translating this tile by $\tau^{-1}$, we get a tile
of $D(x,4000n^2)$ lying between $x$ and $r_{D(x, 4000n^2)}(x)$,
which is impossible. $\square$

\vskip 18pt
\centerline{
\includegraphics[width=1.5in]{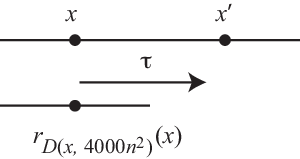}
}
\vskip 6pt
\centerline{Figure 33: Translating $x$ to $x'$}

\noindent {\bf Corollary 9.17.} {\sl Let $x$ be a point in $\tilde E - N_{7000 n^2}(\partial \tilde E)$.
Let $\beta$ be a path starting at $x$ and remaining in $\tilde E - N_{7000 n^2}(\partial \tilde E)$. 
Let $\beta'$ be obtained from $\beta$ by translating each
point in the direction $v_x$. Then $\pi \circ \beta = \pi \circ \beta'$,
where $\pi \colon N(\tilde B) \rightarrow \tilde B$ is the projection map.}

In other words, $\beta$ and $\beta'$ follow the same itinerary through $\tilde B$.
This will be important for us, because curves of this form will form two
sides of a parallelogram which will glue up to form the torus that we
are looking for.

\vskip 6pt
\noindent {\caps 9.8. Completion of the proof}
\vskip 6pt

We are assuming that there is a point $x$ in $\tilde E$ with distance more than
$8000n^2$ from $\partial \tilde E$.

Let $\alpha$ be the geodesic in $D(x, 4000n^2)$ from $x$ to
$x' = r_{D(x, 4000n^2)}(x)$. Let $\beta$ be a geodesic going through $x$
orthogonal to $\alpha$ with length $1000n^2$ in both directions from
$x$. Then $\beta$ remains in $\tilde E - N_{7000n^2}(\partial \tilde E)$.
Let $\beta'$ be the result of translating $\beta$
using the vector $v_x$, so that it runs through $x'$. Then we refer
to the region between $\beta$ and $\beta'$ as a {\sl strip}, and
we denote it by $P$. It is a Euclidean rectangle.

In Section 9.4, the proof divided into two cases: when there is a closed
geodesic in $\tilde E - N_{1000 n^2}(\partial \tilde E)$ with length at most 
$12000n^2$, and when there is not. We initially focused on the case
where there is no such geodesic, and have defined $x$, $\alpha$, $\beta$ and $P$
in this case. But now we want to reintegrate the
two parts of the argument. So, suppose that there is such a geodesic,
which may take to be simple, and call it
$\alpha$. Let $x$ be a point on $\alpha$, and let $\beta$ be a geodesic
through $x$ that is orthogonal to $\alpha$. Suppose that it has length
$1000 n^2$ in both directions from $x$. We proved in Lemma 9.9 that
$N_{1000n^2}(\alpha)$ is isometric to a Euclidean annulus with core curve $\alpha$.
Hence, $\beta$ cuts $N_{1000n^2}(\alpha)$ into a Euclidean rectangle.
We also call this a {\sl strip}, and denote it by $P$.

We now want to emulate the proof of Proposition 9.14, but instead of
starting with a grid, we will start with this strip.

Let $p \colon P \rightarrow \beta$ be orthogonal projection. If $V$ is a finite
union of closed intervals in $\beta$, we say that $p^{-1}(V)$ is {\sl strip-like}.

We will define an increasing sequence of non-negative integers
$m_i$ and a collection of strip-like subsets $P = P_0 \supseteq P_1 \supseteq \dots
\supseteq P_k$ of $P$, with the
following properties.
\item{(1)} There is map $P_i \times [0, m_i] \rightarrow N(\tilde B)$
which is an embedding on $P_i \times [0, m_i)$.
\item{(2)} $P_i \times \{ 0 \} = P_i \subseteq P$.
\item{(3)} The transverse orientation on $P_i$ points into
$P_i \times [0,m_i)$.
\item{(4)} For each point $\{ \ast \}$ in $P_i$, 
$\{ \ast \} \times [0, m_i]$ is a subset of a fibre in $N(\tilde B)$.
\item{(5)} The intersection between $P_i \times [0, m_i]$ and 
$\tilde E$ is $P_i \times ([0,m_i] \cap {\Bbb Z})$.

\noindent This sequence will continue until $P_i \times \{m_i \} \subseteq \tilde E$
has non-empty intersection with $P \subseteq \tilde E$. 

We start with $P_0 = P \times \{ 0 \}$. Set $i$ and $m_0$ to be $0$.
We now apply the following procedure.

\item{(1)} Suppose $m_i > 0$, that $P_i \times [0, m_i] \rightarrow N(\tilde B)$ has been defined
and that it is an embedding on $P_i \times [0, m_i)$.
If $P_i \times \{m_i\}$ is disjoint from $P$, then
add it to this product region. Increase $m_i$ by $1$. Pass to step 2.
If $P_i \times \{m_i\}$ intersects $P$, then terminate this procedure.

\item{(2)} Is $r_{\tilde E}(y)$ defined for all $y \in P_i \times \{ m_i \}$?
If not, then pass to step 3. Otherwise, remain on this step. This 
means that below $P_i \times \{ m_i \}$, there is
another part of $P$. Define this to be $P_i \times \{ m_i +1 \}$.
Between these two surfaces, there is a product region, which we take to be
$P_i \times (m_i, m_i +1)$. Return to step 1.

\item{(3)} In this situation, $r_{\tilde E}(y)$ is not defined for some $y \in P_i \times \{ m_i \}$.
This means that there is at least one cusp of $N(\tilde B)$ directly below some
part of $P_i \times \{ m_i \}$. Extend these cusps vertically into $N(\tilde B)$
until they just touch $P_i \times \{ m_i \}$. Let $C_i$ be the intersection of these extended cusps
with $P_i \times \{ m_i \}$, and let $N(C_i)$ be a thin regular neighbourhood of $C_i$.
Define $P_{i+1}$ to be $P_i - {\rm int}(p^{-1}p(N(C_i))$. This is strip-like.
Let $m_{i+1} = m_i + 1$. Increase $i$ by 1, and pass to step 1.

\vskip 18pt
\centerline{
\includegraphics[width=4.5in]{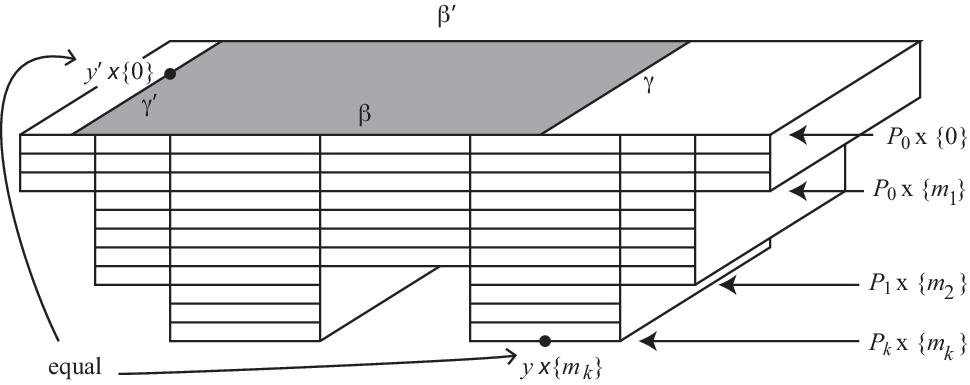}
}
\vskip 6pt
\centerline{Figure 34: Schematic picture of the product regions and strip-like regions}

Now, the total length of $C_1 \cup \dots \cup C_k$ is at most $768n^2$. Therefore,
the length of $p( C_1 \cup \dots \cup C_k)$ is also at most $768n^2$. We therefore 
deduce that when this process terminates $P_{k}$ is non-empty.
It terminates because, for some $y \in P_{k}$, $\{ y \} \times \{ m_k \}$ equals some $y' \in P \times \{ 0 \}$.
Let $\gamma$ be a geodesic starting at $y \times \{ 0 \}$ in the direction of $\alpha$,
and define $\gamma'$ similarly starting at $y' \times \{ 0 \}$. 
Then $\gamma \times [0, m_k] \subset P_k \times [0,m_k]$
forms a product region between $\gamma$ and $\gamma'$, where each fibre in
this product region lies in a fibre in $N(\tilde B)$. Hence, $\gamma$ and $\gamma'$
follow the same itinerary in $\tilde B$.

Let $\delta$ be the subset of $\beta$ lying between $\beta \cap \gamma$ and
$\beta \cap \gamma'$. Define $\delta' \subset \beta'$ similarly. Then
$\gamma \cup \delta \cup \gamma' \cup \delta'$ forms the boundary of
rectangle in $\tilde E$. Opposite sides of this rectangle have the same
image in $\tilde B$. Hence, if we identify opposite sides of this rectangle, the
result is a torus that is carried by $\tilde B$. It is a summand of $\tilde E$
by Lemma 8.2.

This proves Proposition 9.8, which completes the proof of Theorem 7.3 and
hence the main results of this paper. $\square$

\vskip 18pt
\centerline {\caps 10. Final remarks}
\vskip 6pt

\noindent {\caps 10.1. Improving the degree of the polynomials}
\vskip 6pt

We now know that there is a polynomial upper bound on the number
of Reidemeister moves required to turn a diagram of the unknot
or split link into a trivial or disconnected diagram. It is natural to
try to determine the smallest possible degree of such a polynomial.
The result of Hass and Nowik [15] implies that one cannot do better
than a quadratic polynomial. However, the degrees of the polynomials
in Theorems 1.1 and 1.2 are 11.

This can certainly be reduced from 11 to 10,
as follows. In the proof of Theorem 1.4, we started with a reducing
2-sphere with binding weight at most $n 2^{7n^2}$, where $n$ is
the arc index. However, if one starts with a diagram of the link 
having $c$ crossings, then one can find a triangulation of its exterior
using at most $8c$ tetrahedra, and hence find a reducing
sphere with weight at most $c2^{56c}$ with respect to this triangulation. 
One can then compare this
triangulation with that of Dynnikov, and hence find a reducing sphere
with binding weight that is an exponential function of $c$ rather than $c^2$.
If one follows the remainder of the argument of Theorem 1.4, one
finds that one has reduced the degree of the polynomial in Theorem 1.2 by 1 down to 10.
One can do the same for the polynomial in Theorem 1.1.
We have chosen not to pursue this argument here, because it is
somewhat lengthy.

It seems very hard to reduce the degree below 10 using these arguments.

\vskip 6pt
\noindent {\caps 10.2. Further problems}
\vskip 6pt

This paper raises many interesting and difficult questions. We mention some these.

Is there a polynomial time algorithm to recognise the unknot? It is the
author's best guess that there is not, but a proof of such a fact would
be extremely hard.

Can the arguments in this paper be applied to other knot types?
In particular, can one find an upper bound on the number of 
Reidemeister moves required to transform one diagram of a knot into another
that is a polynomial function of the number of crossings in each diagram?
Currently, the only known upper bound on Reidemeister moves for
arbitrary knots, which is due to Coward and the author [5], is much larger than this. It is of the form of 
a tower of exponentials.

\vskip 18pt
\centerline{\caps References}
\vskip 6pt

\item{1.} {\caps I. Agol,} {\sl Knot genus is NP,} Conference presentation (2002)

\item{2.} {\caps D. Bennequin,} {\sl Entrelacements et \'equations de Pfaff,} 
Ast\'erisque 107--108 (1983) 87--161. 

\item{3.} {\caps J. Birman, E. Finkelstein,} {\sl Studying surfaces via closed braids.}
J. Knot Theory Ramifications 7 (1998) 267--334. 

\item{4.} {\caps J. Birman, W. Menasco,} {\sl Studying links via closed braids. IV.
Composite links and split links.} Invent. Math. 102 (1990) 115--139. 

\item{5.} {\caps A. Coward, M. Lackenby}, {\sl An upper bound on Reidemeister moves},
Amer. J. Math (accepted), arXiv:1104.1882.

\item{6.} {\caps P. Cromwell}
{\sl Embedding knots and links in an open book. I. Basic properties.}
Topology Appl. 64 (1995) 37--58. 

\item{7.} {\caps P. Cromwell, I. Nutt,} {\sl Embedding knots and links in an open book. 
II. Bounds on arc index.} Math. Proc. Cambridge Philos. Soc. 119 (1996) 309--319. 

\item{8.} {\caps I. Dynnikov,} {\sl Arc-presentations of links: monotonic simplification.} 
Fund. Math. 190 (2006), 29--76.

\item{9.} {\caps W. Floyd, U. Oertel,}
{\sl Incompressible surfaces via branched surfaces.} Topology 23 (1984) 117--125. 

\vfill\eject
\item{10.} {\caps O. Goldreich}, {\sl Computational complexity. A conceptual perspective},
Cambridge University Press (2008).

\item{11.} {\caps W. Haken}, {\sl Theorie der Normalfl\"achen.} 
Acta Math. 105 (1961) 245--375. 

\item{12.} {\caps W. Haken}, 
{\sl \"Uber das Hom\"oomorphieproblem der 3-Mannigfaltigkeiten. I.} 
Math. Z. 80 (1962) 89--120. 

\item{13.} {\caps J. Hass, J. Lagarias,} {\sl The number of Reidemeister moves needed for unknotting.}
J. Amer. Math. Soc. 14 (2001), no. 2, 399--428

\item{14.} {\caps J. Hass, J. Lagarias, N. Pippenger},
{\sl The computational complexity of knot and link problems.} J. ACM 46 (1999) 185--211.

\item{15.} {\caps J. Hass, T. Nowik,} {\sl Unknot diagrams requiring a quadratic number of 
Reidemeister moves to untangle.} Discrete Comput. Geom. 44 (2010), no. 1, 91--95.

\item{16.} {\caps J. Hass, J. Snoeyink, W. Thurston}, 
{\sl The size of spanning disks for polygonal curves.} Discrete Comput. Geom. 29 (2003) 1--17. 

\item{17.} {\caps G. Hemion}, {\sl On the classification of homeomorphisms of $2$-manifolds and the classification of $3$-manifolds.}
Acta Math. 142 (1979), no. 1-2, 123--155.

\item{18.} {\caps A. Henrich, L. Kauffman}, {\sl Unknotting Unknots.} arXiv:1006.4176 

\item{19.} {\caps W. Jaco, U. Oertel},
{\sl An algorithm to decide if a 3-manifold is a Haken manifold.}
Topology 23 (1984) 195--209.

\item{20.} {\caps W. Jaco, J. Tollefson}, {\sl Algorithms for the complete decomposition
of a closed 3-manifold}, Illinois J. Math. 39 (1995) 358--406.

\item{21.} {\caps G. Kuperberg}, {\sl Knottedness is in NP, modulo GRH.} Adv. Math. (accepted), arXiv:1112.0845 

\item{22.} {\caps S. Matveev}, {\sl Algorithmic topology and classification of 3-manifolds.}
Algorithms and Computation in Mathematics, 9. Springer, Berlin, 2007.

\item{23.} {\caps J. Storer}, {\sl On Minimal-Node-Cost Planar Embeddings}, Networks 14 (1984) 181--212.

\item{24.} {\caps A. Turing,} {\sl Solvable and Unsolvable Problems}, Science News 31 (1954).

\vskip 18pt
\+ Mathematical Institute, University of Oxford, \cr
\+ Radcliffe Observatory Quarter, Woodstock Road, \cr
\+ Oxford OX2 6GG, United Kingdom. \cr

\end